\documentclass[oneside,english]{amsart}
\usepackage[T1]{fontenc}
\usepackage[latin9]{inputenc}
\usepackage{textcomp}
\usepackage{amsthm}
\usepackage{amstext}
\usepackage{esint}

\makeatletter
\numberwithin{equation}{section}
\numberwithin{figure}{section}
\theoremstyle{plain}
\newtheorem{thm}{\protect\theoremname}[section]
  \theoremstyle{plain}
  \newtheorem{cor}[thm]{\protect\corollaryname}
  \theoremstyle{definition}
  \newtheorem{defn}[thm]{\protect\definitionname}
  \theoremstyle{plain}
  \newtheorem{prop}[thm]{\protect\propositionname}
  \theoremstyle{plain}
  \newtheorem{lem}[thm]{\protect\lemmaname}
  \theoremstyle{remark}
  \newtheorem{rem}[thm]{\protect\remarkname}

\def\makebbb#1{
    \expandafter\gdef\csname#1\endcsname{
        \ensuremath{\Bbb{#1}}}
}\makebbb{R}\makebbb{N}\makebbb{Z}\makebbb{C}\makebbb{H}\makebbb{E}\makebbb{H}\makebbb{P}\makebbb{B}\makebbb{Q}\makebbb{E}

\usepackage{babel}

\usepackage{babel}

\makeatother

\usepackage{babel}

\makeatother

\usepackage{babel}
  \providecommand{\corollaryname}{Corollary}
  \providecommand{\definitionname}{Definition}
  \providecommand{\lemmaname}{Lemma}
  \providecommand{\propositionname}{Proposition}
  \providecommand{\remarkname}{Remark}
\providecommand{\theoremname}{Theorem}

\begin{document}

\title{statistical mechanics of permanents, real-Monge-Ampère equations
and optimal transport }

\author{Robert J. Berman}
\begin{abstract}
We give a new probabilistic construction of solutions to real Monge-Ampère
equations in $\R^{n}$ (satisfying the second boundary value problem
with respect to a given target convex body $P)$ which fits naturally
into the theory of optimal transport. More precisely, certain $\beta$-deformed
permanental (bosonic) $N-$particle point processes are introduced
and their empirical measures are, in the large $N-$limit, shown to
converge exponentially in probability to a deterministic measure whose
potential satisfies the real Monge-Ampère equation in question. In
particular, this allows us to represent the solution as a limit, as
$N\rightarrow\infty,$ of explicit integrals over the $N-$fold products
of $\R^{n}$ and it also leads to explicit limit formulas for optimal
transport maps to the convex body $P.$ Connections to the study of
Kähler-Einstein metrics on complex algebraic varieties, and in particular
toric ones, are briefly discussed.
\end{abstract}

\email{robertb@chalmers.se}

\curraddr{Mathematical Sciences - Chalmers University of Technology and University
of Gothenburg - SE-412 96 Gothenburg, Sweden}
\begin{abstract}
\tableofcontents{}
\end{abstract}
\maketitle

\section{Introduction}

In this paper we will develop a probabilistic/statistical mechanical
approach for producing solutions to real \emph{Monge-Ampère equations}
in $\R^{n},$ satisfying the second boundary value problem with respect
to a given target convex body $P.$ This fits naturally into the theory
of\emph{ optimal transport} \cite{vi1,vi2}. In particular, it will
lead to a probabilistic construction of optimal transport plans from
a set $X$ in $\R^{n}$ to a target convex body $P.$ The approach
arose as a ``spin-off effect'' of the authors work on a probabilistic
approach to Kähler-Einstein metrics on complex algebraic varieties
and, more generally, complex Monge-Ampère equations on complex manifolds
(see \cite{be-3} for an out-line of the general complex geometric
setting and \cite{be} for connections to emergent gravity and boson-fermion
correspondences in physics). For example, the\emph{ permanental} random
point processes on $\R^{n}$ introduced below, which are determined
by the target convex body $P,$ are the push-forwards to $\R^{n}$
of \emph{determinantal }point-processes defined on the complex torus
$\C^{*n}$ and the corresponding limiting real Monge-Ampère measures
on $\R^{n}$ are the push-forwards of the corresponding limiting complex
Monge-Ampère measures on $\C^{*n}.$ The push-forward map is the one
induced from the standard identification of $\C^{*n}$ with $T^{n}\times\R^{n}$,
where $T^{n}$ is the real unit-torus and in the case when the convex
body $P$ is a rational polytope the corresponding determinantal point
processes are naturally defined on the corresponding toric projective
algebraic variety $X_{P},$ compactifying $\C^{*n}.$ From this point
of view the optimal transport theory in $\R^{n}$ thus arises as the
``push-forward'' of the pluripotential theory appearing in the complex
setting \cite{ber-bou,b-b-w,bbgz}.

The general complex geometric framework will be considered in detail
elsewhere \cite{be-4}. Accordingly, we will in this paper concentrate
on the corresponding real setting (see however sections \ref{sub:Relation-to-the}
and \ref{sub:K=0000E4hler-Einstein-metrics,-negativ} for some relations
to the complex setting).

\subsection{The Monge-Ampère, optimal transport and permanental point processes }

In their simplest classical form the real Monge-Ampère equations that
are the focus of the present paper are of the form 
\begin{equation}
\det(\frac{\partial^{2}\phi}{\partial x_{i}\partial x_{j}})=e^{\beta\phi}\rho_{0}(x)\label{eq:ma eq intro cassical}
\end{equation}
for a given function $\rho_{0}(x)$ of unit-mass and a given parameter
$\beta\geq0.$ As usual the solution $\phi$ is demanded to be convex,
but to ensure uniquess further growth conditions at infinity have
to be specified. The relevant situation here will be when $\phi$
grows as the support function $\phi_{P}$ of a given $n-$dimensional
convex body $P$ in $\R^{n}.$ In the PDE literature \cite{ba,vi2}
this is sometimes called the \emph{second boundary value problem }for
the equation above and (for $\rho_{0}$ strictly positive) it turns
out to be equivalent to demanding that the gradient $\nabla\phi$
map $\R^{n}$ diffeomorphically onto the interior of $P:$ 
\begin{equation}
\nabla\phi:\R^{n}\rightarrow P\label{eq:gradient image intro}
\end{equation}
Accordingly, $P$ is sometimes referred to as the\emph{ target.} More
generally, we will consider the setting of a given triple $(P,\mu_{0},\phi_{0})$
consisting of a convex body $P$ in $\R^{n}$ a (Borel) measure $\mu_{0},$
whose support will be denoted by $X$ and a \emph{weight function}
$\phi_{0}$ on $\R^{n},$ i.e. a (a possible non-convex) continuous
function $\phi_{0}$ which grows faster than the support function
$\phi_{P}$ of $P$ at infinity in $\R^{n}$ (see section \ref{sub:Weigted-sets-and})
and such that $e^{-\beta\phi_{0}}\mu_{0}$ has finite total mass.
The corresponding Monge-Ampère equation is then 
\begin{equation}
MA(\phi)=e^{\beta(\phi-\phi_{0})}\mu_{0},\label{eq:ma eq weak intro}
\end{equation}
 where $MA(\phi)$ is the Monge-Ampère measure of $\phi,$ in the
sense of Alexandrov (in particular, for $\phi$ smooth its density
is the determinant of the Hessian appearing in the previous equation)
and the condition \ref{eq:gradient image intro} is assumed to hold
in the sense of sub-gradients \cite{gu}. We may also, after a trivial
scaling, assume that $P$ has unit Euclidean volume. 

In the case when $\beta=0$ (and, say, $\phi_{0}=0)$ the corresponding
Monge-Ampère equation, i.e. the equation 
\[
MA(\phi)=\mu_{0},
\]
where now $\mu_{0}$ is assumed to be a probability measure, plays
a central role in the theory of \emph{optimal transport }\cite{vi1,vi2}.
Under appropriate regularity assumptions on the measure $\mu_{0}$
a solution $\phi$ defines a map $T:=\nabla\phi$ from $\R^{n}$ to
$P,$ which coincides with the so called\emph{ optimal transport map,
}defined with respect to the the target measure $\lambda_{P}:=1_{P}dp$
(i.e. the normalized Lesbesgue measure on the convex body $P)$ and
the\emph{ cost function} $c(x,p)=-x\cdot p.$ This means that it\emph{
}minimizes the corresponding total transport \emph{cost functional}
$C(T)$ over all maps $T$ transporting (i.e. pushing forward) $\mu_{0}$
to $\lambda_{P}.$ The precise definitions are recalled in the appendix.
For the moment we just recall that the cost functional $C$ is more
generally defined on the space of all \emph{couplings} $\Gamma$ (also
called \emph{transference plans}) between two given probability measures
$\mu$ and $\nu,$ i.e $\Gamma$ is a measure on $\R^{n}\times\R^{n}$
whose push-forwards to the first and second factor are equal to $\mu$
and $\nu,$ respectively and 
\[
C(\Gamma):=\int_{\R^{n}\times\R^{n}}c(x,p)\Gamma
\]
Fixing $\mbox{\ensuremath{\nu}}=\lambda_{P}$ we will also write $C(\mu)$
for the corresponding \emph{optimal cost funtional} on the space of
all probability measures on $\R^{n},$ i.e. $C(\mu)$ is the minimal
cost to transport the measure $\mu$ to \emph{$\lambda_{P}:$} 
\[
C(\mu)=\inf_{\Gamma}C(\Gamma)
\]
where $\Gamma$ ranges over all couplings between $\mu$ and $\lambda_{P}.$
It will also be important to consider the ``weighted'' cost functional
$C_{\phi_{0}}(\mu)$ defined in terms of the cost function 
\begin{equation}
c_{\phi_{0}}(x,p):=-x\cdot p+\phi_{0}(x)\label{eq:def of cost function intro}
\end{equation}
We will, in particular, show how to recover the solution $\phi$ of
the equation \ref{eq:ma eq weak intro} from the large $N-$limit
of a certain random point process on $\R^{n}$ with $N$ particles,
canonically determined by the given data $(\mu_{0},\phi_{0},P).$
In particular, this will, in the case $\beta=0,$ yield explicit explicit
approximations for optimal transport maps. From the point of view
of equilibrium statistical mechanics the parameter $\beta$ will play
the role of the \emph{inverse temperature }and the approach will involve
the limiting zero temperature case (i.e. $\beta=\infty)$ where, as
explained below, the role of the equation \ref{eq:ma eq weak intro}
is played by a \emph{free boundary value problem} for the Monge-Ampère
equation, which can be equivalently described as a constrained\emph{
convex envelope. }As will be made clear below the results split into
three different ``phases'' of increasing temperature: $\beta=\infty,$
$\beta>0$ and $\beta=0$ (and we will also briefly comment on the
negative temperature case $\beta<0$ in section \ref{sub:K=0000E4hler-Einstein-metrics,-negativ}
).

In the case $\beta=\infty$ we start with a weight function $\phi_{0}$
defined on a closed subset $X$ of $\R^{n}$ and the corresponding
convex envelope $\phi_{e}$ is then defined as a point wise sup of
convex functions $\phi:$ 
\begin{equation}
\phi_{e}(x):=\sup\{\phi(x):\,\,\,\,\phi\leq\phi_{0}\,\,\mbox{on\,\ensuremath{X,\,\,\,\,\nabla\phi\in P\},}}\label{eq:def of convex env intro}
\end{equation}
The Monge-Ampère measure of $\phi_{e}$ which is supported on $X,$
will be denote by $\mu_{e}.$ In fact, even if $X$ is non-compact
it follows from the growth assumption on the weight $\phi_{0}$ that
the support of $\mu_{e}$ is compact. When $\beta\rightarrow\infty$
the solutions of the corresponding Monge-Ampère equations \ref{eq:ma eq weak intro}
indeed converge, to the envelope $\phi_{e},$ where $X$ is the support
of $\mu_{0},$ at least under an appropriate regularity assumption
 (see Prop \ref{prop:conv of solutions as beta tends to infy}). Interestingly,
the Monge-Ampère measure $\mu_{e}$ also admits a natural interpretation
in terms of the theory of optimal transport. Indeed, the measure $\mu_{e}$
minimizes, among all probability measures $\mu$ on $X,$ the optimal
transport cost $C_{\phi_{0}}(\mu)$ defined by the cost function in
formula \ref{eq:def of cost function intro}. Under suitable regularity
assumptions (for example when $X=\R^{n}$ and is $\phi_{0}$ smooth)
the map $T:=\nabla\phi_{e}$ is the corresponding optimal \emph{transport
map} from the support of $\mu_{e}$ to the convex body $P$ (compare
section \ref{sub:Optimal-transport-theory}).

Next, we turn to the definition of the corresponding random point
processes, which will be defined as ``$\beta-$deformations'' of
certain \emph{permanental} random point processes, interpolating between
a Poisson process for $\beta=0$ and a permanental point process at
$\beta=\infty$. First recall that the\emph{ permanent }of a rank
$N$ matrix $A:=(A_{ij})$ is the real number defined by 
\[
\mbox{per}(A):=\sum_{i=1}^{N}\prod_{\sigma\in S_{N}}a_{i,\sigma(i)},
\]
i.e. it is obtained from the definition of the determinant by removing
the sign dependence on the permutation $\sigma.$ Denote by $P_{\Z}$
the intersection of the convex body $P$ with the integer lattice
$\Z^{n}$ and fix an auxiliary ordering $p_{1},...,p_{N}$ of the
elements of $P_{\Z}.$ Then the cost function \ref{eq:def of cost function intro}
determines a function $A(x_{1},....,x_{N})$ on $(\R^{n})^{N}$ with
values in the space of $N$ times $N-$matrices, defined by 
\[
A_{ij}(x_{1},...,x_{N}):=e^{-c(x_{i},p_{j})}
\]
and we will denote by $\mbox{Per}(x_{1},...,x_{N})$ its permanent,
which defines a real-valued function on $(\R^{n})^{N},$ which is
independent of the ordering of the lattice points $p_{1},...,p_{N}.$
Indeed, we may write
\[
\mbox{Per}(x_{1},...,x_{N})=\sum_{\{p_{1},....,p_{N}\}}e^{x_{1}\cdot p_{1}+\cdots x_{2}\cdot p_{N}},
\]
where, for $(x_{1},...,x_{N})$ fixed, the outer sum ranges over all
possible $N!$ choices of $N$ different elements $p_{1},...,p_{N}$
in $P_{\Z}.$ Given a weighted measure $(\mu_{0},\phi_{0})$ as above
we then obtain a symmetric probability measure on $(\R^{n})^{N},$
i.e. a \emph{random point process }on $X$ with $N$ particles, by
letting its density $\rho^{(N)}(x_{1},...,x_{N})$, with respect to
the product measure $\mu_{0}^{\otimes N},$ be proportional to the
weighted permanent above, i.e.

\[
\mu^{(N)}:=\frac{1}{Z_{N}}\left(\mbox{Per}(x_{1},...,x_{N})e^{-(\phi_{0}(x_{1})+\cdots+\phi(x_{N}))}\right)\mu_{0}^{\otimes N}
\]
 where $Z_{N}$ is the normalizing constant that will sometimes write
as $Z_{N,}[\phi_{0}]$ to indicate the dependence on $\phi_{0.}.$
This is an example of a\emph{ permanental random point process} on
$X$ (see the survey \cite{h-k-p} for general properties of permanental
point processes). Here we just recall that in the physics literature
on many particle quantum systems such processes are used to represent
a gas of free \emph{bosons} \cite{n-o} and accordingly the name \emph{boson
point processes} is sometimes also used in the litterature. In the
present setting the corresponding bosons consist of $N$ ``plane
waves'' $e^{ix\cdot k_{i}}$ with\emph{ imaginary }momenta (wave
numbers) \emph{$k_{i}$} or more precisely: $k_{i}=-ip_{i}$ where
$p_{i}$ ranges over the lattice points $P_{\Z}$ of $P.$

The \emph{empirical measure} of the random point process introduced
above, will be denoted by $\delta_{N}.$ This is the random measure
defined by 
\begin{equation}
\delta_{N}:=\frac{1}{N}\sum_{i=1}^{N}\delta_{x_{i}}\label{eq:empirical measure intro}
\end{equation}
We will consider the large $N$ limit which appears when $P$ is replaced
by the sequence $kP$ of scaled convex bodies, for any positive integer
$k$ (so that $N\sim k^{n}$ to the leading order) and the weight
$\phi$ is replace by $k\phi$ and we will be concerned with the $\beta-$deformations
of the permanental point process above, defined by the probability
measure 
\[
\mu_{\beta_{N}}^{(N)}:=\frac{1}{Z_{N,\beta_{N}}}\left(\mbox{Per}(x_{1},...,x_{N})e^{-k(\phi_{0}(x_{1})+\cdots+\phi(x_{N}))}\right)^{\beta_{N}/k}\mu_{0}^{\otimes N}
\]
where the sequence $\beta_{N}$ is assumed to satisfy 

\[
\lim_{N\rightarrow\infty}\beta_{N}=\beta\in]0,\infty]
\]
(strictly speaking we should really write $N=N_{k}$ to indicate the
dependence on $k,$ but we have omitted the subscript $k$ to simplify
the notation). In particlar, when $\beta_{N}=k$ (so that $\beta=\infty)$
we get a sequence of bona fide permanental point process. Note also
that in the case when $\beta_{N}$ is exactly equal to $\beta,$ for
$\beta$ finite, we could as well suppose that the weight function
$\phi_{0}$ vanishes identically, by replacing $\mu_{0}$ with $e^{-\beta\phi_{0}}\mu_{0},$
but, in fact, it will be useful to separate the weight $\phi_{0}$
from the measure $\mu_{0}.$ 
\begin{thm}
\label{thm:main for perman intro}Given data $(\mu_{0},\phi_{0},P,\beta)$
as above the empirical measure $\delta_{N}$ of the corresponding
random point process on the set $X$ (defined as the support of $\mu_{0})$
converges in probability to the deterministic measure $\mu_{\beta}$
defined by $\mu_{\beta}=MA(\phi_{\beta}),$ where $\phi_{\beta}$
denotes the unique convex solution of equation \ref{eq:ma eq weak intro}
(satisfying \ref{eq:gradient image intro}) for $\beta<\infty$ and
for $\beta=\infty$ the function $\phi_{\beta}$ is the convex envelope
\ref{eq:def of convex env intro}. More precisely, the law of the
empirical measure $\delta_{N}$ (i.e. the probability measure $(\delta_{_{N}})_{*}(\mu_{\beta_{N}}^{(N)}))$
admits a large deviation principle (LDP) with rate $N\beta_{N}$ and
rate functional $F_{\beta},$ where 
\[
F_{\beta}(\mu)=C_{\phi_{0}}(\mu)+\frac{1}{\beta}D_{\mu_{0}}(\mu)-C_{\beta},
\]
 where $C_{\phi_{0}}(\mu)$ is the Monge-Kantorovich optimal cost
functional corresponding to the cost function in formula \ref{eq:def of cost function intro},
$D_{\mu_{0}}(\mu)$ is the entropy of $\mu$ relative to the background
measure $\mu_{0}$ and $C_{\beta}$ is the constant ensuring that
the infimum of $F_{\beta}$ is equal to $0.$ 
\end{thm}
The convergence in probability appearing the previous theorem is equivalent
to the fact that the law of the empirical measure converges weakly
to $\delta_{\mu_{\beta}}$, the Dirac measure at $\mu_{\beta}.$ In
turn, the LDP implies, since $\mu_{\beta}$ is the unique minimizer
of the rate functional $F_{\beta},$ that the latter converge is \emph{exponential
}in a sense which may be loosely formulated as follows: denote by
$\mathcal{B}_{\delta}(\mu)$ a ball of radius $\delta,$ centered
at $\mu$ in the space $\mathcal{M}_{1}(X)$ of all probability measure
on $X,$ equipped with a metric defining the weak topology, then 
\[
\mbox{\ensuremath{\mbox{Prob}}(}\frac{1}{N}\sum_{i=1}^{N}\delta_{x_{i}}\in\mathcal{B}_{\delta}(\mu))\sim e^{-\beta_{N}NF_{\beta}(\mu)}
\]
 as $N\rightarrow\infty$ and $\delta\rightarrow0$ (see section \ref{sub:Large-deviations-(proof}
for the precise definition of the LDP). The Poisson case, i.e. when
$\beta=0,$ is the content of Sanov's classical theorem, which in
turn is a generalization of Cramer's theorem for random vectors in
$\R^{n}$\cite{de-ze}. We also point of that the proof of Theorem
\ref{thm:main for perman intro} will give that $C_{\beta}$ is equal
to the following constant only depending on $(X,\phi_{0},\beta)$ 

\[
C_{\beta}:=C(X,\phi_{0},\beta):=\lim_{N\rightarrow\infty}-\frac{1}{N\beta_{N}}\log Z_{N,\beta_{N}}[\phi_{0}]=\inf_{\mu}\left(C_{\phi_{0}}(\mu)+\frac{1}{\beta}D_{\mu_{0}}(\mu)\right),
\]
 where, in the case $\beta=\infty$ the infimum above is taken over
all probability measures $\mu$ supported on the support $X$ of $\mu_{0}.$ 

Before commenting on the proof of the previous theorem we state some
of its corollaries. First, from the convergence of the one-point correlation
measures we obtain the following corollary, which provides a sequence
of\emph{ explicit} approximate solutions to the real Monge-Ampère
equations above.
\begin{cor}
\label{cor:first cor intro}Let $\mu_{0}$ be a weighted measure and
$P$ an $n-$dimensional convex body and fix $\beta>0.$ If $\mu_{0}=\rho_{0}dx$
for a strictly positive function $\rho_{0},$ then 

\[
\phi_{\beta}^{(N)}(x):=\frac{1}{\beta}\log\int_{X^{N-1}}\frac{1}{Z_{N}}(\mbox{Per}(x,x_{2},...,x_{N}))^{\beta/k}(e^{-\beta\phi_{0}}\mu_{0})^{\otimes N-1}
\]
 satisfies \ref{eq:gradient image intro} and converges locally uniformly,
as $N\rightarrow\infty,$ to the unique solution of the second boundary
value problem with target $P$ for the Monge-Ampère equation \ref{eq:ma eq weak intro}. 
\end{cor}
In fact, in a similar way also obtain explicit explicit approximations
to the inhomogeneous Monge-Ampère equation (obtained by setting $\beta=0$
in \ref{eq:ma eq intro cassical}). Formally, this a consequence of
the previous corollary in the limiting case $\beta=0,$ but the proof
proceeds in a somewhat different manner. 
\begin{cor}
\label{cor:sec cor intro}Let $\mu_{0}$ be a probability measure
of the form $\mu_{0}=\rho_{0}1_{X}dx$ such that $X$ is the closure
of a bounded domain whose boundary $\partial X$ is a null set for
Lebesgue measure and assume that $\rho_{0}$ is bounded from below
and above by positive constants on $X.$ Then 
\[
\phi^{(N)}(x):=\frac{1}{k}\int_{\R^{n(N-1)}}\log(\mbox{Per}(x,x_{2},...,x_{N}))\rho_{0}(x_{2})dx_{2}\cdots\rho_{0}(x_{N})dx_{N}-c_{N},
\]
 where $c_{N}$ is the normalizing constant ensuring that $\int_{\R^{n}}\phi^{(N)}\rho_{0}dx=0,$
\textup{converges}, as $N\rightarrow\infty,$\textup{ locally uniformly
to the unique convex function $\phi$ }solving the second boundary
value problem with target $P$ for \textup{the Monge-Ampère equation
$MA(\phi)=\mu_{0}$ with the normalization condition $\int_{\R^{n}}\phi\mu_{0}=0.$
Moreover, $T^{(N)}(x):=\nabla\phi^{(N)}(x):=$ 
\[
=\frac{1}{k}\int_{\R^{n(N-1)}}\frac{\sum_{\sigma\in S_{N}}p_{\sigma(1)}e^{x\cdot p_{\sigma(1)}+x_{2}\cdot p_{\sigma(2)}+\cdots+x_{N}\cdot p_{\sigma(N)}}}{\sum_{\sigma\in S_{N}}e^{x\cdot p_{\sigma(1)}+x_{2}\cdot p_{\sigma(2)}+\cdots+x_{N}\cdot p_{\sigma(N)}}}\rho_{0}(x_{2})dx_{2}\otimes\cdots\rho_{0}(x_{N})dx_{N}
\]
converges point-wise, in the interior of $X,$ to the (Hölder continuous)
optimal map $T$ for the Monge problem of transporting the probability
measure $\mu_{0}$ on $X$ to the uniform probability measure $\lambda_{P}$
on the target convex body $P.$ }
\end{cor}
The existence of the optimal map $T$ in the previous corollary is
due to Brenier \cite{br} and the Hölder regularity was shown by Caffarelli
\cite{ca-2}. As explained in section \ref{sec:General-target-measures}
a variant of the previous setting can also be considered which in
particular applies to more general target measures $\nu.$ 

Coming back to Theorem \ref{thm:main for perman intro} we point out
that the key point in its proof is a rather general argument which,
in a sense, reduces the problem to the case when $\beta=\infty.$
We can then take $\beta_{N}=k,$ so that the corresponding random
point process is exactly permanental. In that case the theorem turns
out to be a rather immediate consequence of the large deviation principle
for\emph{ determinantal} point processes on polarized \emph{complex}
manifolds proved in \cite{be3}, building on \cite{ber-bou,b-b-w}.
More precisely, the latter result applies in the case when $P$ is
a rational convex polytope, the point being that $P$ then determines
a polarized toric variety $(X_{p},L_{P})$ with a projection map to
$\R^{n}$ (see \cite{ber-ber} and references therein). Anyway, we
will give a direct purely ``real'' proof in the present setting
for the case $\beta_{N}=k$ (the key ingredient is Prop \ref{prop:asymp of free energy for permanant}).
As for the reduction to the case $\beta=\infty,$ it is inspired by
some ideas from statistical mechanics and in particular mean field
theory. In physical terms the idea of the argument may be explained
as follows. Imagine that we know the macroscopic ground state (i.e.
the state of zero energy $E$) of a system of a large number $N$
of particles in thermal equilibrium at zero temperature (i.e. at $\beta=\infty)$.
If we can rule out any \emph{first order phase transitions} at zero-temperature
(which essentially means that the macroscopic equilibrium state is
unique), then increasing the temperature (i.e decreasing $\beta)$
leads to a new macroscopic equilibrium state, minimizing the corresponding\emph{
free energy} functional $E-S/\beta,$ where $S$ is the physical entropy
(i.e. $S=-D$ with our sign conventions). Increasing the temperature
thus gives a transition from an ordered zero-temperature macroscopic
state to a disordered positive temperature macroscopic state. 

However, to apply this reasoning to a the present situation we have
to deal with $N-$particle interacting systems with rather general
interactions ($N-$point Hamiltonians). One of the difficulties that
we then have to confront is that the Hamiltonian in question is not
a sum of two-point functions, as opposed to the more standard situation
studied rigorously by Messer-Spohn \cite{m-s}. This is a reflection
of the fact that the corresponding field equations \ref{eq:ma eq intro cassical}
are fully non-linear (i.e. non-linear in the derivative terms). More
generally: the Hamiltonian is not given in the usual \emph{mean field}
form, i.e. it is not of the functional form $H(\delta_{N})$ for some
$N-$independent functional $H$ on the space of probability measures
on $X$ \cite{bo}. Still, as explained in the next section we can
obtain a rather general convergence result of mean field type for
such interacting particle systems, that is hopefully of independent
interest. As will be explained below this turns out to be related
to previous work of Ellis-Have-Turkington \cite{e-h-t}. The result
applies in particular to the setting above where the corresponding
Hamiltonian may (in the non-weighted case) be written as $H^{(N)}(x_{1},...,x_{N}):=-\frac{1}{k}\log\mbox{Per}(x_{1},...,x_{N}),$
i.e. 
\begin{equation}
H^{(N)}(x_{1},...,x_{N})=-\frac{1}{k}\log\sum_{\sigma\in S_{N}}e^{-kNC(\sigma)},\,\,\, C(\sigma):=-\frac{1}{k}(x_{1}\cdot p_{\sigma(1)}+\cdots x_{N}\cdot p_{\sigma(N)})/N\label{eq:ham in intro}
\end{equation}
As a sideremark we note that this form of writing the Hamiltonian
gives a simple heuristic interpretation of the large deviation result
for $\beta_{N}=k;$ the result essentially says that we may, when
$N\rightarrow\infty,$ replace the whole sum over all permutations
with the contribution from the permutation with minimal cost $C(\sigma)$
(compare Remark \ref{rem:log per is cost asympt} and the relation
between optimal transport and its discrete version described in the
appendix). But the actual proof proceeds in a different way, based
on a duality argument, where the cost functional $C(\mu)$ arises
as the Legendre transform of another functional on the space $C_{b}(X)$
of all bounded continuous functions on $X.$

As explained below, for a fixed $\beta,$ the Monge-Ampère measure
$\mu_{\beta}(=MA(\phi_{\beta}),$ where $\phi_{\beta}$ is the solution
of equation \ref{eq:ma eq weak intro} arises as the minimizer of
the corresponding free energy functional. From this point of view
it is interesting to study the behaviour of the measures $\mu_{\beta},$
describing the equilibrium states at inverse temperature $\beta,$
as $\beta$ varies. For example, when the fixed background measure
$\mu_{0}$ is equal (or comparable) to the usual Euclidean measure
$dx$ the support of $\mu_{\beta}$ is all of $\R^{n}.$ However,
the limiting measure obtained when $\beta$ increases to infinity
is always supported on a\emph{ compact} set, whose ``boundary''
appears as the free boundary for a Monge-Ampère equation, as explained
above. This shows that the the phase transition at $\beta=0$ referred
to above is, in the present setting, reminiscent of a liquid-gas phase
transition.

Before continuing we also point out that the main new technical difficulty
which appears in the complex geometric setting out-lined in \cite{be-3},
where the role of the permanent above is played by a Vandermonde type
determinant is that the corresponding Hamiltonian $H^{(N)}(z_{1},...,z_{N})$
is singular (for example when different points merge). See \cite{clmp,k}
for the case of dimension $n=1,$ or more precisely the case when
$X$ is a domain in the complex plane $\C.$ A major simplifying feature
which appears in the case $n=1$ is that the Vandermonde determinant
factorizes completely and the point process in question is hence a
\emph{Coulomb gas,} i.e. $H^{(N)}(z_{1},...,z_{N})$ is then proportional
to a sum of two-point functions of the form $\log|z_{i}-z_{j}|$ (compare
the end of section \ref{sub:A-general-LDP})

\subsection{Interacting particle system in thermal equilibrium}

Let $X$ be a closed set in $\R^{n}$ and assume, for simplicity,
that $X$ is compact (generalizations to the non-compact case of $\R^{n}$
will be considered in section \ref{sec:Permanental-point-processes,}).
Fix a probability measure $\mu_{0}$ supported on $X.$ For a fixed
positive integer $N$ (representing the number of particles) we assume
given an\emph{ $N-$particle Hamiltonian} $H^{(N)},$ i.e. a continuous
function on the $N-$fold product $X^{N}$ which is symmetric, i.e.
invariant under the action of the permutation group $S_{N}.$ The
corresponding\emph{ Gibbs measure} is the probability measure on $X^{N}$
defined by 

\[
\mu_{\beta_{N}}^{(N)}:=e^{-\beta_{N}H^{(N)}}\mu_{0}^{\otimes N}/Z_{N,\beta_{N}},
\]
 where we have also fixed a positive number $\beta_{N}$ (the\emph{
inverse temperature}) and where $Z_{N}$ is the normalizing constant
(\emph{partition function}), i.e. 
\[
Z_{N,\beta_{N}}:=\int_{X^{N}}e^{-\beta_{N}H^{(N)}}\mu_{0}^{\otimes N}
\]
We will also assume that $H^{(N)}$ is\emph{ uniformly Lipschitz continuous
}(in fact, \emph{equicontinuous} will be enough; compare section \ref{sub:Setup-and-assumptions on hamilton}).
In the case when $H^{(N)}$ is differentiable this thus simply means
that there is a constant $L$ independent of $N$ such that 
\[
\sup_{X}\left|\partial_{x_{i}}H^{(N)}(x_{1}....x_{N})\right|\leq L
\]
Given a continuous function $u$ on $X$ we denote by $Z_{N,\beta_{N}}[u]$
the ``tilted'' partition function obtained by replacing $H^{(N)}$
with $H^{(N)}+u,$ where $u(x_{1},..,x_{N}):=\sum u(x_{i}).$ In other
words $Z_{N,\beta_{N}}[u]$ is the (scaled) Laplace transform of the
law of the empirical measure $\delta_{N}$ defined by the Gibbs measure
associated to $H_{N}.$ 
\begin{thm}
\label{thm:interacting part intro}Let $H^{(N)}$ an $N-$particle
Hamiltonian as above. \textup{Assume that there exists a sequence
$\beta_{N}$ of positive real numbers tending to infinity such that
}$-\frac{1}{\beta_{N}N}\log Z_{N,\beta_{N}}[u]$ converges, when $N\rightarrow\infty,$
to a functional $\mathcal{F}(u)$ which is Gateaux differentiable
on $C^{0}(X).$ Then, for any fixed $\beta>0,$ the law of the empirical
measure $\delta_{N}$ for the random point process on $X$ defined
by the corresponding Gibbs measure $\mu_{\beta}^{(N)}$ satisfies
a LDP with speed $\beta N$ and good rate functional 
\[
F_{\beta}(\mu)=E(\mu)+\frac{1}{\beta}D_{\mu_{0}}(\mu)-C_{\beta}
\]
 where $C_{\beta}$ is the normalizing constant and where the functional
$E(\mu),$ which only depends on the support $X$ of $\mu_{0}$, denotes
the Legendre transform of $\mathcal{F}(u)$ (formula \ref{eq:def of E of mu as legendre transform})
and $D_{\mu_{0}}(\mu)$ is the entropy of $\mu$ relative to $\mu_{0}.$
In particular, $\delta_{N}$ converges in probability to the unique
minimizer $\mu_{\beta}$ of $F_{\beta}.$
\end{thm}
Moreover, we will show that the minimizer $\mu_{\beta}$ can be written
as $\mu_{\beta}=d\mathcal{F}_{|u_{\beta}}$ where $u_{\beta}$ is
a continuous function on $X$ solving the following equation of mean
field type: 
\begin{equation}
d\mathcal{F}_{|u}=e^{\beta u}\mu_{0},\label{eq:general mean field eq intro}
\end{equation}
 where $d\mathcal{F}_{|u}$ denotes the measure defined by the Gateaux
differential of $\mathcal{F}$ at $u$ (see Theorem \ref{thm:existence of sol to general mean f}).
We will also obtain a canonical sequence $u_{N}$ of functions uniformly
converging to $u,$ which are the unique solutions of certain ``finite
$N"$ approximations to the mean field type equation above and which
arise as the limiting fixed point of certain iterations (see section
\ref{sub:Balanced-functions-and}). The construction of $u_{N}$ is
inspired by Donaldson's notion of balanced metrics, introduced in
the setting of Kähler-Einstein metrics \cite{do}) and is hopefully
of independent interest.

Theorem \ref{thm:interacting part intro} generalizes the seminal
result of Messer-Spohn \cite{m-s} and is closely related to a previous
result of Ellis-Haven-Turkington \cite{e-h-t} (compare section \ref{sub:An-alternative-proof}).
In the statement of the theorem we have made the rather strong assumptions
of \emph{(a)} compactness of $X$ and \emph{(b)} equicontinuity of
$H^{(N)}.$ In fact, under these assumptions a more direct proof can
be given by essentially reducing the problem to the classical LDP
in the non-interacting case (i.e. $H^{(N)}=0$ or equivalently $\beta=0)$
which is the content of Sanov's theorem (see section \ref{sub:An-alternative-proof}
where a comparison with the work \cite{e-h-t} is also made). But
the main point here is to give a flexible proof that can be adapted
to more general situations where the assumptions (a) and (b) may not
be satisfied. For example, in the setting of Theorem \ref{thm:main for perman intro}
above the the compactness assumption (a) does not hold and in the
complex geometric setting refereed to above the assumption (b) is
not satisfied, since $H^{(N)}$ is even singular (compare section
\ref{sub:A-general-LDP} ). 

The key step in the proof of the previous theorem is to establish
the convergence of the limiting \emph{mean energy} towards the Legendre
transform of $\mathcal{F}:$ 
\begin{equation}
\lim_{N\rightarrow\infty}\int_{X^{N}}\frac{H^{(N)}}{N}\mu^{\otimes N}=\mathcal{F}^{*}(\mu)\label{eq:conv of mean energy intro}
\end{equation}
(Theorem \ref{thm:existence of mean energ}). Then applying the finite
dimensional Gibbs variational principle and some variational arguments
allows us to compute the limit of the corresponding partition functions,
``tildted by $u".$ Differentiating at $u=0$ this is implies in
particular the convergence of the one-point correlation measures towards
the minimizer of the free energy functional. Finally, a suitable application
of the (generalized) Gärtner-Ellis theorem gives the LDP principle
in the theorem. This approach is inspired by the approach of Messer-Spohn
\cite{m-s}. But the key new observation here is that (even if the
assumptions (a) and (b) are not satisfied) the desired LDP holds as
long as the convergence in formula \ref{eq:conv of mean energy intro}
holds (in fact, the lower bounds always holds). More precisely, it
is enough if the convergence holds for all $\mu$ in the subspace
$\{F_{\beta}<\infty\}$ as long as the latter subspace is dense. 

Another point of the approach developed here, as compared to \cite{m-s},
is to separate the case $\beta>0$ from the more subtle case $\beta<0,$
by bypassing the use of the Hewitt-Savage decomposition theorem used
in \cite{m-s}. In fact, if one uses the latter decomposition theorem
then one obtains, in the spirit of \cite{m-s}, a generalization of
Theorem \ref{thm:interacting part intro} to possibly negative $\beta,$
saying that, after passing to a subsequence, the law of the empirical
measures converge weakly to a probability measure on the space of
all probability measures which is concentrated on the minimizers of
the functional $F_{\beta}.$ The point is that in the case $\beta>0$
the functional $F_{\beta}$ is strictly convex and hence admits a\emph{
unique }minimizer, as in Theorem \ref{thm:interacting part intro}.
However, extensions to the non-compact setting of $\R^{n},$ in the
case $\beta<0,$ require refined growth estimates and we thus leave
this case for a separate paper where applications to Kähler-Einstein
geometry on toric varieties will be given \cite{be-4}. For a brief
outline of the case $\beta<0$ in the Monge-Ampère setting in $\R^{n}$
and the relations to the existence problem for Kähler-Einstein metrics
on complex algebraic varieties and phase transitions, see section
\ref{sub:K=0000E4hler-Einstein-metrics,-negativ}.

\subsubsection*{Outline of the paper}

After having introduced the precise assumptions on the $N-$particle
Hamiltonian in section \ref{sec:Existence-of-the} we prove the existence
of the limiting mean energy for measures with continuous potentials.
This assumption is automatically satisfied in the Monge-Ampère setting
and the existence in the general case is proved in section \ref{sec:Generalized-mean-field}
by establishing an approximation result of independent interest, inspired
by Donaldson's notion of balanced metrics. Then in section \ref{sec:The-Large-Deviation}
we go on to establish the LDP in Theorem \ref{thm:interacting part intro}
and we also give an alternative more direct proof of the latter theorem.
In section \ref{sec:Permanental-point-processes,} we introduce the
Monge-Ampère setting leading up to a proof of Theorem \ref{thm:main for perman intro}.
In section \ref{sec:General-target-measures} a variant of the previous
setting is considered which in particular applies to rather general
target measures. A comparison with the complex geometric setting and
an outlook on further developments is given in section \ref{sec:Outlook}.
The paper is concluded with an appendix giving some background on
optimal transport and establishing the (essentially well-known) comparison
principle in the Monge-Ampère setting.

\subsubsection*{Acknowledgments}

Thanks to Ofer Zeitouni for helpful comments on a first draft of this
paper and Bo'az Klartag for stimulating discussions. As pointed out
in the introduction, this paper and in particular the connections
to the theory of optimal transport arose as a ``spin-off effect''
of some of my work in complex geometry and pluripotential theory.
Lacking background in optimal transport theory I apologize for any
omission of accrediting prior results properly. This research was
supported by grants from the European Research Council (ERC starting
grant 307529) and the Swedish Research Council.

\section{\label{sec:Existence-of-the}Existence of the limiting mean energy}

\subsection{\label{sub:Setup-and-assumptions on hamilton}Setup and assumptions
on the Hamiltonian $H^{(N)}$}

Let $(X,d,\mu_{0})$ be a compact metric space with a probability
(Borel) measure $\mu_{0.},$ where $X$ is the support of $\mu_{0}.$
The typical situation is when $X$ is embedded in a Riemannian manifold
and $d$ is induced from the Riemannian metric.

Let $H^{(N)}$ be an \emph{$N-$particle Hamiltonian,} i.e. continuous
function on the $N-$fold product $X^{N}$ which is symmetric, i.e.
invariant under the action of the permutation group $S_{N}.$ We denote
by $d_{X^{N}}$ the following induced distance function on $X^{N}:$
\[
d_{X^{N}}((x_{1},...x_{N}),y_{1},...,y_{N})):=\frac{1}{N}\sum_{i=1}^{N}d_{X}(x_{i},y_{i})
\]
We will assume that $H^{(N)}/N$ is \emph{equicontinuou}s, i.e. given
any $\epsilon>0$ there exists a $\delta>0$ such that $|H^{(N)}(x)-H^{(N)}(y))|/N\leq\epsilon$
if $d_{X^{N}}(x,y)\leq\delta.$ It will often be convenient to make
the stronger assumption that $H^{(N)}/N$ be \emph{uniformly Lipschitz
continuous} in the sense that, there exists a constant $L$ such that
\begin{equation}
|H^{(N)}(x)-H^{(N)}(y))|/N\leq Ld_{X^{N}}(x,y)\label{eq:Lip assumption}
\end{equation}
In particular, the latter property holds if $H^{(N)}(x,x_{2},...,x_{N})$
is Lipschitz continuous on $X$ with Lipschitz constant $L,$ for
any fixed $(x_{2},...,x_{N})\in X^{N-1}.$ 

We will also assume that $H^{(N)}$ has the following\emph{ limiting
properties:} 
\begin{itemize}
\item The following limit exists for any $u\in C^{0}(X):$ 
\begin{equation}
\mathcal{F}(u):=\lim_{N\rightarrow\infty}\inf_{X^{N}}\frac{1}{N}(H^{(N)}+u),\label{eq:def of f of u as limit of ham}
\end{equation}
where we have used the notation $u(x_{1},..,x_{N}):=\sum_{i}u(x_{i}).$
\item The functional $\mathcal{F}$ is Gateaux differentiable on $C^{0}(X)$
(i.e. differentiable along all affine lines in $C^{0}(X))$
\end{itemize}
It will turn out that, from the statistical mechanical point of view,
the first part is equivalent to the ``existence of the free energy''
(compare Lemma \ref{lem:inf is Lbeta norm}).

We will write $X^{(N)}:=X/S_{N},$ where $S_{N}$ is the permutation
group with $N$ elements and equip it with the induced quotient distance
function 
\[
d_{X^{(N)}}((x_{1},...x_{N}),y_{1},...,y_{N})):=\inf_{\sigma\in S_{N}}\frac{1}{N}\sum_{i=1}^{N}d_{X}(x_{i},y_{\sigma(i)})
\]
and denote by $\mathcal{M}_{1}(X)$ the space of all probability measure
$\mu$ on $X.$ Similarly, we denote by $\mathcal{M}_{(N)}(X^{N})$
the $S_{N}-$invariant subspace of $\mathcal{M}_{1}(X^{N}).$ There
is a standard embedding 
\begin{equation}
X^{(N)}\hookrightarrow\mathcal{M}_{1}(X),\,\,\,\,(x_{1},..,x_{N})\mapsto\delta_{N}:=\frac{1}{N}\sum\delta_{x_{i}}\label{eq:def of empricical measure}
\end{equation}
where we will call $\delta_{N}$ the\emph{ empirical measure. }We
equip the space $\mathcal{M}_{\text{1 }}(X)$ with the Wasserstein
1-metric $d_{W}$ determined by $(X,d).$ Then it is well-known that
the embedding above becomes an isometry (see the end of section \ref{sub:Optimal-transport-theory}),
i.e. 

\begin{equation}
d_{W}(\delta_{N}(x),\delta_{N}(x'))=d_{X^{(N)}}(x,x')\label{eq:wasserstein metric isom}
\end{equation}
When $X$ is compact (which we assume in this section) the topology
on $\mathcal{M}_{\text{1 }}(X)$ induced by $d_{W}$ coincides with
the usual topology defining weak convergence of measures.

\subsection{\label{sub:Existence-of-the}Existence of the limiting mean energy}

Given an element $\mu$ in the space $\mathcal{M}_{1}(X)$ of probability
measures on $X$ we define its \emph{thermodynamic energy} $E(\mu)$
as the Legendre transform of $\mathcal{F}$ at $\mu$ (with a somewhat
non-standard sign convention): 
\begin{equation}
E(\mu):=\sup_{u\in C^{0}(X)}(\mathcal{F}(u)-\int_{X}u\mu):=\sup_{u\in C^{0}(X)}\mathcal{F}_{\mu}\label{eq:def of E of mu as legendre transform}
\end{equation}
The key point is the following result which identifies $E$ with the
limiting\emph{ mean energy }(this terminology will be taken up in
section \ref{sec:The-Large-Deviation}) under the assumption that
$\mathcal{F}(u)$ be Gateaux differentiable.
\begin{thm}
\label{thm:existence of mean energ}Let $H^{(N)}$ be an $N-$particle
Hamiltonian satisfying the properties in section \ref{sub:Setup-and-assumptions on hamilton},
then 
\[
\lim_{N\rightarrow\infty}\frac{1}{N}\int_{X^{N}}H^{(N)}\mu^{\otimes N}=E(\mu)
\]
 for any probability measure $\mu$ on $X.$
\end{thm}

\subsection{The proof of Theorem \ref{thm:existence of mean energ}}

Let us start by fixing $\mu$ and $u\in C^{0}(X)$ and rewriting 
\begin{equation}
\frac{1}{N}\int_{X^{N}}H^{(N)}\mu^{\otimes N}=\left(\inf_{X^{N}}\frac{H^{(N)}+u}{N}-\int_{X}u\mu\right)+I_{N}[\mu,u],\label{eq:pr of exist of mean energ first formula}
\end{equation}
 where 
\[
I_{N}[\mu,u]:=\int_{X^{N}}(\frac{H^{(N)}+u}{N}-\inf_{X^{N}}\frac{H^{(N)}+u}{N})\mu^{\otimes N}
\]
Since, trivially, $I_{N}[\mu,u]\geq0$ it follows that the\emph{ lower}
bound in the theorem to be proved always holds: 
\[
E(\mu):=\sup_{u\in C^{0}(X)}\left(\inf_{X^{N}}\frac{H^{(N)}+u}{N}-\int_{X}u\mu\right)\leq\liminf_{N\rightarrow\infty}\frac{1}{N}\int_{X^{N}}H^{(N)}\mu^{\otimes N},
\]
To handle the upper bound we will use the assumed differentiability
of the functional $\mathcal{F}(u).$ Let us start by introducing the
following terminology: 
\begin{defn}
A function $u_{\mu}$ on $X$ is said to be a\emph{ potential of the
measure $\mu$} if $u_{\mu}$ is a maximizer of the functional whose
sup defines $E(\mu),$ i.e. if 
\begin{equation}
E(\mu):=\mathcal{F}(u_{\mu})-\int_{X}u_{\mu}\mu\label{eq:E of mu in terms of potential}
\end{equation}
Since $\mathcal{F}$ is concave and also assumed Gateaux differentiable
this equivalently means that 
\begin{equation}
\mu=d\mathcal{F}_{|u_{\mu}}.\label{eq:mu as deriv in terms of pot}
\end{equation}

\end{defn}
The key upper bound that we will prove may be formulated as the following
\begin{prop}
\label{prop:upper bound}Assume that the functional $\mathcal{F}$
is Gateaux differentiable. Then, for any probability measure $\mu$
admitting a continuous potential we have 
\[
\limsup_{N\rightarrow\infty}\frac{1}{N}\int_{X^{N}}H^{(N)}\mu^{\otimes N}\leq E(\mu)
\]

\end{prop}
In order to prove the previous proposition we start with the following
simple but very useful lemma (which was used in the similar context
of Fekete points in \cite{b-b-w}).
\begin{lem}
\label{lem:conv of abstr fekete}Fix $u_{*}\in C^{0}(X)$ and assume
that $x_{*}^{(N)}\in X^{N}$ is a minimizer of the function $(H^{(N)}+u_{*})/N$
on $X^{N}.$ If the corresponding large $N-$ limit $\mathcal{F}(u)$
exists for all $u\in C^{0}(X)$ and \emph{$\mathcal{F}$ is Gateaux
differentiable at $u_{*},$} then $\delta_{N}(x_{*}^{(N)})$ converges
weakly towards $\mu_{*}:=d\mathcal{F}_{|u_{*}}.$\end{lem}
\begin{proof}
Fix $v\in C^{0}(X)$ and a real number $t.$ Let $f_{N}(t):=\frac{1}{N}(H^{(N)}+u+tv)(x_{*}^{(N)})$
and $f(t):=\mathcal{F}(u+tv).$ By assumption $\lim_{N\rightarrow}f_{N}(0)=f(0)$
and $\liminf_{N\rightarrow}f_{N}(t)\geq f(t).$ Note that $f$ is
a concave function in $t$ (since it is defined as an inf of affine
functions) and $f_{N}(t)$ is affine in $t.$ But then it follows
from the differentiability of $f$ at $t=0$ that $\lim_{N\rightarrow\infty}df_{N}(t)/dt_{|t=0}=df(t)/dt_{|t=0},$
i.e. that 
\[
\lim_{N\rightarrow\infty}\left\langle \delta_{N}(x_{*}^{(N)}),v\right\rangle =\left\langle d\mathcal{F}_{|u},v\right\rangle ,
\]
 which thus concludes the proof of the lemma.
\end{proof}
Let us also recall the following weak form of Sanov's theorem: 
\begin{lem}
\label{lem:conv in law equiv to conv of moments}For any given $\mu\in\mathcal{M}_{\text{1 }}(X)$
we have that $(\delta_{N})_{*}\mu^{\otimes N}$ converges to $\delta_{\mu}$
in the space $\mathcal{M}_{\text{1 }}(X)$ equipped with the weak
topology.\end{lem}
\begin{proof}
In fact, this is an immediate consequence of the following useful
result: if $\mu_{N}\in\mathcal{M}_{(N)}(X^{N})$ then $(\delta_{N})_{*}\mu_{N}\rightarrow\delta_{\mu}$
iff the corresponding $j-$point correlation measures $\int_{X^{N-j}}\mu_{N}$
converge weakly to $\mu^{\otimes j}$ (see Prop 2.2 in \cite{sn}).
\end{proof}

\subsubsection{The proof of Prop \ref{prop:upper bound}. }

By Theorem \ref{thm:exist of cont potentials} below any probability
measure $\mu$ admits a continuous potential $u_{\mu.}$ Given the
measure $\mu\in\mathcal{M}_{\text{1 }}(X)$ we may thus assume that
the potential $u_{\mu}$ is continuous. Since $X$ is compact this
means that for any given $\epsilon>0$ there exists $\delta>0$ such
that $d_{X}(x,x')<2\delta$ implies that $|u_{\mu}(x)-u_{\mu}(x')|\leq\epsilon.$
Moreover, by the assumption on $H^{(N)}/N$ the corresponding property
of $H^{(N)}/N$ holds uniformly in $N$ on $X^{N}.$ In particular,
\begin{equation}
d_{X^{N}}(x,x')\leq2\delta\Rightarrow|(\frac{H^{(N)}}{N}+u_{\mu})(x)-(\frac{H^{(N)}}{N}+u_{\mu})(x')|\leq2\epsilon\label{eq:uniform cont of perturbed H}
\end{equation}
We denote by $B_{\delta}(\mu)$ the inverse image in $X^{N}$ of a
ball of radius $\delta$ centered at $\mu$ in $(\mathcal{M}_{\text{1 }}(X),d_{W}),$
under the map $\delta_{N}.$ Let us first note that 
\begin{equation}
\lim_{N\rightarrow\infty}\int_{X^{N}-B_{\delta}(\mu)}\frac{H^{(N)}+u_{\mu}}{N}\mu^{\otimes N}=0\label{eq:decay of mean energ}
\end{equation}
Indeed, since $\frac{H^{(N)}+u}{N}$ is uniformly bounded this follows
immediately from the previous lemma. Let now $x_{*}^{(N)}$ be a sequence
as in Lemma \ref{lem:conv of abstr fekete}, associated to $u=u_{\mu}.$
By the lemma we have that $\delta_{N}(x_{*}^{(N)})$ converges weakly
to $\mu.$ Hence, taking $N$ sufficiently large ($N\geq N_{\delta})$
we may as well assume that $x_{*}^{(N)}\in B_{\delta}(\mu).$ As a
consequence, if $x$ is any given point in $B_{\delta}(\mu)$ then,
by the isometry property \ref{eq:wasserstein metric isom}, $d_{X^{(N)}}(x,x_{*}^{(N)})\leq2\delta.$
But then \ref{eq:uniform cont of perturbed H} gives 
\begin{equation}
\int_{B_{\delta}(\mu)}(\frac{H^{(N)}+u}{N}-\inf_{X^{N}}\frac{H^{(N)}+u}{N})\leq2\epsilon\label{eq:decays on small ball}
\end{equation}
Combining \ref{eq:decay of mean energ} and \ref{eq:decays on small ball}
we finally deduce that $\limsup_{N\rightarrow}I_{N}[\mu,u_{\mu}]\leq\epsilon,$
for any $\epsilon>0,$ which thus concludes the proof of Prop \ref{prop:upper bound}.

\subsection{\label{sub:The-mean-energy}The mean energy as the Legendre transform
of the limiting free energy}

It will be very useful to obtain an alternative integral expression
for the functional $\mathcal{F}(u).$ To this end we fix a probability
measure $\mu_{0}$ with support $X$ and a sequence $\beta_{N}\rightarrow\infty$
of positive numbers and set 
\begin{equation}
\mathcal{F}^{(N)}(u):=-\frac{1}{\beta_{N}N}\log\int_{X^{N}}e^{-\beta_{N}N(H^{(N)}/N+u)}\mu_{0}^{\otimes N},\label{eq:def of N particle free energy}
\end{equation}
where we have omitted the explicit dependence on $\mu_{0}$ in the
notation (in terms of the notation introduced in section we thus have
the $\mathcal{F}^{(N)}$ may be written as $\mathcal{F}^{(N)}(u)=-\frac{1}{\beta_{N}N}\log Z_{N,\beta_{N}}[u]).$
\begin{lem}
\label{lem:inf is Lbeta norm}Assume that $H^{(N)}/N$ is equicontinuous
(as defined in section \ref{sub:Setup-and-assumptions on hamilton})
and that $u$ is continuous on $X.$ Then, as $N\rightarrow\infty,$
\textup{
\[
\inf_{X^{N}}\frac{1}{N}(H^{(N)}+u)=\mathcal{F}^{(N)}(u)+o(1)
\]
for any sequence $\beta_{N}\rightarrow\infty.$ }\end{lem}
\begin{proof}
Setting $G^{(N)}:=-(H^{(N)}+u)/N,$ then it will be enough to show
that for any $\epsilon>0$ there exists a constant $C_{\epsilon}$
such that 
\[
\frac{1}{\beta_{N}N}\log\int_{X^{N}}e^{\beta_{N}NG^{(N)}}\mu_{0}^{\otimes N}\geq\sup_{X^{N}}G^{(N)}-C_{\epsilon}/\beta_{N}-\epsilon,
\]
To this end we denote by $x_{*}^{(N)}$ a configuration where the
sup in the rhs above is attained and restrict the integration to a
polydisc $\Delta_{\delta}(x_{*}^{(N)})$ of radius $\delta$ centered
at $x_{*}^{(N)}.$ Thus 
\[
\frac{1}{\beta_{N}N}\log\int_{X^{N}}e^{\beta_{N}N(G^{(N)}}\mu_{0}^{\otimes N}\geq\sup_{X^{N}}G^{(N)}+\frac{1}{\beta_{N}N}\log\int_{\Delta_{\delta}(x_{*}^{(N)})}e^{\beta_{N}N(G^{(N)}-G^{(N)}(x_{*}^{(N)}))}
\]
Now, by the uniform continuity assumption, given $\epsilon>0$ we
can take $\delta>0$ such that 
\[
\frac{1}{\beta_{N}N}\log\int_{\Delta_{\delta}(x_{*}^{(N)})}e^{\beta_{N}N(G^{(N)}-G^{(N)}(x_{*}^{(N)})}\geq-\epsilon+\frac{1}{\beta_{N}N}\log\int_{\Delta_{\delta}(x_{*}^{(N)})}\mu_{0}^{\otimes N}
\]
 Since, by assumption, the ``coordinates'' of $x_{*}^{(N)}$ are
contained in the support of $\mu_{0}$ the last term above may, for
any $\delta>0$ be estimated from below by $\frac{1}{\beta_{N}N}\log(C_{\delta})^{N}$
which thus concludes the proof (a similar argument is used in the
proof of Lemma \ref{lem:b-m for convex etc}).
\end{proof}
By the previous Lemma the first point in the assumptions about the
limiting properties of $H^{(N)}$ (in section \ref{sub:Setup-and-assumptions on hamilton})
is thus equivalent to the existence of the limit 
\[
\mathcal{F}(u):=\lim_{N\rightarrow\infty}\mathcal{F}^{(N)}(u)
\]
and in particular the limit only depends on the support $X$ of $\mu_{0}$
and not on $\mu_{0}$ itself.

\section{\label{sec:Generalized-mean-field}Generalized mean field equations
and balanced functions}

In this section we will, among other thing, prove that any probability
measure admits a continuous potential, which was used in the proof
of Theorem \ref{thm:existence of mean energ}. The theorem will be
deduced from a general result of independent interest (Theorem \ref{thm:conv of balanced})
which will allow us to write $u_{\mu}$ as a $C^{0}-$limit of equicontinuous
functions $u_{N}$ maximizing the functional 
\begin{equation}
\mathcal{F}_{\mu}^{(N)}(u):=\mathcal{F}^{(N)}(u)-\int u\mu\label{eq:N-particle free energy with mu}
\end{equation}
In fact, this argument will give a more precise result saying that
$u_{\mu}$ is in a certain subspace $\mathcal{P}(X)$ of $C^{0}(X)$
which, in the setting of the real Monge-Ampère operator may be identified
with the space of all convex functions whose gradient image is contained
in the given target convex body $P.$ In the latter setting the functions
$u_{N}$ essentially coincide with Donaldson's $\mu-$balanced metrics
(see \cite{do,bbgz}). 

This approach will also lead to a natural setting for formulating
global equations, that we will refer to as\emph{ generalized mean
field equations,} which generalize the real Monge-Ampère equation
\ref{eq:ma eq weak intro}.

\subsection{Setup}

To simplify the exposition of the proof it will be convenient to assume
that $H^{(N)}$ is uniformly Lipschitz continuous with Lipschitz constant
$L,$ but the proof under the more general assumption of uniform continuity
is essentially the same. Let $L(X)$ be the space of all Lipschitz
continuous functions on $X$ with Lipschitz constant $L.$ We let
$\mathcal{P}_{N}(X)$ be the subspace of $L(X)$ of all functions
$u$ such that there exists a finite measure $\nu$ on $X^{N}$ such
that 
\[
u(x)=\frac{1}{\beta_{N}}\log\int_{y\in X^{N-1}}e^{-\beta_{N}H^{(N)}(x,y)}\nu
\]
Next, we define $\mathcal{P}(X)$ as the closure in $C^{0}(X)$ of
the union of all spaces $\mathcal{P}_{N}(X)$ as $N$ ranges over
all positive integers. By construction we thus have
\[
\mathcal{P}_{N}(X)\subset\mathcal{P}(X)\subset L(X)\subset C^{0}(X)
\]

\subsection{\label{sub:Balanced-functions-and}Balanced functions and their large
$N-$limit}

Let $\pi_{N}$ be the operator 
\begin{equation}
\pi_{N}:\, C^{0}(X)\rightarrow\mathcal{P}_{N}\label{eq:mapping property}
\end{equation}
defined by 
\begin{equation}
\pi_{N}(u)(x):=\frac{1}{\beta_{N}}\log\frac{1}{Z_{N}[u]}\int_{y\in X^{N-1}}e^{-(\beta_{N}H^{(N)}(x,y)+u(y))}\mu_{0}^{\otimes(N-1)},\label{eq:def of finte dim proj}
\end{equation}
 where, as usual, $Z_{N}[u]:=\int_{X^{N}}e^{-(\beta_{N}H^{(N)}(x)+u(x))}\mu_{0}^{\otimes N}.$
The definition of $\pi_{N}(u)$ is made so that $e^{\beta_{N}(\pi_{N}(u)-u)}\mu_{0}$
is a probability measure (in fact, $e^{\beta_{N}(\pi_{N}(u)-u)}\mu_{0}$
coincides with the one-point correlation measure $\mu_{1}^{(N)}$
of the corresponding Gibbs measure; compare section \ref{sub:Setup:-the-Boltzmann-Gibbs}). 
\begin{defn}
A function $u_{N}$ on $X$ is said to be\emph{ balanced with respect
to $(\mu_{0},\beta_{N})$ }if $\pi_{N}(u_{N})=u_{N}.$ \end{defn}
\begin{prop}
\label{prop:existence of balanced}Given a probability measure $\mu_{0}$
on $X$ there exists, for any integer $N,$ a function $u_{N}\in\mathcal{P}_{N}(X)$
which is balanced \emph{with respect to $(\mu_{0},\beta_{N}).$ }Moreover,
$u_{N}$ maximizes the functional $\mathcal{F}_{\mu_{0}}^{(N)}(u)$
 on $C^{0}(X)$ and is uniquely determined mod $\R.$\end{prop}
\begin{proof}
First observe that, by definition, $\pi_{N}(u+c)=\pi_{N}(u)+c$ for
any constant $c$ and hence $\pi_{N}$ descends to a map on $L(X)/\R.$
By the Arzelà-Ascoli theorem the latter space is a compact subspace
of the quotient Banach space $C^{0}(X)/\R$ (where $C^{0}(X)$ is
equipped with the usual $C_{0}-$norm $\left\Vert u\right\Vert :=\sup_{X}|u|).$
The existence of $u_{N}$ now follows from the Schauder fixed point
theorem applied to the continuous operator $\pi_{N}$ acting on the
compact convex space $L(X)/\R$ (the continuity is an immediate consequence
of the explicit expression \ref{eq:def of finte dim proj}). Note
that, by the mapping property \ref{eq:mapping property} $u_{N},$
is in fact contained in the subspace $\mathcal{P}_{N}(X)$ of $L(X).$
To conclude the proof of the proposition we observe that a direct
calculation reveals that the differential of $\mathcal{F}^{(N)}$
is given by the following formula: 
\[
d(\mathcal{F}^{(N)})_{|u}=e^{\beta_{N}(\pi_{N}(u)-u)}\mu_{0}
\]
Hence, if $\pi_{N}(u_{N})=u_{N}$ then $d(\mathcal{F}^{(N)})_{|u_{N}}=\mu_{0}$
which equivalently means that $u_{N}$ is a critical point of the
functional $\mathcal{F}_{\mu_{0}}^{(N)}.$ The maximization property
then follows directly from the fact that $\mathcal{F}_{\mu_{0}}^{(N)}$
is concave (since $\mathcal{F}^{(N)}$ is, as follows directly from
the concavity of log). Finally, since $\mathcal{F}_{\mu_{0}}^{(N)}$
is in fact\emph{ strictly} concave mod $\R$ (by the strict concavity
of log) the uniqueness of a critical point modulo additive constants
follows. \end{proof}
\begin{thm}
\label{thm:conv of balanced}Let $\mu_{0}$ be a probability measure
with compact support $X$ and $H^{(N)}$ an $N-$particle Hamiltonian
satisfying the assumptions in section \ref{sub:Setup-and-assumptions on hamilton}.
Let $\mu$ be another probability measure on $X$ with the same support
$X$ as $\mu_{0}.$ Then there exists (after perhaps passing to a
subsequence) a sequence of functions $u_{N}$ in $\mathcal{P}_{N}(X)$
which are balanced with respect to \emph{$(\mu,\beta_{N})$ and such
that $u_{N}\rightarrow u_{\mu}$ in $C^{0}(X),$ where $u_{\mu}$
is a potential for $\mu,$ i.e. 
\[
d\mathcal{F}_{|u_{\mu}}=\mu
\]
}\end{thm}
\begin{proof}
Let us first prove the result for $\mu=\mu_{0}.$ Let $u_{N}$ be
a sequence of functions which are balanced wrt \emph{$(\mu,\beta_{N}).$}
Since $\pi_{N}(u+c)=\pi_{N}(u)+c$ we may as well assume that $u_{N}$
is\emph{ normalized }in the sense that $u_{N}(x_{0})=0$ at some fixed
point $x_{0}.$ By the Arzelà-Ascoli theorem there exists, after perhaps
passing to a subsequence, a function $u_{\infty}\in L(X)$ such that
$u_{N}\rightarrow u_{\infty}$ in $C^{0}-$norm. Let us next show
the following 
\begin{equation}
\mbox{Claim:\,\ \ensuremath{u_{\infty}}\ensuremath{\mbox{ maximizes \ensuremath{\mathcal{F}_{\mu_{0}}}on \ensuremath{C^{0}(X)}}}}\label{eq:claim u infty maxim}
\end{equation}
First, by the previous proposition, $u_{N}$ maximizes $\mathcal{F}_{\mu}^{(N)}$
and hence if $u$ is a fixed element in $C^{0}(X)$ we get 
\[
\mathcal{F}_{\mu_{0}}(u):=\lim_{N\rightarrow\infty}\mathcal{F}_{\mu_{0}}^{(N)}(u)\leq\mathcal{F}_{\mu_{0}}^{(N)}(u_{N})
\]
Next, by construction $u_{N}\leq u_{\infty}+\delta_{N},$ where $\delta_{N}$
is a sequence of positive numbers tending to zero and hence 
\[
\mathcal{F}_{\mu_{0}}^{(N)}(u_{N})\leq\mathcal{F}_{\mu_{0}}^{(N)}(u_{\infty}+\delta_{N})=\mathcal{F}_{\mu_{0}}^{(N)}(u_{\infty}),
\]
 where, by definition, the rhs above converges to $\mathcal{F}_{\mu_{0}}(u_{\infty})$
as $N\rightarrow\infty.$ This proves the claim above. Since $\mathcal{F}$
is Gateaux differentiable it follows that the differential $d\mathcal{F}_{\mu_{0}}$
vanishes at $u_{\mu}:=u_{\infty},$ which translates to the equation
in the theorem. 
\end{proof}
Before turning to the proof of \ref{thm:exist of cont potentials}
we state the following corollary of the previous theorem (or rather
its proof), which gives yet another formula for the energy functional
$E$ defined by \ref{eq:def of E of mu as legendre transform}.
\begin{cor}
\label{cor:energy from limit of balanced}Let $\mu$ be a probability
measure with support $X$ and $\beta_{N}$ a sequence tending to infinity.
Denote by $u_{N}\in\mathcal{P}_{N}(X)$ a sequence of functions which
are balanced wrt $(\mu,\beta_{N}).$ Then 
\[
E(\mu)=\lim_{N\rightarrow\infty}\sup_{u\in\mathcal{P}_{N}(X)}\left(F_{\mu}^{(N)}(u)-\int u\mu\right)=\lim_{N\rightarrow\infty}F_{\mu}^{(N)}(u_{N})-\int u_{N}\mu
\]

\end{cor}

\subsection{Existence of continuous potentials }

We can now prove the existence of continuous potentials:
\begin{thm}
\label{thm:exist of cont potentials}Let $H^{(N)}$ be an $N-$particle
Hamiltonian satisfying the assumptions in section \ref{sub:Setup-and-assumptions on hamilton}
and denote by $\mathcal{F}(u)$ the corresponding limiting functional
(formula \ref{eq:def of f of u as limit of ham}). Then any probability
measure $\mu$ on $X$ admits a continuous potential $u_{\mu}.$ Equivalently,
$u_{\mu}$ solves the equation 
\[
d\mathcal{F}_{|u}=\mu
\]
Moreover, if $H^{(N)}$ admits a uniform Lipschitz constant $L,$
then so does $u_{\mu}.$ \end{thm}
\begin{proof}
In the case when $\mu$ has the same support as $\mu_{0}$ (i.e. the
set $X)$ the result follows immediately from the previous theorem
(since the limiting functional $\mathcal{F}$ only depends on the
support of $X,$ by Lemma \ref{lem:inf is Lbeta norm}). In the general
case we instead first apply Theorem \ref{thm:conv of balanced} to
$\mu_{\epsilon}:=(1-\epsilon)\mu+\epsilon\mu_{0}$ and obtain (normalized)
potentials $u_{\epsilon}$ of $\mu_{\epsilon}$ in $L(X).$ After
passing to a subsequence we may assume that $u_{\epsilon}\rightarrow u$
in $C^{0}-$norm. The proof is now concluded by noting that $u$ maximizes
the functional $\mathcal{F}_{\mu},$ as proved by a slight modification
of the proof of the claim appearing in the proof of Theorem \ref{thm:conv of balanced}.
\end{proof}

\subsection{The generalized mean field equations}

In this section we will show that the minimizer of the free energy
functional $F_{\beta}$ can be obtained from the solutions of the
equation \ref{eq:general mean field eq intro} appearing in the introduction
of the paper, which can be seen as a generalization of the real Monge-Ampère
equation \ref{eq:ma eq weak intro}. In fact, this part of the argument
was carried out in the more general setting of\emph{ complex} Monge-Ampère
equations in \cite{berm2}, which is analytically considerably more
involved.

We first define 

\[
\mathcal{G}_{\mu_{0},\beta}(u):=\mathcal{F}(u)-\frac{1}{\beta}\log\int e^{\beta u}\mu_{0}
\]
whose critical point equation is 
\begin{equation}
d\mathcal{F}_{|u}=\frac{e^{\beta u}\mu_{0}}{\int e^{\beta u}\mu_{0}}\label{eq:general of normalized exp m-a eq}
\end{equation}
By the Gateaux differentiability and concavity of $\mathcal{F}(u)$
we have that $u$ satisfies the previous equation iff $u$ maximizes
$\mathcal{G}_{\mu_{0},\beta}.$ The equation above is invariant under
the additive action of $\R$ on $C^{0}(X)$ and may hence be formulated
as an equation on $C^{0}(X)/\R$ which in turn is equivalent to the
following equation on $C^{0}(X):$
\begin{equation}
d\mathcal{F}_{|u}=e^{\beta u}\mu_{0}\label{eq:general of normalized exp m-a eq-2}
\end{equation}
The latter equation is the critical point equation for 
\[
\mathcal{\tilde{\mathcal{G}}}_{\mu_{0},\beta}(u):=\mathcal{F}(u)-\frac{1}{\beta}\int e^{\beta u}\mu_{0}
\]

\begin{thm}
\label{thm:existence of sol to general mean f}There exists a unique
continuous solution $u$ to the equation \ref{eq:general of normalized exp m-a eq-2}.
The corresponding probability measure $\mu:=d\mathcal{F}_{|u}$ is
the unique minimizer of the corresponding free energy functional $F_{\beta}$
on $\mathcal{M}_{1}(X).$ \end{thm}
\begin{proof}
The existence can be obtained using a variant of the balanced functions
used in the proof of Theorem \ref{thm:exist of cont potentials}.
One simply replaces the integration measure $\mu_{0}$ in the previous
definitions with the measure $\mu_{u}:=e^{\beta u}\mu_{0}$ (followed
by a suitable normalization). For example, one sets 
\begin{equation}
\pi_{N}(u)(x):=\frac{1}{\beta_{N}}\log\frac{(1-\beta/\beta_{N})}{Z_{N}[u]}\int_{y\in X^{N-1}}e^{-(\beta_{N}H^{(N)}(x,y)+u(y))}(\mu_{u}^{\otimes(N-1)}),\label{eq:def of finte dim proj-1}
\end{equation}
 where now $Z[u]:=\int_{X^{N}}e^{-(\beta_{N}H^{(N)}(x)+u(x))}\mu_{u}^{\otimes N}.$
Then a direct calculation gives 
\[
d\mathcal{F}_{|u}^{(N)}=e^{\beta_{N}(\pi_{N}(u)-u)}\mu_{u}
\]
where now $\mathcal{F}^{(N)}$ has been defined wrt the integration
measure $\mu_{u}$ and hence $\pi_{N}(u)=u$ iff $d\mathcal{F}_{|u}^{(N)}=\mu_{u}$
i.e. iff $u$ is a critical point of $\mathcal{G}_{\mu_{u},\beta}(u)$
or equivalently (by concavity) $u$ is a maximizer of $\mathcal{G}_{\mu_{u},\beta}(u).$
Moreover, this time we can directly apply the Banach fixed point theorem
on $C^{0}(X)$ to obtain a fixed point of $\pi_{N}.$ Indeed, $\pi_{N}$
defines a contraction mapping on $C^{0}(X)$ equipped with the sup-norm:
\[
\left\Vert \pi_{N}(u)-\pi_{N}(v)\right\Vert \leq(1-\beta/\beta_{N})\left\Vert u-v)\right\Vert .
\]
To see this note that, for any constant $c,$ $\pi_{N}(u+c)=\pi_{N}(u)+(1-\beta/\beta_{N})$
and if $u\leq w$ then $\pi_{N}(u)\leq\pi_{N}(w).$ Taking $c:=\left\Vert u-v)\right\Vert $
and $w:=u+c$ thus proves the upper bound in the contraction property
above and hence, by symmetry, the lower bound as well (interchanging
$u$ and $v).$ 

Finally, to prove that $\mu_{*}:=d\mathcal{F}_{|u}$ minimizes the
free energy functional $F_{\beta}$ on $\mathcal{M}_{1}(X)$ we take
$\mu$ in $\mathcal{M}_{1}(X)$ with finite entropy and consider the
affine segment $\mu_{t}:=(1-t)\mu_{*}+t\mu.$ By basic properties
of Legendre transforms (compare Lemma \ref{lem:maxim is differ}),
if $u_{\mu}$ is a potential for $\mu,$ then $-u_{\mu}$ is a sub-differential
for $E$ at $\mu$ and in particular 
\[
\frac{dE(\mu_{t})}{dt}_{t=0}\geq\int\left(-u\right)(\mu-\mu_{*})
\]
Recall also the well-known fact that the differential of the relative
entropy $D_{\mu_{0}}$ (at a point $\mu$ in the convex set where
it is finite) is represented by $\log(\mu/\mu_{0})$ (using for example
that $D_{\mu_{0}}$ is the Legendre transform of $u\mapsto\log\int e^{u}\mu_{0}$).
In particular, 

\[
\frac{dF_{\beta}(\mu_{t})}{dt}_{|t=0}\geq\int\left(-u+\frac{1}{\beta}\log(\mu_{*}/\mu_{0})\right)(\mu-\mu_{*})=0
\]
and hence by convexity $F_{\beta}(\mu)\geq F_{\beta}(\mu_{*}),$ as
desired. 
\end{proof}

\section{\label{sec:The-Large-Deviation}The Large Deviation Principle for
Gibbs measures}

\subsection{\label{sub:Setup:-the-Boltzmann-Gibbs}Setup: the Gibbs measure $\mu^{(N)}$
associated to the Hamiltonian $H^{(N)}$}

Let $X$ be topological space assumed to be compact (occasionally
we will also consider cases where $X$ is non-compact, in particular
in the Monge-Ampère setting) . A\emph{ random point process} with
$N$ particles is by definition a probability measure $\mu^{(N)}$
on the $N-$particle space $X^{N}$ which is symmetric, i.e. invariant
under permutations. Its \emph{one point correlation measure }$\mu_{1}^{(N)}$
(or \emph{first marginal}) is the probability measure on $X$ defined
as the push forward of $\mu_{1}^{(N)}$ to $X$ under the map $X^{N}\rightarrow X$
given by projection onto the first factor (or any factor, by symmetry):
\[
\mu_{1}^{(N)}:=\int_{X^{N-1}}\mu^{(N)}
\]
(similarly, the $j-$point correlation measure $\mu_{j}^{(N)}$ is
defined as the $j$ th marginal, i.e. the push forward to $X^{j}).$
In the following we will denote by $\mathcal{M}_{1}(Y)$ the space
of all probability measures on a space $Y$ and we will be particularly
concerned with the case when $Y=X^{N}.$ In the latter case we will
usually use the notation $\mu_{N}$ for (not necessarily symmetric)
elements of $\mathcal{M}_{1}(\mu_{N})$ and reserve the notation $\mu^{(N)}$
for specific Gibbs measures defined as below.

\subsubsection{The canonical Gibbs ensembles associated to the Hamiltonian $H^{(N)}$}

Fix a back-ground probability measure $\mu_{0}$ with support $X$
and let $H^{(N)}$ be a given $N-$particle Hamiltonian, i.e. a symmetric
continuous and bounded function on $X^{N}$ satisfying the assumptions
in section \ref{sub:Setup-and-assumptions on hamilton}. Also fixing
a positive number $\beta$ the corresponding\emph{ Gibbs measure}
is the symmetric probability measure on $X^{N}$ defined as 
\[
\mu_{\beta}^{(N)}:=e^{-\beta H^{(N)}}\mu_{0}^{\otimes N}/Z_{N},
\]
 where the normalizing constant 
\[
Z_{N,\beta}:=\int_{X^{N}}e^{-\beta H^{(N)}}\mu_{0}^{\otimes N}
\]
is called the\emph{ ($N-$particle) partition function. }Occasionally
we will simplify the notation by omitting the subscript $\beta.$

\subsection{\label{sub:Mean-entropy,-energy}Mean entropy, energy and free energy}

First we recall the general definition of the \emph{relative entropy
}(or the \emph{Kullback\textendash{}Leibler divergence}) of two measures
$\nu_{1}$ and $\nu_{2}$ on a space $Y:$ if $\nu_{1}$ is absolutely
continuous with respect to $\nu_{2},$ i.e. $\nu_{1}=f\nu_{2},$ one
defines 
\[
D(\nu_{1},\nu_{2}):=\int_{Y}\log(\nu_{1}/\nu_{2})\nu_{1}
\]
 and otherwise one declares that $D(\mu):=\infty.$ Note the sign
convention used: $D$ is minus the \emph{physical} entropy. 

Next, we define the\emph{ mean entropy} (relatively $\mu_{0}^{\otimes N})$
of a probability measure $\mu_{N}$ on $X^{N}$ (i.e. $\mu_{N}\in\mathcal{M}_{1}(X^{N}))$
as 
\[
D^{(N)}(\mu_{N}):=\frac{1}{N}D(\mu_{N},\mu_{0}^{\otimes N}).
\]
When $N=1$ we will simply write $D(\mu):=D^{(1)}(\mu)=D(\mu,\mu_{0}).$
On the other hand the \emph{mean energy} of $\mu_{N}$ is defined
as 
\[
E^{(N)}(\mu_{N}):=\frac{1}{N}\int_{X^{N}}H^{(N)}\mu_{N}
\]

Finally, the \emph{mean (Gibbs) free energy functional} on $\mathcal{M}_{1}(X^{N})$
is now defined as \emph{
\[
F^{(N)}:=E^{(N)}+\frac{1}{\beta}D^{(N)}
\]
}

Next, we will collect some basic general lemmas. First we have the
following simple special case of the well-known sub-additivity of
the entropy.
\begin{lem}
\label{lem: mean entropy}The following properties of the entropy
hold: 
\begin{itemize}
\item $D(\nu_{1},\nu_{2})\geq0$ with equality iff $\nu_{1}=\nu_{2}$
\item \textup{For a product measure on $X^{N}$ 
\[
D^{(N)}(\mu^{\otimes N})=D(\mu)
\]
}
\item More generally, \textup{
\[
D^{(N)}(\mu_{N})\geq D(\mu_{N,1}),
\]
 where $\mu_{N,1}$ is the corresponding first marginal (one point
correlation measure) on $X.$}
\end{itemize}
\end{lem}
The proof of the previous lemma uses only the (strict) concavity of
the function $t\mapsto\log t$ on $\R$ (see for example \cite{k}).
The latter (strict) concavity also immediately gives the following
\begin{lem}
\label{lem:(Gibbs-variational-principle).}(Gibbs variational principle).
Fix $\beta>0.$ Given a function $H^{(N)}$ on $X^{N}$ and a measure
$\mu_{0}$ on $X,$ the corresponding\emph{ free energy functional
$F^{(N)}$ on $\mathcal{M}_{1}(X^{N})$} attains its minimum value
on the corresponding Gibbs measure $\mu_{\beta}^{(N)}$ and only there.
More precisely, \emph{
\[
\inf_{\mathcal{M}_{1}(X^{N})}\beta F^{(N)}=\beta F^{(N)}(\mu^{(N)})=-\frac{1}{N}\log Z_{N}
\]
}\end{lem}
\begin{proof}
We recall the simple proof: since $\log(ab)=\log a+\log b,$ we have
\[
F^{(N)}(\mu_{N})=\frac{1}{\beta N}\int_{X^{N}}\log(\mu_{N}/e^{-\beta H^{(N)}}\mu_{0}^{\otimes N})\mu_{N}:=\frac{1}{\beta N}\int_{X^{N}}\log(\mu_{N}/\mu^{(N)})\mu_{N}-\frac{1}{N\beta}\log Z_{N}
\]
which proves the lemma using Jensen's inequality (i.e. the first point
in the previous lemma). 
\end{proof}
Note that the same argument applies if $\beta<0,$ since we can simply
replace $H^{N}$ with $-H^{N}$ in the previous argument (as long
as the corresponding partition function $Z_{N}$ is finite).

\subsection{Convergence of the one-point correlation measure towards the minimizer
of the free energy}

Given a background measure $\mu_{0}$ on $X$ and a Hamiltonian $H^{(N)}$
satisfying the assumptions in section \ref{sub:Setup-and-assumptions on hamilton}
we define, for any $\beta>0,$ the corresponding \emph{free energy
functional} as the following functional on the space $\mathcal{M}_{1}(X):$

\begin{equation}
F_{\beta}(\mu):=E(\mu)+D_{\mu_{0}}(\mu)/\beta,.\label{eq:def of free energy functional for beta}
\end{equation}
 where $E(\mu)$ is the thermodynamical energy defined by \ref{eq:def of E of mu as legendre transform}.
We start by proving a weak version of Theorem \ref{thm:interacting part intro},
stated in the introduction. 
\begin{thm}
\label{thm:conv of one-point cor measur}Let $H^{(N)}$ be an $N-$particle
Hamiltonian on $X^{N}$ satisfying the assumptions in section and
fix a positive number \textbf{$\beta.$ }Then the one-point correlation
measure $\mu_{1}^{(N)}$ of the corresponding Gibbs measure converges
weakly towards a probability measure $\mu_{*}$ on $X$ which is the
unique minimizer of the free energy functional $F_{\beta}$ on $\mathcal{M}_{1}(X).$
Moreover, \emph{
\[
\lim_{N\rightarrow\infty}-\frac{1}{\beta N}\log Z_{N,\beta}=\inf_{\mathcal{M}_{1}(X)}F_{\beta}
\]
}
\end{thm}
Take a measure $\mu\in\mathcal{M}_{1}(X).$ By the Gibbs variational
principle (i.e. the previous lemma) we have \emph{
\[
F^{(N)}(\mu^{(N)})\leq F^{(N)}(\mu^{\otimes N})=\frac{1}{N}\int H^{(N)}\mu^{\otimes N}+\frac{1}{\beta}D(\mu),
\]
}where we have used the first point in Lemma \ref{lem: mean entropy}
in the last equality. Hence, applying the upper bound in Theorem \ref{thm:existence of mean energ}
gives\emph{ 
\[
F^{(N)}(\mu^{(N)})\leq\inf_{\mu\in\mathcal{M}_{1}(X)}(E(\mu)+\frac{1}{\beta}D(\mu))
\]
 }To obtain a lower bound we apply the second point in Lemma \ref{lem: mean entropy}
to get \emph{
\[
\frac{1}{N}\int H^{(N)}\mu^{(N)}+\frac{1}{\beta}D(\mu_{1}^{(N)})\leq F^{(N)}(\mu^{(N)})
\]
}Next, we fix $u\in C^{0}(X)$ and rewrite the first term in the lhs
above as follows: \emph{
\[
\frac{1}{N}\int H^{(N)}\mu^{(N)}=\int\frac{1}{N}(H^{(N)}+u)\mu^{(N)}-\int_{X}u\mu_{1}^{(N)}
\]
}and hence replacing the first integral in the rhs with its infimum
gives \emph{
\[
\inf_{X^{N}}\frac{1}{N}(H^{(N)}+u)-\int_{X}u\mu_{1}^{(N)}\leq\frac{1}{N}\int H^{(N)}\mu^{(N)},
\]
 }which by definition means that, if $\mu_{*}$ denotes a weak limit
point of the sequence $\mu_{1}^{(N)},$ then \emph{
\[
\mathcal{F}(u)-\int_{X}u\mu_{*}\leq\lim\inf_{N\rightarrow\infty}\frac{1}{N}\int H^{(N)}\mu^{(N)},
\]
}Taking the sup over all $u\in C^{0}(X)$ thus gives \emph{
\[
E(\mu_{*})\leq\lim\inf_{N\rightarrow\infty}\frac{1}{N}\int H^{(N)}\mu^{(N)}.
\]
 }Next, since $D$ is lower-semi continuous we have \emph{
\[
D(\mu_{*})\leq\lim\inf_{N\rightarrow\infty}D(\mu_{1}^{(N)})
\]
}All in all, this means that \emph{
\[
F(\mu_{*})\leq\lim\inf_{N\rightarrow\infty}F^{(N)}(\mu^{\otimes N})\leq\lim\sup_{N\rightarrow\infty}F^{(N)}(\mu^{\otimes N})\leq\inf_{\mu\in\mathcal{M}_{1}(X)}F_{\beta}(\mu)
\]
}But then it must be that all the inequalities in the previous line
are actually equalities. In particular, $\mu_{*}$ is a minimizer
of $F.$ Finally, by construction $E$ is convex on $\mathcal{M}_{1}(X)$
and since $D(\mu)$ is strictly convex on the convex subset of $\mathcal{M}_{1}(X)$
where it is finite it follows that $F_{\beta}$ is also strictly convex
on the subset where it is finite. As a consequence $F$ has a unique
minimizer which thus must coincide with the limit point $\mu_{*}.$

\subsection{\label{sub:Large-deviations-(proof}Proof of the Large Deviation
Principle (Theorem \ref{thm:interacting part intro})}

Let us recall the general definition of a LDP due to Donsker and Varadhan
(see for example the book \cite{de-ze}):
\begin{defn}
\label{def:large dev}Let $\mathcal{M}$ be a Polish space, i.e. a
complete separable metric space.

$(i)$ A function $I:\mathcal{\, M}\rightarrow[0,\infty]$ is a \emph{rate
function} iff it is lower semi-continuous. It is a \emph{good} \emph{rate
function} if it is also proper.

$(ii)$ A sequence $\Gamma_{k}$ of measures on $\mathcal{M}$ satisfies
a \emph{large deviation principle} with \emph{speed} $r_{k}$ and
\emph{rate function} $I$ if

\[
\limsup_{k\rightarrow\infty}\frac{1}{r_{k}}\log\Gamma_{k}(\mathcal{F})\leq-\inf_{\mu\in\mathcal{F}}I(\mu)
\]
 for any closed subset $\mathcal{F}$ of $\mathcal{M}$ and 
\[
\liminf_{k\rightarrow\infty}\frac{1}{r_{k}}\log\Gamma_{k}(\mathcal{G})\geq-\inf_{\mu\in G}I(\mu)
\]
 for any open subset $\mathcal{G}$ of $\mathcal{M}.$
\end{defn}

\subsubsection{Functional analytic framework and Legendre-Fenchel transforms}

Let $X$ be a Polish space and denote by $C_{b}(X)$ the space of
all bounded continuous functions on $X$ and by $\mathcal{M}(X)$
the space of all signed finite Borel $\mu$ measures on $X.$ We will
write the corresponding integration pairing as 
\[
\left\langle u,\mu\right\rangle :=\int_{X}u\mu
\]
We equip $\mathcal{M}(X)$ with the weak topology generated by $C_{b}(X)$.
Then we may identify $C_{b}(X)$ with the topological dual $\mathcal{M}(X)^{*}$
of $\mathcal{M}(X),$ i.e. with the space of all linear continuous
functions on $\mathcal{M}(X)$ \cite{de-ze}. We will be mainly concerned
with the subspace $\mathcal{M}_{1}(X)$ of all probability measures
on $X$ which is a convex subset of $\mathcal{M}(X)$ (and compact
iff $X$ is compact). This latter space is a locally convex topological
vector space. As such it admits a good duality theory (see section
4.5.2 in \cite{de-ze}): given a functional $\Lambda$ on the vector
space $C_{b}(X)$ its \emph{Legendre(-Fenchel) transform} is the following
functional $\Lambda^{*}$ on $\mathcal{M}(X):$ 
\[
\Lambda^{*}(\mu):=\sup_{u\in C_{b}(X)}(\Lambda(u)-\left\langle u,\mu\right\rangle )
\]
 Conversely, if $H$ is a functional on the vector space $\mathcal{M}(X)$
we let 
\[
H^{*}(u):=\inf_{\mu\in\mathcal{M}(X)}(H(\mu)+\left\langle u,\mu\right\rangle )
\]
 Note that we are using rather non-standard sign conventions. In particular,
$\Lambda^{*}(\mu)$ is always \emph{convex }and \emph{lower semi-continuous
(lsc),} while $H^{*}(u)$ is \emph{concave} and \emph{upper-semicontinuos
(usc). }As a well-known consequence of the Hahn-Banach separation
theorem we have the following fundamental duality relation (Lemma
4.5.8 in \cite{de-ze}): 
\begin{equation}
\Lambda=(\Lambda^{*})^{*}\label{eq:dualiy}
\end{equation}
 iff $\Lambda$ is concave and usc. We also recall the following standard
\begin{lem}
\label{lem:maxim is differ}Assume that $\Lambda$ is a functional
on $C_{b}(X)$ which is finite, lsc, concave and Gateaux differentiable
(i.e differentiable along lines). Then, for a fixed $u\in C_{b}(X)$
the differential $d\Lambda_{|u}$ is the unique minimizer of the following
functional on $\mathcal{M}(X):$ 
\begin{equation}
\mu\mapsto\Lambda^{*}(\mu)+\left\langle u,\mu\right\rangle \label{eq:functional in lemma max diff}
\end{equation}
 (and the minimum value equals $\Lambda(u)).$ Conversely, if the
latter functional admits a unique minimizer $\mu_{u}$ for any $u\in C_{b}(X)$,
then $\Lambda(u)$ is Gateaux differentiable and $d\Lambda_{|u}=\mu_{u}.$
\end{lem}
Recall that, in general, if $\Gamma$ is a probability measure on
a topological vector space $\mathcal{M}$ then its\emph{ Laplace transform
}is the following functional defined on the topological dual $\mathcal{M}^{*}$
of $\mathcal{M}$: 
\[
\widehat{\Gamma}[u]:=\int_{\mu\in\mathcal{M}}\Gamma e^{-\left\langle u,\mu\right\rangle }
\]
(assuming that the integral is finite). 

We will use the following abstract form of the Gärtner-Ellis theorem
(see \cite{de-ze} , Cor 4.6.14, p. 148) and references therein about
the different versions of this theorem) which can be seen as an infinite
dimensional version of the method of stationary phase, i.e. the Laplace
method.
\begin{thm}
\label{thm:(Abstract-G=0000E4rtner-Ellis-Theorem}(abstract Gärtner-Ellis
Theorem ). Let $\mathcal{M}$ be a locally convex Hausdorff topological
vector space and $\Gamma_{k}$ a sequence of Borel measures on $\mathcal{M}$
which is exponentially tight with respect a sequence $r_{k}$ of positive
numbers and such that the Laplace transforms $\widehat{\Gamma}_{k},$
seen as functionals on the dual $\mathcal{M}^{*},$ satisfy 
\[
-\frac{1}{r_{k}}\log\widehat{\Gamma}_{k}[r_{k}u]\rightarrow\Lambda[u]
\]
 for any $u$ in $\mathcal{M}^{*}$ where the functional $\Lambda$
is Gateau differentiable on $\mathcal{M}^{*}.$ Then $\Gamma_{k}$
satisfies a LDP with speed $r_{k}$ and with a rate functional $H:=\Lambda^{*}$
on $\mathcal{M},$ i.e. the rate function  $H$ is the Legendre-Fenchel
transform of $\Lambda.$
\end{thm}
The definition of exponential tightness will be recalled in the proof
of Theorem \ref{thm:main for perman intro} (in case $X$ is compact
this condition is automatically satisfied).

\subsubsection{End of proof of Theorem \ref{thm:interacting part intro}}

Given $u\in C^{0}(X)$ we let $f_{N}[u]:=-\frac{1}{N\beta}\log\int\mu^{(N)}e^{-N\beta\left\langle u,\delta_{N}\right\rangle },$
i.e. the scaled logarithm of the Laplace transform at $u$ of the
probability measure $\Gamma:=(\delta_{N})_{*}\mu^{(N)}$ on $\mathcal{M}(X).$
It may also be written as 
\[
\Lambda_{N}[u]=-\frac{1}{N\beta}\log\frac{Z_{N,\beta}[u]}{Z_{N,\beta}[0]},
\]
 where 
\begin{equation}
Z_{N,\beta}[u]=\int_{X^{N}}e^{-\beta(H+u)}\mu_{0}^{\otimes N}\label{eq:perturbed parition function}
\end{equation}
By the previous theorem applied to $H_{u}^{(N)}:=H^{(N)}+u$ we have
\[
\Lambda_{N}[u]\rightarrow\Lambda(u):=\inf_{\mu\in\mathcal{M}_{1}(X)}(E(\mu)+\int_{X}u\mu+\frac{1}{\beta}D(\mu))-C:=\inf_{\mu\in\mathcal{M}_{1}(X)}(F_{\beta}(\mu)+\int_{X}u\mu)-C:
\]
Next, we observe that the infimum in the rhs above is up to a harmless
additive constant, by definition, the Legendre transform of the extended
free energy functional $\mu\mapsto F_{\beta}(\mu),$ defined by the
expression \ref{eq:def of free energy functional for beta} on $\mathcal{M}_{1}(X)$
and set to be equal to $\infty$ on the complement of $\mathcal{M}_{1}(X)$
in $\mathcal{M}(X).$ But, as explained in the proof of the previous
theorem, the corresponding functional is strictly convex on the subspace
where it is finite. But then it follows from Lemma \ref{lem:maxim is differ}
(and the duality relation $F=\Lambda^{*})$ that $\Lambda$ is Gateaux
differentiable and hence we can apply the Gärtner-Ellis theorem to
conclude.

\subsection{\label{sub:The-ambient-point}The ambient point of view}

To better see the connection to the Monge-Ampère setting, to be considered
in the following section, we consider the following general ``ambient
setting''. Start with a (possible non-compact) Riemannian manifold
$(Y,g)$ - to be referred to as the ``ambient space'' (which in
the Monge-Ampère will be equal to $\R^{n}$ equipped with the Euclidean
metric) and a Hamiltonian $H^{(N)}$ on $Y$ satisfying the Lipschitz
assumption in section \ref{sub:Setup-and-assumptions on hamilton},
with a Lipschitz constant $L,$ defined with respect to the Riemannian
metric $g.$ Fix a reference element $\phi_{0}$ in the corresponding
space $\mathcal{P}(Y)$ (to be referred to as the space of ``ambient
potentials'') and assume that 
\begin{itemize}
\item For any function $u\in C_{b}(Y)$ the following limit exists 
\[
\mathcal{F}(u):=\lim_{N\rightarrow\infty}\inf_{Y^{N}}\frac{1}{N}(H^{(N)}+\phi_{0}+u),
\]
 
\item The functional $\mathcal{F}(u)$ is Gateaux differentiable.
\end{itemize}
Also, on the subspace $\mathcal{P}_{+}(Y)$ of all functions $\phi$
in $\mathcal{P}(Y)$ such $\phi-\phi_{0}$ is bounded we can then
define an operator 
\[
M:\,\,\,\mathcal{P}(Y)\rightarrow\mathcal{M}_{1}(Y),\,\,\,\,\, M(\phi)=d\mathcal{F}_{|\phi-\phi_{0}}
\]
(which can be seen as a generalized Monge-Ampère operator). 

Now, given a sequence $\beta_{N}$ of positive numbers any choice
of probability measure $\mu_{0}$ on $Y$ with compact support $X$
induces a random point processes on $X$ defined by the Gibbs measure
of the restricted Hamiltonian (compare section \ref{sub:Setup:-the-Boltzmann-Gibbs}).
Note that, by construction, there is a natural extension map from
$\mathcal{P}_{N}(X)$ to $\mathcal{P}_{N}(Y).$ Under suitable assumptions
one can then show that the sequence of balanced metrics appearing
in the proof of Theorem \ref{thm:existence of sol to general mean f}
 converge to a function $\phi$ satisfying the following generalized
mean field equations on the ambient space $Y:$

\begin{equation}
M(\phi)=\frac{e^{\beta_{N}\phi}\mu_{0}}{\int e^{\beta_{N}\phi}\mu_{0}}\label{eq:general of normalized exp m-a eq-1}
\end{equation}
We will not develop this general theory further, as the corresponding
results will be obtained directly in the particular Monge-Ampère setting
considered in section \ref{sub:Existence-and-extremal}. A remarkable
feature of the latter setting is that the corresponding operator $M(\phi)$
is a\emph{ local} operator, since it coincides with the Monge-Ampère
measure which (up to regularization of $\phi)$ is a differential
operator.

\subsection{\label{sub:An-alternative-proof}An alternative proof of the LDP
for equicontinuous Hamiltonians on compact spaces}

In this section we assume again that $X$ is a compact metric space
equipped with a probability measure $\mu_{0}$ whose support is equal
to $X.$ Let us start by recalling a large deviation result from \cite{e-h-t},
which in the present setting may be formulated as follows.
\begin{thm}
\label{thm:-ellis et al}\cite{e-h-t} Let $X$ be a compact metric
space, $H_{N}$ a sequence of bounded continuous symmetric functions
on $X^{N}$ and $U$ a bounded continuous function on $\mathcal{M}_{1}(X)$
such that $\sup_{X^{N}}\left|H^{(N)}/N-U(\delta_{N})\right|\rightarrow0,$
as $N\rightarrow\infty.$ Then, for any real number $\beta$ the measure
$(\delta_{N})_{*}(e^{-\beta H^{(N)}}\mu_{0}^{\otimes N})$ on $\mathcal{M}_{1}(X)$
satisfies a LDP with speed $N$ and rate functional $\beta E+D_{\mu_{0}}.$ 
\end{thm}
As will be explained below the proof can be reduced to the special
case when $H^{(N)}=0,$ which is the content of Sanov's classical
result (this is slightly different than the reduction to Sanov's theorem
in \cite{e-h-t}, which uses a Laplace principle). Using Theorem \ref{thm:-ellis et al}
we will in this section give a more direct proof of Theorem \ref{thm:interacting part intro}
stated in the introduction. In fact, we will obtain the following
more general form of the latter theorem \ref{thm:interacting part intro}.
\begin{thm}
\label{thm:general form of interacting particle in compact case}Let
$X$ be a compact metric space and $H^{(N)}/N$ a sequence of symmetric
functions on $X^{N}$ which is uniformly bounded and equicontinuous
(in the sense of section \ref{sub:Setup-and-assumptions on hamilton}).
Let $\beta_{N}$ be a sequence of positive number tending to infinity
and assume that $(\delta_{N})_{*}(e^{-\beta_{N}H^{(N)}}\mu_{0}^{\otimes N})$
satisfies an LDP with rate functional $E(\mu)$ and speed $\beta_{N}N.$
Then, for $\beta$ any fixed (possibly negative) number $(\delta_{N})_{*}(e^{-\beta H^{(N)}}\mu_{0}^{\otimes N})$
satisfies a LDP with speed $\beta N$ and rate functional $E+D_{\mu_{0}}/\beta.$ 
\end{thm}
The previous result combined with the Gärtner-Ellis theorem immediately
gives Theorem \ref{thm:interacting part intro} stated in the introduction.
Note that in the setting of the latter theorem the functional $E$
is automatically convex. In order to reduce the proof of the previous
theorem to Theorem \ref{thm:-ellis et al} we will invoke the following 
\begin{lem}
There exists a continuous and bounded function $U$ on $\mathcal{M}_{1}(X)$
such that $\lim_{j\rightarrow\infty}\sup_{X^{N}}|H_{N_{j}}/N_{j}-(\delta_{N_{j}})^{*}U|=0.$
In fact, $U$ can be taken to have the same modulus of continuity
as the sequence $H^{(N)}/N.$\end{lem}
\begin{proof}
First we may, using the map \ref{eq:def of empricical measure} $\delta_{N}$
embed $X^{(N)}$ isometrically in $\mathcal{M}_{1}(X),$ equipped
with the Wasserstein 1-metric. Then $H^{(N)}/N$ extends to a continuous
function $U_{N}$ defined on all of $\mathcal{M}_{1}(X)$ preserving
the modulus of continuity. Accepting this for the moment the existence
of $U$ follows immediately from the Arzelà-Ascoli theorem applied
to the sequence $U_{N}$ on the compact space $\mathcal{M}_{1}(X).$
As for the extension property if follows from general considerations:
let $\mathcal{K}$ be a subspace of a metric space$(\mathcal{M},d)$
and $u_{N}$ a sequence of functions on $\mathcal{K}$ which we for
simplicity assume is uniformly Lipschitz continuous (the general case
is proved in a similar manner). Setting $U_{N}(y):=\inf_{x\in\mathcal{K}}(u(x)+d(x,y))$
then gives the desired extension.
\end{proof}
To prove Theorem \ref{thm:general form of interacting particle in compact case}
we start by showing that the following
\[
\mbox{claim:\,\,}E(\mu)=U(\mu)
\]
holds, for any $\mu$ such that $D_{\mu_{0}}(\mu)\neq0$ and in particular
the convergence in the previous lemma must hold for the full sequence
indexed by $N.$ To prove the claim first note that the assumed LDP
implies, by a general result for LDPs (see Theorem 4.1.18 in \cite{de-ze}),
that, setting $\Gamma_{N,\beta_{N}}:=(\delta_{N})_{*}(e^{-\beta_{N}H^{(N)}}\mu_{0}^{\otimes N}),$
\[
\lim_{\delta\rightarrow0}\liminf_{N\rightarrow\infty}\frac{1}{\beta_{N}N}\log\int_{B_{\delta}(\mu)}\Gamma_{N,\beta_{N}}=\lim_{\delta\rightarrow0}\limsup_{N\rightarrow\infty}\frac{1}{\beta_{N}N}\log\int_{B_{\delta}(\mu)}\Gamma_{N,\beta_{N}}=E(\mu),
\]
 where $B_{\delta}(\mu)$ is the ball of radius $\delta$ centered
at $\mu,$ defined with respect to any metric compatible with the
weak topology on the space $\mathcal{M}_{1}(X)$ and since $X$ is
assumed compact we can take the metric to the Wasserstein 1-metric.
Now, by Sanov's theorem 
\[
\lim_{\delta\rightarrow0}\liminf_{N\rightarrow\infty}\frac{1}{N}\log\int_{B_{\delta}(\mu)}(\delta_{N})_{*}\mu_{0}^{\otimes N}\geq-D_{\mu_{0}}(\mu)
\]
Hence, if $D_{\mu_{0}}(\mu)<\infty$ the previous lemma gives (again
using that the map \ref{eq:def of empricical measure} defined by
$\delta_{N}$ is an isometry) that 
\[
U(\mu)+\frac{1}{\beta_{N}}O(1)=E(\mu),
\]
 where $O(1)$ denotes a bounded term and since $\beta_{N}\rightarrow\infty$
this proves the claim above. We can now apply Theorem \ref{thm:-ellis et al}
to deduce the LDP with respect to a fixed $\beta$ with the rate functional
$\beta U+D$ and speed $N,$ which according to the previous claim
coincides with $\beta E+D,$ as desired (using that $D_{\mu_{0}}(\mu)=\infty$
iff $F_{\beta}(\mu)=\infty).$ Alternatively, repeating the arguments
above with $\beta_{N}$ replaced with $\beta$ immediately gives,
using the previous lemma, that 
\[
\lim_{\delta\rightarrow0}\liminf\frac{1}{\beta N}\log_{N\rightarrow\infty}\int_{B_{\delta}(\mu)}\Gamma_{N,\beta}=\lim_{\delta\rightarrow0}\limsup\frac{1}{\beta N}\log_{N\rightarrow\infty}\int_{B_{\delta}(\mu)}=U+\frac{1}{\beta}D,
\]
 which coincides with $E+\frac{1}{\beta}D.$ Finally, by Theorem 4.1.11
and Lemma 1.2.18 in \cite{de-ze} the desired LDP holds follows.

\section{\label{sec:Permanental-point-processes,}Permanental point processes,
the Monge-Ampère and optimal transport }

\subsection{\label{sub:Setup function spaces in convx}Setup}

Consider $\R^{n}$ equipped with the Euclidean scalar product that
we will denote by a dot. The corresponding Euclidean coordinates will
be denoted $x=(x_{1},...,x_{n}).$ We will identify the dual linear
space with $\R^{n}$ equipped with the coordinates $p=(p_{1},...,p_{N})$
so that the corresponding duality pairing may be written as $\left\langle x,p\right\rangle =x\cdot p.$
We fix once and for all a convex body $P$ in $\R^{n}.$ Without loss
of generality we may assume that $0$ is contained in the interior
of $P$ and that $P$ has unit Euclidean volume. The support function
of $P$ will be denoted by $\phi_{P},$ i.e. $\phi_{P}(x):=\sup_{p\in P}\left\langle x,p\right\rangle .$
The set of all bounded continuous function $u$ on $\R^{n}$ will
be denoted by $C_{b}(\R^{n}).$

Given a convex function $\phi$ on $\R^{n}$ we will denote by $MA(\phi)$
its Monge-Ampère measure (in the sense of Alexandrov), i.e. if $E$
is a Borel set then $(MA(\phi)(E)$ is defined as the Lebesgue measure
of the image of $E$ under the sub-gradient $\nabla\phi$ of $\phi$
\cite{gu} (viewed as a multivalued map from $\R^{n}$ to $\R^{n})$
\cite{gu}. Following \cite{ber-ber} we will denote by $\mathcal{P}(\R^{n})$
the space of all convex functions $\phi$ on $\R^{n}$ such that $\phi-\phi_{P}$
is bounded from above and by $\mathcal{P}_{+}(\R^{n})$ the subspace
of all $\phi$ such that $\phi-\phi_{P}$ is bounded. 

Given a (possible non-convex) upper-semi continuous (usc) function
$\phi$ will write $\phi^{*}$ for its Legendre transform, i.e. 
\begin{equation}
\phi^{*}(p):=\sup_{x\in\R^{n}}\left\langle x,p\right\rangle -\phi(x)\label{eq:def of leg trans of functions}
\end{equation}
 By basic properties of the Legendre transform we have that, if $\phi$
is in $\mathcal{P}(\R^{n}),$ then the subgradient image $(\nabla\phi)(\R^{n})$
is contained in $P$ and $\phi^{*}$ is equal to $\infty$ on the
complement of $P.$ Similarly, under the Legendre transform the space
$\mathcal{P}_{+}(\R^{n})$ corresponds to the space of all bounded
convex functions on $P$ (see \cite{ber-ber} and references therein).
Finally, following \cite{ber-ber} we will say that $\phi$ $\mathcal{P}(\R^{n})$
has\emph{ finite energy} if $\phi^{*}$ is integrable on $P$, i.e.
if 
\[
\mathcal{E}(\phi):=-\int_{P}\phi^{*}(p)dp>-\infty
\]
and we will denote the space of all finite energy convex functions
by $\mathcal{E}_{P}^{1}(\R^{n}).$ We thus we have the following (strict)
inclusions: 
\[
\mathcal{P}_{+}(\R^{n})\subset\mathcal{E}_{P}^{1}(\R^{n})\subset\mathcal{P}(\R^{n})
\]
Let us finally recall the following basic compactness property of
the space $\mathcal{P}(\R^{n})$ which is an immediate consequence
of the Arzelà-Ascoli theorem (compare \cite{ber-ber})
\begin{prop}
\textup{\label{prop:compactness of space of convex}Let $\phi_{j}$
be a sequence of }\textup{\emph{normalized}}\textup{ functions in
$\mathcal{P}(\R^{n}),$ which by definition means that $\sup_{\R^{n}}(\phi_{j}-\phi_{P})=0$
or equivalently that $\phi_{j}(0)=0.$ Then, perhaps after passing
to a subsequence, $\phi_{j}$ converges locally uniformly to a normalized
function $\phi$ in $\mathcal{P}(\R^{n}).$}
\end{prop}

\subsubsection{\label{sub:Weigted-sets-and}Weighted sets and measures}

By definition a\emph{ weighted set} $(X,\phi_{0})$ consists of a
closed set $X$ in $\R^{n}$ and a\emph{ weight function} $\phi_{0}$
on $X,$ i.e. a continuous function $\phi_{0}$ on $X$ such that
$\phi-\phi_{P}\rightarrow\infty$ as $|x|\rightarrow\infty$ in $X$
(in particular if $X$ is compact then latter growth condition is
vacuous). Occasionally we will identify $\phi_{0}$ with a function
on $\R^{n}$ by letting the extension be identically equal to $\infty$
on the complement of $X.$ Given a measure $\mu_{0}$ on $X$ we say
that $(\mu_{0},\phi_{0})$ is a\emph{ weighted measure} if $\phi_{0}$
is a weight function on the support $X$ of $\mu_{0}$ and $e^{\beta(\phi_{P}-\phi_{0})}\mu_{0}$
has finite total mass for any positive number $\beta.$

Given a closed set $X$ we define the corresponding \emph{projection
operator} $\Pi_{X}$ from the space of weights on $X$ to the space
$\mathcal{P}_{+}(\R^{n})$ by the following convex envelope 
\[
\Pi_{X}(\phi_{0})(x)=\sup_{\phi\in\mathcal{P}(\R^{n})}\{\phi(x):\,\,\,\,\phi\leq\phi_{0}\,\,\mbox{on\,\ensuremath{X\}}}
\]
(if $\phi$ is continuous then $\Pi_{X}(\phi)$ is in $\mathcal{P}(\R^{n})$
and if $\phi$ is a weight function then $\Pi_{X}(\phi)$ is in $\mathcal{P}_{+}(\R^{n}),$
since $\phi_{P}-C$ is a candidate for the sup above, if $C$ is sufficiently
large). In particular, if $(X,\phi_{0})$ is a weighted set we will
occasionally write $\phi_{e}:=\Pi_{X}(\phi_{0})$ and $\mu_{e}:=MA(\Pi_{X}(\phi_{0})).$
From the growth assumption on $\phi_{0}$ it follows that the \emph{incidence
set} 
\[
D_{0}:=\{\Pi_{X}(\phi_{0})=\phi_{0}\}
\]
 is compact and, since, by general properties of free convex envelopes,
$MA(\Pi_{X}(\phi))$ is always contained in $D_{0}$ we conclude that
$\mu_{e}$ has compact support contained in $X.$ We also note that
``the Legendre transform doesn't see the projection $\Pi_{X}$''
in the following sense 
\begin{equation}
(\Pi_{X}(\phi_{0}))^{*}=\phi_{0}^{*}\,\,\,\,\mbox{on\,\ensuremath{P}, }\label{eq:lef does not see proj}
\end{equation}
(which is a special case of formula \ref{eq:sup same as projec} in
Lemma \ref{lem:b-m for convex etc} below). We note that $\Pi_{X}(\phi_{0}))^{*}$
is bounded since $\Pi_{X}(\phi_{0})$ is in $\mathcal{P}_{+}(\R^{n}).$

\subsubsection{The generalized permanental point processes }

Set 
\[
\mbox{Per}(x_{1},...,x_{N_{k}}):=\mbox{Per}(e^{x_{i}\cdot p_{j}}),
\]
 where $p_{j}$ ranges over the $N$ lattice points in $kP$ (compare
the notation in the introduction of the paper). For a given sequence
$\beta_{N}$ such that $\lim_{N\rightarrow\infty}\beta_{N}=\beta\in]0,\infty]$
we let $\mu_{\beta_{N}}^{(N)}$ be the probability measure on $X^{N}$
defined by 
\[
\mu_{\beta_{N}}^{(N)}:=\frac{1}{Z_{N,\beta_{N}}}\mbox{Per}(x_{1},...,x_{N_{k}})^{\beta_{N}/k}e^{-\beta_{N}(\phi(x_{1})+\cdots+\phi(x_{N}))}\mu_{0}^{\otimes N},
\]
 where $Z_{N,\beta_{N}}$ is the normalizing constant (which is finite
by the growth assumption on $\phi_{0};$ see below). Setting 
\begin{equation}
H_{\phi_{0}}^{(N)}:=-\frac{1}{k}\mbox{Per}(x_{1},...,x_{N_{k}})+\phi(x_{1})+\cdots+\phi(x_{N})\label{eq:def of hamilt perfurbed by weight}
\end{equation}
the probability measure $\mu_{\beta_{N}}^{(N)}$ thus becomes the
Gibbs measure, at inverse temperature $\beta_{N},$ determined by
the $N-$particle Hamiltonian $H_{\phi_{0}}^{(N)}.$

\subsection{Large $N$ asymptotics and variational principles}

Since the logarithm of any convex combination of functions of the
form $e^{x\cdot p}$ for $p\in P$ is in the space $\mathcal{P}(\R^{n})$
it follows that for $(x_{2},...x_{N})\in X^{N}$ fixed $\psi(x):=-H^{(N)}(x,x_{2}...,x_{N})$
defines an element in $\mathcal{P}(\R^{n})$ and similarly when the
other variables are fixed. In particular, the partial gradients satisfy
\begin{equation}
\nabla_{x_{i}}H^{(N)}\in-P\label{eq:gradient of hamiltonian is in p}
\end{equation}
 for any $N$ and any $x_{i}\in X.$ Since $P$ is assumed bounded
$H^{(N)}/N$ is uniformly Lipschitz continuous \ref{sub:Setup-and-assumptions on hamilton}
and when $X$ is compact $H_{\phi_{0}}^{(N)}/N$ is thus equicontinuous
in the sense of section \ref{sub:Setup-and-assumptions on hamilton}.
In order to be able to handle the non-compact case we recall that
the support of $\mu_{e}:=MA(\Pi_{X}\phi_{0})$ has compact support
and it is thus contained in a large ball $B_{R}.$ 
\begin{lem}
\label{lem:b-m for convex etc}Let $(\mu_{0},\phi_{0})$ be a weighted
measure and $X$ the support of $\mu_{0}.$ Then 
\[
\Pi_{X\cap B_{R}}\phi_{0}=\Pi_{X}\phi_{0}
\]
and for any $\phi_{k}$ in $\mathcal{P}(\R^{n})$ we have 
\begin{equation}
\sup_{X}e^{k(\phi_{k}-\phi_{0})}=\sup_{\R^{n}}e^{k(\phi_{k}-\Pi_{X}\phi_{0})}=\sup_{D_{\phi_{0}}}e^{k(\phi_{k}-\phi_{0})}\label{eq:sup same as projec}
\end{equation}
Moreover, for any $\epsilon>0,$ there is $C_{\epsilon}>0$ (independent
of $\phi_{k})$ such that 

\[
\sup_{X}e^{k(\phi_{k}-\phi_{0})}\leq C_{\epsilon}e^{k(\Pi_{X}\phi_{0}-\phi_{0})}e^{\epsilon k}\int_{X\cap B_{R}}e^{k(\phi_{k}-\phi_{0})}\mu_{0}
\]
 \end{lem}
\begin{proof}
By definition $\Pi_{X\cap B_{R}}\varphi\geq\Pi_{X}\varphi$ and $\Pi_{X\cap B_{R}}\varphi\leq\varphi$
on the support of $MA(\Pi_{X}\varphi).$ Since the latter set is contained
in the incidence set $D$ (see above) this means that $\Pi_{X\cap B_{R}}\varphi\leq\Pi_{X}\varphi$
a.e. wrt $MA(\Pi_{X}\varphi)$ and hence the inequality holds everywhere
according to the domination principle for $MA$ (see the appendix).
This shows that $\Pi_{X\cap B_{R}}\varphi=\Pi_{X}\varphi.$ Next,
note that the first equality in \ref{eq:sup same as projec} follows
directly from the extremal definition of $\Pi_{X}$ and the second
one follows from the domination principle, as in the previous equality.
To prove the inequality in the lemma we first observe that, when $X$
is\emph{ compac}t we have 
\[
\sup_{X}e^{k(\phi_{k}-\phi_{0})}\leq C_{\epsilon}e^{\epsilon k}\int_{X}e^{k(\phi_{k}-\phi_{0})}\mu_{0},
\]
using that $\phi_{k}$ is uniformly continuous on $X$ with a constant
of continuity only depending on $P$ (in fact, since $\nabla\phi_{k}$
is in $P$ we even have a Lipschitz constant only depending on $P).$
Indeed, given $\epsilon>0$ and $\phi_{k}$ we simply estimate $\int_{X}e^{k(\phi_{k}-\phi_{0})}\mu_{0}$
from below by the integral over a small poly-disc $\Delta_{\delta}(x_{k})$
of radius $\delta$ centered at the point $x_{k}\in X$ where the
sup of $\phi_{k}$ is attained. By the Lipschitz property $\delta$
can be chosen so that $e^{k(\phi_{k}-\phi_{0})}\geq e^{k((\phi_{k}-\phi_{0})(x_{k})-\epsilon)}$and
hence the desired estimate holds with the constant $C_{\delta}:=\inf_{x\in X}\mu_{0}(\Delta_{\delta}(x)),$
which is strictly positive for any $\delta.$ Indeed, by the definition
the support $X$ of a measure $\mu_{0}$ the continuous function $\mu_{0}(\Delta_{\delta}(x))$
on $X$ is point-wise strictly positive, hence globally strictly positive
on $X,$ by compactness). 

In the general case when $X$ may be non-compact we can apply the
previous inequality to $X\cap B_{R}$ and set $\psi:=\phi_{k}-\frac{1}{k}\log(C_{\epsilon}e^{\epsilon k}\int_{X\cap B_{R}}e^{k(\phi_{k}-\phi_{0})}\mu_{0}$
so that $\psi\leq\phi_{0}$ on $X\cap B_{R}.$ Then, by the extremal
definition of $\Pi_{X\cap B_{R}}\phi_{0},$ we get $\psi\leq\Pi_{X\cap B_{R}}\phi_{0}=\Pi_{X}\phi_{0}$
on all of $\R^{n}$ which concludes the proof of the lemma.\end{proof}
\begin{prop}
\label{prop:asymp of free energy for permanant}Let $\mu_{0}$ be
a measure on $\R^{n}$ with support $X$ and $\phi_{0}$ a weight
function on $\R^{n}$ such that $e^{\beta_{*}(\phi_{P}-\phi_{0})}\mu_{0}$
has finite total mass for some non-negative number $\beta.$ Then
\[
\lim_{N\rightarrow\infty}\frac{1}{kN}\log\int_{X^{N}}\mbox{Per}_{k}(x_{1},...,x_{N_{k}})e^{-k\phi_{0}((x_{1})+\cdots+\phi_{0}((x_{1}))}\mu_{0}^{\otimes N}=\int_{P}(\Pi_{X}\phi_{0})^{*}dp=\int_{P}\phi_{0}^{*}dp
\]
and the same limit holds when the integral over $X^{N}$ is replaced
with a sup. Equivalently, for any $u\in C_{b}(\R^{n})$ we have 
\[
-\lim_{N\rightarrow\infty}\frac{1}{kN}\log\int_{X^{N}}e^{-k(H_{\phi_{0}}^{(N)}+u)}\mu_{0}^{\otimes N}=\lim_{N\rightarrow\infty}\frac{1}{N}\inf_{X_{N}}(H_{\phi_{0}}^{(N)}+u)=-\int_{P}(\Pi_{X}\phi_{0})^{*}dp
\]
\end{prop}
\begin{proof}
First note that 
\[
\frac{1}{kN}\log\int\mbox{Per}_{k}(x_{1},...,x_{N_{k}})e^{-k\phi_{0}((x_{1})+\cdots)}\mu_{0}^{\otimes N}=\frac{1}{N}\sum_{p\in P_{\Z/k}}\frac{1}{k}\log\int_{\R^{n}}e^{k(p\cdot x-\phi_{0}(x))}\mu_{0}(x)+\frac{\log N!}{Nk}
\]
Next, we will prove the\emph{ lower bound }in the proposition. First
by the inequality in the previous lemma we may as well, when letting
$k\rightarrow\infty,$ replace the term inside the sum with $v(p):=\sup_{X\cap B_{R}}(p\cdot x-\phi_{0}(x))$
to get 
\[
\frac{1}{kN}\log\int\mbox{Per}_{k}(x_{1},...,x_{N_{k}})e^{-k\phi_{0}((x_{1})+\cdots)}\mu_{0}^{\otimes N}\geq\frac{1}{N}\sum_{p\in P_{\Z/k}}v(p)+o(1)
\]
(also using Stirling's formula, which gives that $\frac{\log N!}{Nk}\rightarrow0).$
Moreover, by the previous we may as well replace the sup over $X\cap B_{R_{k}}$
in the definition of $v(p)$ with a sup over $X$ to get $v=(\Pi\phi_{0})^{*}$
which is bounded on $P.$ Hence, since, 
\begin{equation}
\frac{1}{N}\sum_{p\in P_{\Z/k}}\delta_{p}\rightarrow1_{P}dp\label{eq:weak conv of measure}
\end{equation}
weakly on $P,$ as $k\rightarrow\infty,$ this concludes the proof
of the lower bound:

\[
\lim_{k\rightarrow\infty}\frac{1}{kN}\log\int_{X^{N}}\mbox{Per}_{k}(x_{1},...,x_{N_{k}})e^{-k\phi_{0}((x_{1})+\cdots+\phi_{0}((x_{1}))}\mu_{0}^{\otimes N}\geq\int_{P}(\Pi_{X}\phi_{0})^{*}dp(=\int_{P}(\phi_{0})^{*}dp)
\]
To handle the\emph{ upper bound} let us first, to fix ideas, assume
that $\mu_{0}$ has finite total mass $M.$ Then we can trivially
estimate
\[
\int_{\R^{n}}e^{k(p\cdot x-\phi_{0}(x))}\mu_{0}(x)\leq M\sup_{X}e^{k(p\cdot x-\phi_{0}(x))}
\]
 and conclude as before. To handle the general case we split the integral
over $\R^{n}$ according to the decomposition $\R^{n}=B_{R}\cup B_{R}^{c}$
to get 
\[
\int_{\R^{n}}e^{k(p\cdot x-\phi_{0}(x))}\mu_{0}(x)\leq M\sup_{X\cap B_{R}}e^{k(p\cdot x-\phi_{0}(x))}+\int_{B_{R}^{c}}e^{k(P_{X}\phi_{0}-\phi_{0})}\mu_{0}(x)\cdot C_{\epsilon}e^{\epsilon k}\int_{X\cap B_{R}}e^{k(\phi_{k}-\phi_{0})}\mu_{0},
\]
 using the inequality in the previous lemma in the estimate of the
second term, Now, $\Pi_{X}\phi_{0}-\phi_{0}\leq\phi_{P}-\phi_{0}+C$
and hence, since $\phi_{0}-\phi_{P}\rightarrow\infty$ at $\infty$
in $\R^{n}$ we get $\Pi_{X}\phi_{0}-\phi_{0}\leq(1-\delta)(\phi_{P}-\phi_{0})$
for $\delta$ a sufficiently small number (taking $R$ above sufficiently
large but fixed). Accordingly, assuming that $e^{\beta_{*}(\phi_{P}-\phi_{0})}\mu_{0}$
has finite mass $M'$ for some positive number $\beta_{*}$ we get
for $k$ sufficiently large that 
\[
\int_{\R^{n}}e^{k(p\cdot x-\phi_{0}(x))}\mu_{0}(x)\leq\sup_{X\cap B_{R}}e^{k(p\cdot x-\phi_{0}(x))}(M+MC_{\epsilon}e^{\epsilon k}\int_{X\cap B_{R}}\mu_{0})\leq C'_{\epsilon}e^{\epsilon k}\sup_{X\cap B_{R}}e^{k(p\cdot x-\phi_{0}(x))}
\]
Hence, we get, just as before, that 
\[
\lim_{k\rightarrow\infty}\frac{1}{kN}\log\int_{X^{N}}\mbox{Per}_{k}(x_{1},...,x_{N_{k}})e^{-k\phi_{0}((x_{1})+\cdots+\phi_{0}((x_{1}))}\mu_{0}^{\otimes N}\leq\int_{P}(\Pi_{X}\phi_{0})^{*},
\]
 which concludes the proof of the asymptotics for the integrals in
the theorem. Finally, applying the previous lemma to $\phi_{k}(x)=-\log H^{(N)}(x,x_{2},...,x_{N})$
for any choice of $(x_{2},...,x_{N})\in X^{N-1},$ etc, one coordinate
a time, and arguing as above, also shows that the integral over $X^{N}$
may as well, asymptotically, be replaced with a sup over $X^{N}\cap B_{R}^{N},$
which in turn coincides with the sup over $X^{N}$ itself.
\end{proof}

\subsection{Functionals on convex functions and probability measures}

Fixing a weighted set $(X,\phi_{0})$ we now define, following the
notation in section \ref{sub:Setup-and-assumptions on hamilton} a
functional $\mathcal{F}(u)$ on $C_{b}(\R^{n})$ as the limiting functional
appearing in Prop \ref{prop:asymp of free energy for permanant}:

\begin{equation}
\mathcal{F}(u):=-\int_{P}(\phi_{0}+u)^{*}(p)dp\label{eq:def of functional f of u in ma setting}
\end{equation}
The connection to the Monge-Ampère operator appears as follows (compare
\cite{ber-ber}). Defining 
\[
\mathcal{E}(\phi):=-\int_{P}\phi^{*}(p)dp
\]
we have $\mathcal{E}(\phi_{P})=0$ and 

\[
d\mathcal{E}(\phi)=MA(\phi)
\]
 for any $\phi\in\mathcal{P}(\R^{n})$ such that $\mathcal{E}(\phi)>\infty,$
which by definition means that $\phi$ is in the space $\mathcal{E}_{P}^{1}(\R^{n})$
of all functions with \emph{finite energy.} In particular, integrating
along affine lines in $\mathcal{P}_{+}(\R^{n})$ the functional $\mathcal{E}$
may be written as the following energy type functional 
\begin{equation}
\mathcal{E}(\phi):=\int_{0}^{1}(\phi-\phi_{P})MA(\phi_{0}(1-t)+t\phi)dt,\label{eq:def of energy as mixed in global conv}
\end{equation}
which after expansions and integration over $t$ can be written as
a mixed Monge-Ampère expression (anyway, we will not use this representation).
We may now rewrite \ref{eq:def of functional f of u in ma setting}
as 
\[
\mathcal{F}(u)=\mathcal{E}(\Pi_{X}(\phi_{0}+u)
\]
The following proposition is the key result in the variational approach
to Monge-Ampère equations: 
\begin{prop}
\label{prop:diff theorem}\cite{ber-ber}The functional $\mathcal{F}(u)$
is Gateaux differentiable on $C_{b}(\R^{n})$ and 
\[
d\mathcal{F}_{|u}=MA(\Pi_{X}(\phi_{0}+u))
\]
Similarly, if $\phi_{0}$ has finite energy then the corresponding
statement also holds.\end{prop}
\begin{rem}
In the present real setting the previous proposition can be obtained
from basic properties of the Legendre transform (see \cite{ber-ber}),
but it also holds in the more general complex setting (see Theorem
B in \cite{ber-bou}), where the proof is based on the complex analog
of the following ``orthogonality relation''
\[
\int_{X}MA(\Pi_{X}(\phi))(\Pi_{X}(\phi)-\phi)=0,
\]
 i.e. $\Pi_{X}(\phi)=\phi$ almost everywhere with respect to $MA(\Pi_{X}(\phi)).$
\end{rem}
Given a weighted set $(X,\phi_{0})$ we next define the \emph{weighted
energy} $E_{\phi_{0}}(\mu)$ as the Legendre transform of the functional
$\mathcal{F}$ 
\[
E_{\phi_{0}}(\mu):=\sup_{u\in C_{b}^{0}(X)}(\mathcal{E}(\Pi_{X}(\phi_{0}+u)-\int_{X}u\mu)
\]
 i.e. by the formula \ref{eq:def of E of mu as legendre transform},
but replacing $C^{0}(X)$ with the space $C_{b}(X)$ of bounded continuous
functions on $X$ (recall that we are using a different sign convention
than in the Legendre transform on $\R^{n}$ defined by \ref{eq:def of leg trans of functions}).
\begin{prop}
\label{prop:energy in terms of ma potential }Let $\mu$ be a probability
measure on $\R^{n}$ supported on a the closed set $X.$ Then 
\begin{equation}
E_{\phi_{0}}(\mu):=\sup_{\phi\in\mathcal{P}_{+}(\R^{n})}(\mathcal{E}(\phi)-\int(\phi-\phi_{0})\mu)\label{eq:energy as sup over convex}
\end{equation}
and $E_{\phi_{0}}(\mu)$ is finite iff there exists $\phi_{\mu}\in\mathcal{P}(\R^{n})$
such that 
\begin{equation}
E_{\phi_{0}}(\mu)=\mathcal{E}(\phi_{\mu})-\int(\phi_{\mu}-\phi_{0})\mu,\label{eq:energy in terms of pot}
\end{equation}
which in turn is equivalent to $\phi_{\mu}$ being a \emph{potential
of $\mu$} with finite energy, i.e. $\phi_{\mu}\in\mathcal{E}_{P}^{1}(\R^{n})$
is uniquely determined mod $\R$ by the Monge-Ampère equation 
\[
MA(\phi_{\mu})=\mu
\]
\end{prop}
\begin{proof}
To prove the first formula first observe that since $\phi:=\Pi_{X}(\phi_{0}+u)\leq\phi_{0}+u$
we immediately get the upper bound in \ref{eq:energy as sup over convex}.
To get the lower bound we plug in the following function in the definition
of $E_{\phi_{0}}(\mu):$ $u:=\phi_{\mu}-\phi_{0},$ where $MA(\phi_{\mu})=\mu$
(compare below). By the domination principle for the Monge-Ampère
operator (see the appendix) $\Pi_{X}\phi_{\mu}=\phi_{\mu}$ and hence
$\mathcal{E}(\phi_{\mu})-\int(\phi_{\mu}-\phi_{0})\mu\leq E_{\phi_{0}}(\mu).$
Finally, writing $\phi_{\mu}$ as a decreasing sequence of elements
in $\mathcal{P}_{+}(\R^{n})$ concludes the proof of formula \ref{eq:energy as sup over convex},
by basic continuity properties of $\mathcal{E}$ (see \cite{ber-ber}).
As for formula \ref{eq:energy in terms of pot} it follows immediately
from the variational construction of potentials in \cite{ber-ber}
(compare the proof of Prop \ref{prop:existence of solution to ma mean field}
below). 
\end{proof}
We will often omit the subscript $\phi_{0}$ in the notation for the
weighted energy. Note that, when $\mu$ has compact support we can
decompose 
\begin{equation}
E_{\phi_{0}}(\mu)=E_{0}(\mu)+\int\phi_{0}\mu,\label{eq:decomp formula for energy}
\end{equation}
 where $E_{0}(\mu)$ is independent of $\phi_{0}.$ In fact, the formula
applies to any $\mu$ such that $E_{\phi_{0}}(\mu)<\infty$ (for example,
using an approximation argument; compare \cite{berm2} for the complex
setting). 

The next proposition gives the relation to the theory of optimal transport:
\begin{prop}
\label{prop:energy as cost}Let $\mu$ be a probability measure on
$\R^{n}.$ Then 
\[
E_{\phi_{0}}(\mu)=C_{\phi_{0}}(\mu)
\]
 where $C_{\phi_{0}}(\mu)$ is the Monge-Kantorovich cost functional
defined with respect to the target measure $\lambda_{P}:=1_{P}dp$
on the target convex body $P$ and the cost function $c(x,p)=-x\cdot p+\phi_{0}(x)$
(see section \ref{sub:Optimal-transport-theory}). \end{prop}
\begin{proof}
By definition, 
\[
E(\mu)=\sup_{u\in C_{b}(X)}(\int_{P}-(\Pi_{X}(\phi_{0}+u))^{*}(p)dp-\int u\mu)
\]
Now set $v:=(\Pi_{X}(\phi_{0}+u))^{*}$ which is a bounded function
(as explained in section \ref{sub:Weigted-sets-and}). Using the extremal
property of the Legendre transform we may rewrite the previous line
as 
\[
E(\mu)=\sup_{v\in C_{b}(P),u\in C_{b}(\R^{n})}(\int_{P}-vdp-\int u\mu),
\]
 where the sup ranges over all $u$ and $v$ such that $-v-u\leq c(x,p),$
where $c(x,p)=-x\cdot p+\phi_{0}(x).$ According to the general Kantorovich
duality theorem \cite{vi2} this means that $E_{\phi_{0}}(\mu)=C_{\phi_{0}}(\mu)$
if the following condition on the cost function is satisfied: $-c(x,y)\leq f(x)+g(p)$
for some functions $f$ and $g$ such that $f\in L^{1}(X,\mu)$ and
$g\in L^{1}(P,dp).$ In the present setting we have $-c(x,p)\leq\Pi_{X}(\phi_{0}){}^{*}$
on $X\times P,$ where we recall that $\Pi_{X}(\phi_{0}){}^{*}$ is
bounded on $P.$ Hence we may take $f(x)=0$ and $g(p)=\Pi_{X}(\phi_{0}){}^{*}$
which thus concludes the proof of the alternative formula.
\end{proof}
Given a weighed measure $(\mu_{0},\phi_{0})$ we now define, following
the notation in previous sections, the corresponding free energy functional
on $\mathcal{M}_{1}(X),$ where $X$ is the support of $\mu_{0},$
by 
\begin{equation}
F_{\beta}=E_{\phi_{0}}+\frac{1}{\beta}D_{\mu_{0}},\label{eq:free energy in ma setting}
\end{equation}
 for any given $\beta\in]0,\infty].$
\begin{prop}
The free energy functional $F_{\beta}$ defines a good rate functional
on $\mathcal{M}_{1}(X).$\end{prop}
\begin{proof}
Both functionals $E$ and $D$ are lsc, since they may be realized
as Legendre transforms and to show that $F_{\beta}$ is a good rate
functional on $\mathcal{M}(X)$ we thus only need to verify its properness
(which is automatic when $X$ is compact). This could be proven directly,
but anyway it is a general fact that a functional is a good rate functional
if it is the rate functional of a LDP for which exponential tightness
holds (see Lemma 1.2.18 in \cite{de-ze}) and the latter properties
will be established in the proof of Theorem \ref{thm:main for perman intro}.
\end{proof}

\subsection{Proof of Theorem \ref{thm:main for perman intro}}

The case when $\beta=\infty$ follows immediately from Prop \ref{prop:asymp of free energy for permanant}
and the Gärtner-Ellis theorem (compare the proof in the complex setting
considered in )\cite{berm2}. We hence consider the case when $\beta<\infty$
and to simplify the exposition we assume that $\beta_{N}/k=\beta,$
but the proof in the general case is essentially the same. 

Let us first consider the case when $X$ is compact. Setting $\beta_{N}=k$
Proposition \ref{prop:asymp of free energy for permanant} then allows
us to apply the general Theorem \ref{thm:interacting part intro}
to deduce the LDP in Theorem \ref{thm:main for perman intro}. The
fact that the potential of the minimizer $\mu_{\beta}$ solves the
Monge-Ampere equation \ref{eq:ma eq weak intro} is proved in section
\ref{prop:existence of solution to ma mean field}. But it can, for
$X$ compact, also be seen as a special case of Theorem \ref{thm:existence of sol to general mean f}.
Indeed, the latter results translates into $\mu=MA(\Pi_{X}u),$ for
some continuous function $u$ on $X$ satisfying $MA(\Pi_{X}(\phi_{0}+u)=e^{\beta u}\mu_{0},$
where by the limiting construction used in the proof of \ref{thm:existence of sol to general mean f}
$\phi_{0}+u$ is the restriction to $X$ of a function in the space
$\mathcal{P}(\R^{n}).$ But then it follows immediately that $\Pi_{X}(\phi_{0}+u)=\phi_{0}+u$
on $X$ and hence $\phi:=\Pi_{X}(\phi_{0}+u)$ solves the Monge-Ampère
equation \ref{eq:ma eq weak intro} and $\mu_{\beta}=MA(\phi),$ as
desired.

In the general non-compact case we can use a variant of the localization
argument used in the proof of Prop \ref{prop:asymp of free energy for permanant}
to verify that the tightness assumption is satisfied and to reduce
the problem to the compact set $B_{R}\cap X.$ Indeed, let us first
check the validity of the analog of Theorem \ref{thm:conv of one-point cor measur}.
For the upper bound we still get, just as before, 
\[
F^{(N)}(\mu^{(N)})\leq\inf_{\mu\in\mathcal{M}_{1}(\R^{n})_{c}}(E(\mu)+\frac{1}{\beta}D(\mu)),
\]
 where $\mathcal{M}_{1}(\R^{n})_{c}$ denotes the space of all\emph{
compactly supported} probability measures on $\R^{n}.$ Anyway, by
a simple approximation argument, where $\mu$ gets replaced with $1_{B_{R}}\mu/\mu(B_{R})$
for a sequence of $R\rightarrow\infty$ we may as well take the infimum
appearing in the right hand side above over all of $\mathcal{M}_{1}(\R^{n})$
(see \cite{berm2} for the complex case)

To prove the lower bound we first observe that the first marginals
$\mu_{1}^{(N)}$ define a\emph{ tight} sequence. Indeed, by the inequality
in Lemma \ref{lem:b-m for convex etc} we have that 
\[
\mu_{1}^{(N)}\leq Ce^{\beta(\pi_{X}(\phi)-\phi_{0}))}\mu_{0}\leq C'e^{\beta(\phi_{P}-\phi_{0})}\mu_{0}
\]
By assumption, the measure appearing in the rhs above has finite total
mass on $\R^{n}$ and hence the tightness of the sequence $\mu_{1}^{(N)}$
follows. But then it follows from the standard weak pre-compactness
of tight sequences that $\mu_{1}^{(N)}$ has a limit point $\mu_{*}$
in $\mathcal{M}_{1}(\R^{n}).$ Moreover, repeating the argument for
the lower bound in the proof of Theorem \ref{thm:conv of one-point cor measur}
and using the asymptotics in \ref{prop:asymp of free energy for permanant}
for the infimum of $(H_{\phi_{0}}^{(N)}+u),$ for any $u\in C_{b}(\R^{n}),$
gives 
\[
E(\mu_{*})+D(\mu_{*})/\beta\leq F^{(N)}(\mu^{(N)})\leq\inf_{\mu\in\mathcal{M}_{1}(\R^{n})}(E(\mu)+D(\mu)/\beta),
\]
 But then it follows from Prop \ref{prop:existence of solution to ma mean field}
that $\mu_{*}$ coincides with the unique minimizer of $E(\mu)+D(\mu)/\beta$
and that its potential satisfies the desired equation. 

Finally, to prove the LDP in the non-compact case we just have to
verify the exponential tightness of the corresponding sequence $\Gamma_{k}.$
More precisely, that we need to prove is the space $\mathcal{M}_{1}(X)$
may be exhausted by compact subsets $\mathcal{F}_{\alpha}$ for $\alpha>0$
such that $\limsup_{k\rightarrow\infty}\log(\Gamma_{k}(\mathcal{M}_{1}(X)-\mathcal{F}_{\alpha})/\beta N_{k}<-\alpha.$
To prove this we let $\mathcal{F}_{\alpha}$ be the set of all measures
$\mu$ on $\mathcal{M}_{1}(X)$ such that $\int(\phi_{0}-\Pi_{X}\phi_{0})\mu\leq3\alpha.$
Since, by assumption $\phi_{0}-\Pi_{X}\phi_{0}\rightarrow\infty$
at infinity in $\R^{n}$$,$ the set $\mathcal{F}_{\alpha}$ is indeed
compact. By definition
\[
\Gamma_{k}(\mathcal{M}_{1}(X)-\mathcal{F}_{\alpha})=\frac{1}{Z_{N}}\int_{\{\phi_{0}-\Pi_{X}\phi_{0}>3\alpha N\}}\mbox{Per}_{k}(x_{1},...,x_{N_{k}})^{\beta/k}\mu_{0}^{\otimes N_{k}}
\]
Now, applying Lemma \ref{lem:b-m for convex etc} $N$ times (i.e.
one ``coordinate at time'' ) the density in the previous integral
may be estimated from above by $C_{\epsilon}^{N_{k}}e^{\epsilon N_{k}\beta}e^{-\beta(\phi_{0}-P_{F}\phi_{0})}$
for some fixed small $\epsilon>0$ (taken so that $\epsilon<\alpha/2).$
Hence, decomposing 
\[
e^{-\beta(\phi_{0}-\Pi_{X}\phi_{0})}=e^{-\frac{1}{2}\beta(\phi_{0}-\Pi_{X}\phi_{0})}e^{-\frac{1}{2}\beta(\phi_{0}-\Pi_{X}\phi_{0})}\leq e^{-\frac{1}{2}\beta N3\alpha}e^{-\frac{1}{2}\beta(\phi_{0}-\phi_{P})}C^{\beta}
\]
and integrating wrt $\mu_{0}^{\otimes N_{k}}$ (and using that $e^{-\beta_{*}(\phi_{0}-\phi_{P})}\mu_{0}$
is assumed to have finite total mass wrt any positive number $\beta_{*}$
) finishes the proof of the exponential tightness.
\begin{rem}
\label{rem:log per is cost asympt}In the case when $X$ is compact
we can directly combine Proposition \ref{prop:asymp of free energy for permanant}
with the Gärtner-Ellis theorem to deduce that $e^{-kH^{(N)}}\mu_{0}^{\otimes N}$
satisfies a LDP principle with speed $kN$ and rate functional $E(\mu).$
Then, arguing as in section \ref{sub:An-alternative-proof} (and using
that, by Prop \ref{prop:energy as cost}, $E=C)$ it follows that
\[
\sup_{X^{N}}\left|-\frac{1}{k}\log\mbox{Per \ensuremath{(x_{1},...,x_{N})}}-C(\frac{1}{N}(\delta_{x_{1}}+\cdots+\delta_{x_{N}}))\right|=0,
\]
 where $C$ denotes the Monge-Kantorovich total cost functional associated
to the cost function $c(x,p)=-\left\langle x,p\right\rangle .$ 
\end{rem}

\subsubsection{Proof of Cor \ref{cor:first cor intro}}

By Theorem \ref{thm:main for perman intro} the corresponding one-point
correlation measures $\rho_{1}^{(N)}\mu_{0}$ converge weakly to $\mu_{\beta}.$
Since, the latter measure may be written as $\mu_{\beta}=e^{\beta(\phi-\phi_{0})}\mu_{0}$
for $\phi$ the unique solution in $\mathcal{P}(\R^{n})$ to equation
\ref{eq:ma eq weak intro} this means that 
\[
\rho_{1}^{(N)}\mu_{0}\rightarrow e^{\beta(\phi-\phi_{0})}\mu_{0}
\]
 weakly. Now $\phi^{(N)}:=\frac{1}{\beta}\log\rho_{1}^{(N)}-\phi_{0}$
is a sequence in $\mathcal{P}(\R^{n})$ such that $\int e^{\beta\phi^{(N)}}\mu_{0}=1$
and hence, by the inequality in Lemma \ref{lem:b-m for convex etc},
fixing a point $x_{0}$ in $X$ gives that $\phi^{(N)}(x_{0})$ is
uniformly bounded from above. Hence, by Prop \ref{prop:compactness of space of convex},
we may, after perhaps passing to a subsequence, assume that $\phi^{(N)}\rightarrow\phi_{*}$
locally uniformly for some function $\phi_{*}$ in $\mathcal{P}(\R^{n}).$
But then it follows from the convergence above that $\phi_{*}=\phi$
almost everywhere wrt $\mu_{0}$ and hence everywhere if the support
of $X$ is all of $\R^{n}.$ Since, we may repeat the same argument
for any subsequence of $\phi^{(N)}$ it follows from the uniqueness
of the solution $\phi$ that the whole sequence $\phi^{(N)}$ converges
to $\phi,$ as desired.

\subsubsection{Proof of Cor \ref{cor:sec cor intro}}

By assumption the support $X$ of $\mu_{0}:=\rho_{0}1_{X}dx$ is compact
and it is the closure of its interior. First note that $-\phi^{(N)}(x)$
represents the differential of $E_{N}(\mu_{0}):=\int_{X^{N}}H^{(N)}(\mu_{0})^{\otimes N}.$
More precisely, if we fix a smooth probability density $\rho$ with
compact support in the interior of $X$ and set $\mu_{t}:=\mu_{0}+t(\rho-\rho_{0})dx,$
for $t$ sufficiently small, i.e. $|t|\leq\epsilon$ such that $\mu_{t}\in\mathcal{M}_{1}(X),$
then a direct calculation gives $dE_{N}(\mu_{t})/dt_{|t=0}=-\int_{X}\phi^{(N)}(x)(\rho-\rho_{0})dx.$
Next, we note that, since $X$ is compact $E(\mu_{t})<\infty,$ where
$E(\mu)$ denotes the unweighted energy (i.e. $\phi_{0}=0)$ and $\lim_{N\rightarrow\infty}E_{N}(\mu_{t})=E(\mu_{t})$
(by Theorem \ref{thm:existence of mean energ}). Hence, by Lemma \ref{lem:maxim is differ}
and the uniqueness mod $\R$ of potentials it follows that $E(\mu_{t})$
is differentiable and $dE(\mu_{t})/dt_{|t=0}=-\int_{X}\phi_{\mu_{0}}(x)(\rho-\rho_{0})dx,$
where $\phi_{\mu_{0}}$ is any potential for $\mu_{0}$ which we may
as well take to be the one uniquely determined by the normalization
condition $\int\phi_{\mu_{0}}\mu_{0}=0.$ Since $E(\mu_{t})$ is moreover
convex it then follows from basic properties of convex functions that
$\lim_{N\rightarrow\infty}dE_{N}(\mu_{t})/dt=dE(\mu_{t})/dt$ (compare
the proof of Lemma \ref{lem:conv of abstr fekete}), i.e. 
\begin{equation}
\lim_{N\rightarrow\infty}\int_{X}\phi^{(N)}(x)(\rho-\rho_{0})dx=\int_{X}\phi_{\mu_{0}}(x)(\rho-\rho_{0})dx\label{eq:proof of cor two}
\end{equation}
Now, by assumption, $\int\phi^{(N)}\mu_{0}=0$ and thus we may, by
Prop \ref{prop:compactness of space of convex}, after perhaps passing
to a subsequence, assume that $\phi^{(N)}\rightarrow\phi_{*}$ in
$\mathcal{P}(\R^{n})$ such that $\int\phi_{*}\mu_{0}=0.$ But since
$\phi_{\mu_{0}}$ satisfies the same normalization condition and \ref{eq:proof of cor two}
holds for any $\rho$ as above we conclude that $\phi_{*}=\phi_{\mu_{0}}$
on the interior of $X.$ In particular, $MA(\phi_{*})=\mu_{0}$ on
the interior of $X$ and since $\int_{\R^{n}}MA(\phi_{*})\leq1$ and
$\int_{\R^{n}}\mu_{0}=1$ it follows that $MA(\phi_{*})=\mu_{0}$
on all of $\R^{n}$ and hence, by uniqueness of normalized potentials,
 $\phi_{*}=\phi_{\mu_{0}}$ on all of $\R^{n},$ as desired. 

Next, it follows form the regularity result in \cite{ca-2} that $\phi$
is $C^{1,\alpha}-$smooth in the interior of $X$ for some $\alpha>0.$
We briefly recall the argument: first the solution $\phi$ in $\mathcal{P}(\R^{n})$
of the Monge-Ampère equation $MA(\phi)=\mu_{0}$ (in the sense of
Alexandrov) has, by assumption, a Monge-Ampère measure $MA(\phi)$
which is absolutely continuous wrt Lebesgue measure and vanishes on
the complement of $X.$ Hence, its restriction to the interior of
$X$ defines a Brenier solution on $X,$ i.e. it satisfies the MA-equation
on the interior of $X$ in the weak sense of Brenier and the almost
everywhere defined map $\nabla\phi$ from the interior of $X$ to
$P$ is almost everywhere surjective. But as shown in \cite{ca-2}
any Brenier solution $\phi$ (which is in fact uniquely determined)
is $\mathcal{C}^{1,\alpha}-$smooth in the interior of $X$ (the point
is that, as shown in \cite{ca-2}, it is strictly convex and then
the regularity follows from \cite{ca}). Finally, since $\phi$ is
convex and differentiable on the interior of $X$ the previous convergence
of $\phi^{(N)}$ implies, by basic properties of convex functions,
point-wise convergence everywhere on the interior of $X$ for the
corresponding gradients.

\subsection{\label{sub:Existence-and-extremal}Variational properties of Monge-Ampère
equations and regularity }

In this section we will establish some properties of the Monge-Ampère
equations and convex envelopes studied above, most of which can be
reduced to essentially known results.
\begin{prop}
\label{prop:ma of envelope minimizes energy}Let $(X,\phi_{0})$ be
a weighted set and denote by $\phi_{e}:=\Pi_{X}\phi_{0}$ the corresponding
convex envelope. Then
\begin{itemize}
\item $\mu_{e}:=MA(\phi_{e})$ is the unique minimizer of the functional
$\mu\mapsto E_{\phi_{0}}(\mu)$ on $\mathcal{M}_{1}(X).$ 
\item If $X=\R^{n}$ and $\phi_{0}$ is smooth, then $\phi_{e}$ is locally
$\mathcal{C}^{1,1},$ i.e. $\nabla\phi_{e}$ is locally a Lipschitz
map and 
\[
\mu_{e}=1_{D_{0}}\det(\frac{\partial^{2}\phi}{\partial x_{i}\partial x_{j}})dx.
\]
 where $D_{0}$ is the compact set where $\phi_{e}=\phi_{0}$( and
the density above is point-wise defined a.e. on $D_{0}).$
\end{itemize}
\end{prop}
\begin{proof}
\emph{The minimizing property:} By construction $E_{\phi_{0}}(\mu)$
is the Legendre transform of the functional $\mathcal{F}_{X}(u):=\mathcal{E}(\Pi_{X}(\phi_{0}+u)$
and hence, since the latter functional is Gateaux differentiable,
it follows from general properties of Legendre transforms (Lemma \ref{lem:maxim is differ})
that $E_{\phi_{0}}(\mu)$ admits a unique minimizer, which is given
by the differential $d\mathcal{E}(\Pi_{X}(\phi_{0}+u)$ at $u=0.$
Finally, by Prop \ref{prop:diff theorem} the latter differential
is equal to $MA(\Pi_{X}\phi_{0}).$ 

\emph{Regularity:} This is a special case of the regularity results
for the generalized Lelong class obtained in \cite{be0}, modeled
on the classical approach of Bedford-Taylor. \end{proof}
\begin{prop}
\label{prop:existence of solution to ma mean field}Let $\mu_{0}$
be a measure on $\R^{n}$ and $\phi_{0}$ a continuous function on
$\R^{n}$ such that $\phi_{0}-\phi_{P}$ is proper and assume that
$\int e^{\beta(\phi_{P}-\phi_{0})}\mu_{0}<\infty.$ Then 
\begin{itemize}
\item The Monge-Ampère equation \ref{eq:ma eq weak intro} admits a solution
$\phi_{\beta}$ in $\mathcal{P}_{\text{ }}(\R^{n})$ of finite energy,
i.e. $\mathcal{E}(\phi_{\beta})>-\infty.$ 
\item Any two solutions of full Monge-Ampère mass coincide up to an additive
constant
\item The probability measure $\mu_{\beta}:=MA(\phi_{\beta})$ is the unique
minimizer of the corresponding free energy functional $F_{\beta}.$
\item If $\phi_{0}$ is smooth and $\mu_{0}=\rho_{0}dx$ for a strictly
positive smooth function $\rho_{0},$ then the solution $\phi_{\beta}$
is also smooth 
\end{itemize}
\end{prop}
\begin{proof}
\emph{Existence:} The existence of a weak solution is well-known,
at least when $\mu_{0}$ is absolutely continuous wrt $dx$ \cite{ba},
but it may be illuminating to give a variational proof of the general
case, in the spirit of the present paper. The point is that, following
the variational approach in \cite{ber-ber}, we just need to verify
that the following coercivity estimate holds: 
\[
\mathcal{D}(\phi):=-\mathcal{E}(\phi)+\frac{1}{\beta}\log\int e^{\beta(\phi-\phi_{0})}\mu_{0}\geq-\mathcal{E}(\phi)+(\phi-\phi_{0})(x_{0})
\]
 where $x_{0}$ is a fixed point in the support $X$ of $\mu_{0}.$
But this follows immediately from the inequality 
\[
\sup_{X}e^{\beta(\phi-\phi_{0})}\leq C\int_{X\cap B_{R}}e^{\beta(\phi-\phi_{0})}\mu_{0},
\]
 which is a direct consequence of the inequality in Lemma \ref{lem:b-m for convex etc}.
With the coercivity inequality in place the solution $\phi$ may be
obtained as a minimizer of the functional $\mathcal{D}.$ Indeed,
since $\mathcal{D}(\phi+c)=\mathcal{D}(\phi)$ we can take a sequence
of functions $\phi_{j}\in\mathcal{P}_{+}(\R^{n})$ such that $\mathcal{D}(\phi_{j})$
converges to the supremum of $\mathcal{D}$ on $\mathcal{P}(\R^{n})$
and such that $\phi_{j}$ is normalized, i.e. $(\phi_{j}-\phi_{0})(x_{0})=0).$
By the Arzelà-Ascoli theorem $\phi_{j}$ converges, after perhaps
passing to a subsequence, locally uniformly to $\phi_{\infty}$ in
$\mathcal{P}(\R^{n}).$ Moreover, by the coercivity estimate above
$\phi_{\infty}$ has finite energy. Defining
\[
\mathcal{\tilde{D}}(\phi):=-\mathcal{E}(\Pi_{X}\phi)+\frac{1}{\beta}\log\int e^{\beta(\phi-\phi_{0})}\mu_{0}
\]
 we have $\mathcal{\tilde{D}}(\phi)\geq\mathcal{\tilde{D}}(\Pi_{X}\phi)$
and hence, if $u$ is any given smooth compactly supported function
on $\R^{n}$ then the function $t\mapsto\mathcal{\tilde{D}}(\phi+tu)$
on $\R$ has a maximum at $t=0$ and is differentiable by Prop \ref{prop:diff theorem}.
Hence, its derivative at $t=0$ vanishes and the formula for the differential
in Prop \ref{prop:diff theorem} then shows that $\phi_{\infty}$
satisfies the desired equation up to a multiplicative normalization
factor which can be removed by adding a constant to $\phi_{\infty}.$ 

\emph{Uniqueness:} The uniqueness of finite energy solutions can be
shown by convexity arguments, but, anyway, the general case follows
form the comparison principle for $\mbox{MA .}$Indeed, if $u$ and
$v$ are in $\mathcal{E}(X)$ then the comparison principle says 
\[
\int_{\{u<v\}}MA(v)\leq\int_{\{u<v\}}MA(u)
\]
But if $u$ and $v$ are solutions of equation \ref{eq:ma eq weak intro}
then it must be that $u=v$ a.e wrt the measure $\mu_{0}$ and hence
$MA(u)=MA(v)=\mu_{0},$ which implies that $u-v$ is constant, by
the uniqueness of normalized potentials of any probability measure. 

\emph{Minimizing property:} This is proved precisely as in the proof
of Theorem \ref{thm:existence of sol to general mean f}.

\emph{Regularity:} This is proved exactly as in \cite{ber-ber}, using
Caffarelli's interior regularity results. Briefly, since $\phi$ has
finite energy the image of the corresponding subgradient map to $P$
is surjective in the almost everywhere sense. But this implies that
$\phi$ is proper (using that we may assume that $0$ is in an interior
point in $P),$ i.e. the sublevel sets $\Omega_{R}:=\{\phi<R\}$ are
bounded convex domains exhausting $\R^{n}.$ On $\Omega_{R}$ the
function $u:=\phi-R$ defines a function in $C_{0}(\bar{\Omega}_{R}),$
vanishing on the boundary to which we may apply the regularity results
in \cite{ca,ca-2,ca0,ca1} to deduce that $u$ is smooth (see \cite{ber-ber}
for the complete argument).\end{proof}
\begin{rem}
Note that if there is a solution $\phi$ in $\mathcal{P}_{\text{+ }}(\R^{n})$
then necessarily $\int e^{\beta(\phi_{P}-\phi_{0})}\mu_{0}\leq C\int e^{\beta(\phi-\phi_{0})}\mu_{0}=C\int MA(\phi)\leq C<\infty.$ \end{rem}
\begin{prop}
\label{prop:conv of solutions as beta tends to infy}Let $(\mu_{0},\phi_{0})$
be a weighted measure and let $X$ be the support of $\mu_{0}.$ Denote
by $\phi_{\beta}$ be the unique solution in in $\mathcal{P}_{\text{ }}(\R^{n})$
to the corresponding equation \ref{eq:ma eq weak intro}. Then $\phi_{\beta}$
converges, as $\beta\rightarrow\infty,$ locally uniformly to the
envelope $\Pi_{X}\phi_{0}(:=\phi_{e})$ iff $\mu_{e}$ has finite
entropy with respect to $\mu_{0}.$ In particular, this is the case
if $X=\R^{n}$ and $\phi_{0}$ is smooth.\end{prop}
\begin{proof}
Let us first verify that the family $\mu_{\beta}:=MA(\phi_{\beta})$
is tight. By the inequality in Lemma \ref{lem:b-m for convex etc}
we have, since $\int e^{\beta(\phi_{\beta}-\phi_{0})}\mu_{0}=1,$
that 
\begin{equation}
(\phi_{\beta}-\phi_{0})\leq C/\beta+\phi_{e}-\phi_{0}\leq C/\beta+C'+\phi_{P}-\phi_{0}.\label{eq:proof of phi beta converges to env}
\end{equation}
 Now, by assumption, $\phi_{P}-\phi_{0}\rightarrow-\infty$ and hence
there exists $\delta>0$ such that 
\[
x\in X-B_{R}\implies(\phi_{\beta}-\phi_{0})\leq(1-\delta)(\phi_{P}-\phi_{0})
\]
for some large ball $B_{R}$ (where we may assume that $(\phi_{P}-\phi_{0})<0).$
But then $\int_{X-B_{R}}e^{\beta(\phi_{\beta}-\phi)}\mu_{0}\leq\int_{X-B_{R}}e^{\beta(1-\delta)(\phi_{P}-\phi)}\mu_{0}=\epsilon_{R},$
where $\epsilon_{R}\rightarrow0,$ as $R\rightarrow\infty$ (since,
by assumption, the integral is finite for some $\beta).$ Thus the
family $\mu_{\beta}$ is tight as desired. In particular, the family
admits a limit point $\mu_{\infty}$ in $\mathcal{M}_{1}(X)$ and
we next show that it coincides with $\mu_{e},$ the unique minimizer
of the energy functional $E$ on $\mathcal{M}_{1}(X).$ This can be
proved by following the argument given in the complex case in \cite{berm2}
(Theorem 3.13). But here we note that a simpler argument can be given
in the real setting. Let us first assume that $D_{\mu_{0}}(\mu_{e})<\infty.$
Since $\mu_{\beta}$ minimizes the functional $F_{\beta}$ we then
get 
\[
E(\mu_{e})=\lim_{\beta\rightarrow\infty}F_{\beta}(\mu_{e})\geq\limsup_{\beta\rightarrow\infty}F_{\beta}(\mu_{\beta})=E(\mu_{\infty}),
\]
 using that $E$ is lower semi-continuous and that $D_{\mu_{0}}(\mu_{\beta})\leq C$
(by \ref{eq:proof of phi beta converges to env}) to get that last
inequality. Hence, $\mu_{\infty}=\mu_{e}=MA(\phi_{e})$ by Prop \ref{prop:ma of envelope minimizes energy}.
In other words $MA(\phi_{\beta})\rightarrow MA(\phi_{e})$ weakly.
But after passing to a subsequence we may assume that $\phi_{\beta}\rightarrow\phi_{\infty}$
for some $\phi_{\infty}\in\mathcal{P}(\R^{n})$ and since $\phi_{\beta}$
has finite energy and hence full Monge-Ampère mass it follows that
$MA(\phi_{\beta})\rightarrow MA(\phi_{\infty})=\mu_{e}$ (compare
\cite{ber-ber}). But by the uniqueness mod $\R$ of potentials it
then follows that $\phi_{\infty}=\phi_{e}+C.$ Finally, to see that
$C=0$ we set $u:=\phi_{\infty}-\phi_{0}$ which is continuous and
bounded from above. Hence, by the equations for $\phi_{\beta}$ we
have that $0=\lim_{\beta\rightarrow\infty}(\log\int e^{\beta u}\mu_{0})/\beta=\sup_{X}u,$
which forces $C=0,$ using that $\sup_{X}(\Pi_{X}\phi_{0}-\phi_{0})=0$
(indeed, the incidence set $D_{\phi_{0}}$ is non-zero, since it contains
the support of the probability measure $MA(\Pi_{X}\phi_{0})).$ To
get the converse statement we assume that $\phi_{\beta}\rightarrow\phi_{e}.$But
then $\mu_{\beta}\rightarrow\mu_{e}$ and since $D_{\mu_{0}}$ is
lsc and $D_{\mu_{0}}(\mu_{\beta})\leq C$ it thus follows that $D_{\mu_{0}}(\mu_{e})<\infty.$
\end{proof}

\section{\label{sec:General-target-measures}General target measures and random
allocation of the target points}

The proof of Theorem \ref{thm:main for perman intro} in fact applies
to a more general setting where $-x\cdot p$ is replaced by a function
$c(x,p)$ and $\lambda_{P}$ is replaced with a probability measure
$\nu$ on $\R^{n}.$ One furthermore needs to fix a sequence $\beta_{N}^{*}$
(playing the role of $k)$ such that 
\[
\beta_{N}^{*}\rightarrow\infty,\,\,\,\,\frac{\log N!}{N\beta_{N_{*}}}\rightarrow\infty
\]
and a sequence of $N-$tuples $p^{(N)}:=(p_{1}^{(N)},...,p_{N}^{(N)})$
such that 
\begin{equation}
\frac{1}{N}\sum_{i=1}^{N}\delta_{p_{i}^{(N)}}\rightarrow\nu,\label{eq:weak approx of nu}
\end{equation}
 weakly as $N\rightarrow\infty.$ Given this data one then replaces
$\mbox{Per}(x_{1},...,x_{N})$ with the permanent
\begin{equation}
\mbox{Per}_{c}(x_{1},...,x_{N}):=\mbox{Per}(e^{-\beta_{N_{*}}c(x_{i},p_{j})})_{i,j\leq N}\label{eq:per wrt c}
\end{equation}
and sets 
\[
\mu_{\beta_{N}}^{(N)}:=\frac{\mbox{Per}_{c}(x_{1},...,x_{N})^{\beta_{N}/\beta_{N^{*}}}}{Z_{N,\beta_{N}}}\mu_{0}^{\otimes N}
\]
Under suitable regularity assumptions on $c(x,y$) and $\nu$ the
previous proof of Theorem \ref{thm:main for perman intro} generalize
verbatim to this more general setting (for example, the proof applies
if $c(x,p)$ is continuous and uniformly Lipschitz wrt $x$ as $p$
ranges over the support of $\nu$ and $\nu$ is absolutely continuous
with respect to $\lambda_{P}).$ The deterministic measure $\mu_{\beta}$
appearing in the limit then coincides with the minimizer of the corresponding
free energy functional defined with respect to the cost functional
$C(\mu,\nu)$ (but for a general cost $c(x,p)$ it can not be directly
linked to a Monge-Ampère equation). The key point is that the proof
of Prop \ref{prop:asymp of free energy for permanant} still applies
if one replaces the Legendre transform $\phi^{*}(p)$ with $\phi^{*,c}(p):=\sup_{x\in\R^{n}}c(x,p)-\phi(p)$
and the measure $\lambda_{P}$ with $\nu$ (up to signs this is the
same transform as the one appearing in \cite{g-mc}). Then the Kantorovich
duality argument used in the proof of Prop \ref{prop:energy as cost}
shows that the corresponding functional $E(\mu$) coincides with the
optimal cost functional $C(\mu,\nu),$ determined by $c(x,p)$ (see
the appendix).

\subsection{Random allocation of target points and quenched variables}

Let us in particular consider the case when we still have $c(x,p)=-x\cdot p,$
but replacing $\lambda_{P}$ with a general target measure $\nu$
absolutely continuous with respect to $\lambda_{P},$ where, as before,
$P$ denotes a given convex body. As before we also assume given a
measure $\mu_{0}$ on $\R^{n}$ and for simplicity we will assume
that its support $X$ is compact and that the weight $\phi_{0}$ vanishes
identically. We will denote by $MA_{\nu}$ the modified Monge-Ampère
operator associated to $\nu,$ i.e. for any Borel set $E$ we set
$MA_{\nu}(\phi)(E):=\nu((\nabla\phi)(E)).$ 

To the data $(\mu_{0},\nu,\beta)$ we may now associate the following
Monge-Ampère equation: 
\begin{equation}
MA_{\nu}(\phi)=e^{\beta\phi}\mu_{0},\label{eq:weak ma eq with nu}
\end{equation}
 assuming as before that $\nabla\phi$ maps $\R^{n}$ into $P.$ If
$\nu=1_{P}e^{-\psi_{0}(p)}dp$ and $\phi$ is smooth the previous
equation just means that 
\[
MA(\phi)e^{-\psi_{0}(\nabla\phi)}=e^{\beta\phi}\mu_{0}
\]
In particular, in the case $\beta=0$ the corresponding equation appears
in the optimal transport problem defined with respect to the target
measure $\nu.$ We will denote the corresponding optimal cost function
by $C(\mu,\nu),$ which, as before, will be considered as a functional
of $\mu.$ If one would also fix a sequence of $p^{(N)}$ approximating
$\nu$ in the sense of \ref{eq:weak approx of nu}, then, as explained
above, the previous results apply to this more general setting. However,
it is also interesting to see that there is a variant of this setting
which does not depend on fixing a sequence of $p^{(N)}$ and to which
we next turn. The idea is to view all the previously defined objects,
such as $H^{(N)},$ $\mu_{\beta_{N}}^{(N)}$ etc as random variables
on $(P^{N},\nu^{\otimes N}).$ This means that we view the variables
$p_{i}$ appearing in 
\[
\mbox{Per}(x_{1},...,x_{N},p_{1},...,p_{N}):=\mbox{Per}(e^{-\beta_{N_{*}}(x_{i}\cdot p_{j})})_{i,j\leq N}
\]
as independent random variables identically distributed according
to probability measure $\nu.$ In other words we perform a \emph{random
allocation} of the $p_{i}:$s according to the measure $\nu$ (which
is somewhat related to the setting considered in \cite{h-s}). In
the terminology appearing in the mathematics of disordered systems
we thus view the variables $(p_{1},...,p_{N})$ as \emph{quenched}
(i.e. frozen); compare \cite{bo}. 

We will denote by $\E$ expectations defined with respect to the ensemble
$(P^{N},\nu^{\otimes N}).$ The previous arguments can then be adapted
to prove the following variant of Theorem \ref{thm:main for perman intro}
(or rather Theorem \ref{thm:conv of one-point cor measur}):
\begin{thm}
\label{thm:conv of one point cor with random alloc}Assume given data
$(\mu_{0},\nu,\beta)$ as above and a sequence $\beta_{N}\rightarrow\beta>0.$
Then the following weak convergence of measures on $X$ holds, as
$N\rightarrow\infty:$ 
\[
\E(\int_{X^{N-1}}\mu_{\beta_{N}}^{(N)})\rightarrow MA(\phi_{\beta}),
\]
 where $\phi_{\beta}$ is a solution to the equation \ref{eq:weak ma eq with nu}.
\end{thm}
Note that we may equivalently view $\E(\int_{X^{N-1}}\mu_{\beta_{N}}^{(N)})$
as the expectation of the empirical measure $\frac{1}{N}\sum_{i=1}^{N}\delta_{x_{i}}$
with respect to the following probability measure on $X^{N}\times P^{N}:$
\[
\gamma_{\beta_{N}}^{(N)}:=\frac{\mbox{Per}(x_{1},...,x_{N},p_{1},...,p_{N})^{\beta_{N}/\beta_{N^{*}}}}{\int_{X^{N}}\mbox{Per}(x_{1},...,x_{N},p_{1},...,p_{N})^{\beta_{N}/\beta_{N^{*}}}\mu_{0}^{\otimes N}}\mu_{0}^{\otimes N}\otimes\nu^{\otimes N}
\]
The convergence of the expected one-point correlation measures in
the previous theorem yields for any $\beta>0,$ just as before, a
sequence of explicit\emph{ }approximate solutions $\phi_{N,\beta}$
to the real Monge-Ampère equations \ref{eq:weak ma eq with nu} (but
now integrating wrt $\mu_{0}^{\otimes(N-1)}\otimes\nu^{\otimes N}).$
For the case $\beta=0$ we also obtain a variant of Cor \ref{cor:sec cor intro},
which may be formulated as follows:
\begin{cor}
\label{cor:general quenched conv towards transport map}Assume given
two probability measures $\mu_{0}$ and $\mu_{1}$ of the form $\mu_{i}=\rho_{i}1_{X_{i}}dx_{i}$
such that $X_{i}$ is the closure of a bounded domain whose boundary
$\partial X_{i}$ is a null set for Lebesgue measure and assume that
$\rho_{i}$ is bounded from below and above by positive constants
on $X.$ Assume also that $X_{1}$ is convex. Set 
\[
\phi_{N}(x_{1}):=\frac{1}{\beta_{N}^{*}}\int_{X^{N-1}\times P^{N}}\log\mbox{Per}(x_{1},...,x_{N},p_{1},...,p_{N})\mu_{0}^{N-1}\otimes\mu_{1}^{N}-c_{N},
\]
where $c_{N}$ is the normalizing constant ensuring that $\int_{\R^{n}}\phi^{(N)}\mu_{0}=0.$
\textup{Then $\phi_{N}$ converges}, as $N\rightarrow\infty,$\textup{
locally uniformly to the unique convex function $\phi$ }solving the
equation \ref{eq:weak ma eq with nu} for $\beta=0$ and such that
$\nabla\phi$ maps $X_{0}$ almost surjectively onto $X_{1}.$\textup{
Moreover, $T^{(N)}=\nabla\phi^{(N)}$ converges point-wise, in the
interior of $X_{0},$ to the (Hölder continuous) optimal map $T$
for the Monge problem of transporting the probability measure $\mu_{0}$
to $\mu_{1},$ where the optimality is defined with respect to the
cost function $c(x,p)=|x-p|^{2}.$ }
\end{cor}

\subsubsection{Proof of Theorem \ref{thm:conv of one point cor with random alloc}}

As the arguments are similar to the previous ones, we will be rather
brief. Let us first show that 
\begin{equation}
\lim_{N\rightarrow\infty}\E(-\frac{1}{\beta N}\log Z_{N,\beta})=\inf_{\mathcal{\mu\in M}_{1}(X)}\left(C(\mu,\nu)+D_{\mu_{0}}(\mu)/\beta\right)\label{eq:conv of quenced free energies}
\end{equation}
To this end we first observe that $-\frac{1}{\beta N}\log Z_{N,\beta}$
viewed as a function of the quenched variables $(p_{1},...,p_{N})$
is Lipschitz continuous in each coordinate $p_{i}$ with a uniform
Lipschitz constant $L$ (which is proportional to the diameter of
$X).$ Indeed, this follows immediately from the fact that for any
fixed $(x_{1},..,x_{N})$ the Hamiltonian $H^{N}/N$ has the corresponding
Lipschitz property, which in turn follows from the fact that $\partial c(x,p)/\partial p=x$
is uniformly bounded, since we have assumed that $X$ is compact.
Moreover, by the same argument $-\frac{1}{\beta N}\log Z_{N,\beta}$
is uniformly bounded. Now, using (a weak form of) Sanov's theorem
we may replace the integration over $P^{N}$ with the integral over
a ball $B_{\delta}(\nu)$ of a fixed small radius $\delta$ centered
at $\nu$ in the space $\mathcal{M}_{1}(P)$ of all probability measures
on $P.$ Then we pick a sequence $p^{(N)}\in P^{N}$ approximating
$\nu$ in the sense of \ref{eq:weak approx of nu} and in particular
$\delta^{(N)}(p^{(N)})\in B_{\delta}(\nu)$ for $N$ sufficiently
large. The point is that, as explained above, along this sequence
we have 
\begin{equation}
\lim_{N\rightarrow\infty}(-\frac{1}{\beta N}\log Z_{N,\beta})=\inf_{\mathcal{\mu\in M}_{1}(X)}\left(C(\mu,\nu)+D_{\mu_{0}}(\mu)/\beta\right)\label{eq:proof of thm quenched one point}
\end{equation}
Finally, by the uniform Lipschitz estimate and the uniform bound on
$-\frac{1}{\beta N}\log Z_{N,\beta}$ the oscillation of $-\frac{1}{\beta N}\log Z_{N,\beta}$
on $B_{\delta}(\nu)$ is bounded by a uniform constant times $\delta$
and hence the previous convergence implies the convergence in \ref{eq:conv of quenced free energies}
by first letting $N\rightarrow\infty$ and then $\delta\rightarrow0.$

Now, fixing a continuous function $u$ on $X$ we can repeat the previous
argument with $H^{(N)}$ replaced with $H^{(N)}+u$ to get 
\[
\Lambda_{N}[u]:=\E(-\frac{1}{\beta N}\log Z_{N,\beta}[u])\rightarrow\Lambda[u]:=\inf_{\mathcal{\mu\in M}_{1}(X)}\left(C(\mu,\nu)+D_{\mu_{0}}(\mu)/\beta+\int u\mu\right)
\]
Next we observe that the measure $\E(\int_{X^{N-1}}\mu_{\beta_{N}}^{(N)})$
on $X$ represents the differential at $u=0$ of the functional $\Lambda_{N}[u]$
on $C^{0}(X).$ But the latter functional is convex and converges
to $\Lambda$ whose differential at $0$ is the unique minimizer $\mu_{\beta}$
of the functional on $\mathcal{M}_{1}(X)$ appearing in the rhs above
\ref{eq:proof of thm quenched one point}. But then it follows, as
before, by basic convex analysis that $\mu_{\beta}$ represents the
differential of $\Lambda$ at $u=0$ and that $\E(\int_{X^{N-1}}\mu_{\beta_{N}}^{(N)})$
converges to $\mu_{\beta},$ as desired.

\subsection{Proof of Cor \ref{cor:general quenched conv towards transport map}}

Switching the order of integration we can write 
\[
\E(E^{(N)}(\mu))=\frac{1}{N}\int_{X^{N}}\left(\int_{P^{N}}\frac{1}{\beta_{N}^{*}}\log\sum_{\sigma\in S_{N}}e^{\beta_{N}^{*}(x_{1}\cdot p_{\sigma(1)}+\cdots+x_{N}\cdot p_{\sigma(N)})}\nu^{\otimes N}\right)\mu^{\otimes N}
\]
Now we can proceed precisely as in the proof of the previous theorem
by localizing the integration over $P^{N}$ to a ball $B_{\delta}(\nu)$
of a fixed small radius $\delta$ centered at $\nu$ in the space
$\mathcal{M}_{1}(P)$ of all probability measures on $P$ and picking
a sequence $p^{(N)}\in P^{N}$ approximating $\nu$ in the sense of
\ref{eq:weak approx of nu}. The point is that, as explained above,
along this sequence we have $E_{p^{(N)}}^{(N)}(\mu)\rightarrow C(\mu,\nu).$
Finally, by the Lipschitz uniform estimate for $H^{(N)}/N$ the oscillation
of $H^{(N)}/N$ on $B_{\delta}(\nu)$ is bounded by a uniform constant
times $\delta$ and hence 
\[
\lim_{N\rightarrow\infty}\E(E^{(N)}(\mu))=\lim_{N\rightarrow\infty}E_{p^{(N)}}^{(N)}(\mu)=C(\mu,\nu)
\]
 Now the proof is concluded by differentiating with respect to $\mu,$
precisely as in the proof of Cor \ref{cor:sec cor intro}. Note that
the convexity of $P$ is crucial to get the Hölder regularity of the
transport map, as explained in \cite{ca-2}.

\section{\label{sec:Outlook}Outlook}

\subsection{\label{sub:Relation-to-the}Relation to the complex setting, determinantal
point processes and toric varieties }

\subsubsection{The toric setting}

In this section we come back to the original setting where $\nu=\lambda_{P}$.
Consider the complex torus $\C^{*n}$ and denote by $T^{n}$ the corresponding
real unit-torus in $\C^{*n},$ which acts, in the standard way, holomorphically
on $\C^{*n}.$ Denote by $\mbox{Log }$be the map 
\[
\mbox{Log}:\,\,\C^{*n}\rightarrow\R^{n}
\]
 from $\C^{*n}$ to $\R^{n}$ defined by $x:=\mbox{Log}(z):=x,$ where
$x$ is the vector whose $j$ th coordinate is the log of the squared
absolute value of the $j$ th coordinate of $z.$ The fibers are thus
the orbits of the real torus $T^{n}$ on $\C^{*n}.$ The definition
is made so that, if $p\in\Z^{n},$ then $|z^{p}|^{2}=e^{p\cdot x}$
in multiindex notation. 

Denote by $\Delta^{(N)}(z_{1},...,z_{N})$ the \emph{Vandermonde determinant}
on $(\C^{*n})^{N}$ determined by the convex body $P,$ i.e. 
\[
\Delta^{(N)}(z_{1},...,z_{N})=\det(z_{i}^{p_{j}}),
\]
where $p_{j}$ ranges over the $N$ lattice points in $kP.$ Fixing
a measure $\tilde{\mu}_{0}$ on $\C^{*n}$ and a continuous function
$\tilde{\phi}_{0}(z)$ of suitable growth one obtains, for any sequence
$\beta_{N}$ of positive numbers, a probability measure $\tilde{\mu}^{(N)}$
on $(\C^{*n})^{N}$ by setting 
\begin{equation}
\tilde{\mu}_{\beta_{N}}^{(N)}:=\frac{1}{Z_{N,\beta_{N}}}|\Delta^{(N)}(z_{1},...,z_{N})|^{2\beta_{N}/k}e^{-k(\tilde{\phi}_{0}(z_{1})+\cdots+\tilde{\phi}_{0}(z_{N}))}\tilde{\mu}_{0}^{\otimes N}\label{eq:gibbs meas in complex toric setting}
\end{equation}
For $\beta_{N}=k$ this is defines a\emph{ determinantal point process
}(see \cite{h-k-p}\emph{ }for general properties of such processes
and \cite{be3} for large deviation results for these particular determinantal
processes). The relation to the present paper stems from the simple
observation that if the background data $(\tilde{\mu}_{0},\tilde{\phi}_{0})$
is $T^{n}-$invariant, then the push-forward of the corresponding
determinantal point process is a\emph{ permanental} point process.
More precisely, in the case $\beta_{N}=k,$ the push-forward of the
corresponding probability measure $\tilde{\mu}^{(N)}$ is precisely
the permantental probability measure $\tilde{\mu}^{(N)}$ studied
in the present paper, determined by the weighted measure $(\mu_{0},\phi_{0}),$
where $\tilde{\phi}_{0}=\mbox{Log}^{*}\phi_{0}$ and $\mbox{Log}_{*}\tilde{\mu}_{0}=\mu_{0}$
(i.e. abusing notation slightly $\phi(x)=\tilde{\phi}(z)$ and $\tilde{\mu}=\mu\wedge d\theta,$
where $d\theta$ denotes the invariant probability measure on $T^{n}).$
To see this just note that expanding $\Delta^{(N)}(z_{1},...,z_{N})$
as an alternating sum over the permutations $\sigma$ in $S_{N}$
and using Parseval's formula for $x$ fixed to carry out the integration
over the corresponding torus fiber $\mbox{Log \ensuremath{^{-1}(\{x\})}}$
immediately gives 
\begin{equation}
(\mbox{Log })_{*}\left(\Delta^{(N)}(z_{1},...,z_{N})(\mu_{0}\otimes d\theta)^{\otimes N}\right)=\mbox{Per}(x_{1},...,x_{N})(\mu_{0})^{\otimes N}.\label{eq:push det is per}
\end{equation}
Moreover, it is also interesting to see that there is a pluripotential
analog of this deteminantal/permanental correspondence: denoting by
$E(\tilde{\mu})$ the\emph{ pluricomplex energy} of a probability
measure $\tilde{\mu}$ on $\C^{*n}$ (defined with respect to the
reference weight $\mbox{Log}^{*}\phi_{P})$ one gets, if $\tilde{\mu}$
is $T^{n}-$invariant, that 
\begin{equation}
E(\tilde{\mu})=C(\mu),\label{eq:pluri energy as cost}
\end{equation}
 where, as before, $C(\mu)$ is the optimal cost functional functional
defined with respect to the cost function $c(x,p):=-x\cdot p$ and
the target measure $\lambda_{P}.$ This follows immediately from Prop
\ref{prop:energy as cost} combined with the essentially well-known
fact that, when $\tilde{\mu}$ is $T^{n}-$invariant, $E(\tilde{\mu})$
can be expressed in terms of the Legendre transform of the Monge-Ampère
potential of $\tilde{\mu}.$ The key point is the basic fact that
if $\tilde{\phi}$ is a $T^{n}-$invariant plurisubharmonic function
(i.e. $\partial\bar{\partial}\tilde{\phi}\geq0)$ then $\phi$ is
convex and 
\[
\mbox{Log }_{*}MA_{\C}(\tilde{\phi})=MA(\phi),
\]
 where $MA_{\C}$ denotes the\emph{ complex} Monge-Ampère operator,
i.e. 
\[
MA_{\C}(\psi):=(\frac{i}{2\pi}\partial\bar{\partial}\psi)^{n}/n!\left(=c_{n}\det(\frac{\partial\psi}{\partial z_{i}\partial\bar{z_{j}}})\right)
\]
Alternatively, the relation \ref{eq:pluri energy as cost} follows
from the correspondence \ref{eq:push det is per} by combining the
large deviation principle for the determinantal point processes in
\cite{be3}, where $E$ appears as the rate functional, applied to
the toric case, with the large deviation principle for the corresponding
permanental point-process proved in the present paper. Strictly speaking
the results in \cite{be3} only apply when $P$ is a rational polytope,
but the proofs are essentially the same in the general case (compare
the setting in \cite{be0}). The point is that when $P$ is a rational
polytope it defines a toric variety $X_{P}$ with an ample line bundle
$L_{P}$ to which the results in \cite{be3} can be applied. Briefly,
the toric variety $X_{P},$ which is an equivariant compactification
of $\C^{*n},$ may be defined as the projective algebraic variety
obtained as the closure in complex projective space $\P^{N}$ of the
affine algebraic variety in $\C^{N}$ defined by the image of the
map 
\[
\C^{*n}\rightarrow\C^{N},\,\,\,\, z\mapsto(z^{p_{1}},...,z^{p_{N}}),
\]
 for $k$ sufficiently large and $L_{P}$ is the restriction of the
hyperplane line bundle on $\P^{N}$ (see \cite{ber-ber} and references
therein). 

It should be stressed that, unless $\beta_{N}=k,$ the push-forward
under the map Log of the probability measure $\tilde{\mu}_{\beta_{N}}^{(N)}$
on $(\C^{*n})^{N}$ is \emph{not} equal to the corresponding probability
measure $\mu_{\beta_{N}}^{(N)}$ on $(\R^{n})^{N}.$ Still, one would
expect that this is true in an asymptotic sense, as $N\rightarrow\infty.$

\subsubsection{The general complex geometric setting and Kähler-Einstein geometry}

The general complex geometric setting of Gibbs measures of the form
\ref{eq:gibbs meas in complex toric setting} and the relation to
the Kähler-Einstein geometry will be studied in detail elsewhere \cite{be-4}
(for outlines see \cite{be,be-3} ). Here we will just give a brief
impressionistic view of the setting. The general geometric background
data consists of a pair $(\mu_{0},\phi_{0})$ where $\mu_{0}$ is
a measure on the $n-$dimensional complex manifold $X$ and $\phi_{0}$
is a metric on an ample line bundle $L\rightarrow X$ (more precisely,
we will denote by $\phi_{0}$ the collection of local functions such
that $e^{-\phi_{0}}$ represents the metric with respect to given
local trivializations of $L$). To this data one may associate a sequence
of Gibbs measure of the form \ref{eq:gibbs meas in complex toric setting},
but with $\Delta^{(N_{k})}$ replaced with any generator of the deteminant
line $\Lambda^{N}H^{0}(X,L^{\otimes k}),$ where $H^{0}(X,L^{\otimes k})$
denotes the $N-$dimensional space of global holomorphic sections
with values in the $k$ th tensor power of $L.$ When $\beta_{N}=\beta$
the corresponding mean field type equations are then of the form 
\begin{equation}
MA_{\C}(\phi)=e^{\beta(\phi-\phi_{0})}\mu_{0}\label{eq:ma eq on complex manifold with background}
\end{equation}
 for a positively curved metric $\phi$ on the line bundle $L.$ The
relation to Kähler-Einstein geometry stems from the fact when $L$
is taken as the\emph{ canonical line bundl}e $K_{X}:=\Lambda^{n}(T^{*}X)$
any metric $\phi_{0}$ determines a measure $\mu_{0}=e^{+\phi_{0}}dz\wedge d\bar{z}$
and the equation \ref{eq:ma eq on complex manifold with background}
is then\emph{ instrically} defined for $\beta=1$ (i.e. independent
of $\phi_{0})$. In fact, as is well-known the equation is then equivalent
to the Einstein equation with cosmological constant $\Lambda=-1$
for the Kähler metric $\omega:=\frac{i}{2\pi}\partial\bar{\partial}\phi$
on $X,$ i.e. the equation 
\begin{equation}
\mbox{\ensuremath{\mbox{Ric }}}\omega=\Lambda\omega,\label{eq:einstein eq with cosm cst}
\end{equation}
 where $\mbox{\ensuremath{\mbox{Ric }}}\omega$ denotes the Ricci
curvature of $\omega.$ Moreover, the corresponding Gibbs measure
is then also intrinsically defined by $X.$ Similarly, if $X$ is
a \emph{Fano manifold, }i.e. the anti-canonical line bundle, $K_{X}^{-1}:=\Lambda^{n}(TX)$
is ample, then we can take $L=K_{X}^{-1}.$ Any metric $\phi_{0}$
on $K_{X}^{-1}$ determines a well-defined measure $\mu_{0}=e^{-\phi_{0}}dz\wedge d\bar{z}$
on $X$ and for $\beta=-1$ the corresponding equation \ref{eq:ma eq on complex manifold with background}
coincides with the Einstein equation for $\omega$ with cosmological
constant $\Lambda=+1.$ However, in this setting the corresponding
Gibbs measure will \emph{not} be well-defined in general since the
partition function may diverge (the reason is that the corresponding
integrand is then locally of the form $1/|f_{k}(z_{1},...,z_{N})|^{2/k}$
for a holomophic function $f_{k}$ and the integrability properites
are thus reflected in the singularities of the hypersurface cut out
by $f_{k}).$ It is then natural to define a statistical mechanical
notion of stability of a Fano manfold $X$ called\emph{ Gibbs stability}
by demanding that the partition function be finite for $k$ sufficently
large \cite{be,be-3}. This should be thought of as a probabilistic
version of other notions of algebro-geometric stability, such as K-stability,
appearing in Kähler-Einstein geometry. Interestingly, Gibbs stability
admits a purely algebro-geometric interpretation saying that the anti-canonical
incidence divisor in $X^{N_{k}}$ has, for $k$ sufficently large,
mild singularities in the sense of the Minimal Model Program (or more
precisely, Kawamata log terminal singularities). There are also variations
of the notions of Gibbs stability, which for example, are needed when
$X$ admits non-trivial holomorphic vector fields (since $X$ will
never be Gibbs stable then). In the following sections we will out-line
some concrete relations between the general complex geometric setting
and the present one.

\subsection{\label{sub:A-general-LDP}A general LDP for Gibbs measures and Coulomb
type gases}

Let us consider the following general setting: $X$ is a topological
space equipped with a (Borel) measure $\mu_{0}$ and $H^{(N)}$ is
a sequence of symmetric functions on $X^{N}.$ We also assume given
a sequence $\beta_{N}\rightarrow\infty$ such that $Z_{N,\beta_{N}}:=\int_{X^{N}}e^{-\beta_{N}H^{(N)}}\mu_{0}^{\otimes N}$
is finite for any $N.$ Then the corresponding Gibbs measures $\mu_{\beta_{N}}^{(N)}$
are well-defined. For simplicity we will assume that $X$ is compact
and that $\mu_{0}$ has finite mass and may hence (up to a harmless
scaling) be assumed to be a probability measure. Let us also assume
that the assumptions in the Gärtner-Ellis theorem hold, i.e. that
there exists a Gateaux differentiable functional $\mathcal{F}(u)$
on $C^{0}(X)$ such that 
\[
-\lim_{N\rightarrow\infty}\frac{1}{\beta_{N}N}\log Z_{N,\beta_{N}}[u]=\mathcal{F}(u)
\]
 By the Gärtner-Ellis theorem it then follows that law of the empirical
measure, i.e. $(\delta_{N})_{*}(\mu_{\beta_{N}}^{(N)}),$ satisfies
a LDP on $\mathcal{M}_{1}(X)$ with speed $\beta N$ and rate functional
equal to $E$ (up to an additive normalizing constant), where, as
before, the functional $E(\mu)$ denotes the Legendre transform of
$\mathcal{F},$ i.e. $E=\mathcal{F}^{*}.$

Replacing the sequence $\beta_{N}$ with a fixed positive number $\beta$
it is natural to ask under what conditions the corresponding Gibbs
measures $\mu_{\beta}^{(N)}$ satisfy a LDP with rate functional equal
to $E+D_{\mu_{0}}/\beta$ (up to an additive normalizing constant)?
By Theorem \ref{thm:interacting part intro} it is enough to assume
equicontinuity of $H^{(N)},$ but the proof given above actually applies
in a considerable more general setting and, loosely speaking, shows
that the result holds as long as a certain chaoticity property holds
(a more direct proof under the equicontinuity setting is given in
section \ref{sub:An-alternative-proof} below). In order to formulate
this properly let us denote by $\mu_{\beta_{N}}^{(N),u}$ the ``tilted''
Gibbs measures obtained by replacing $H^{(N)}/N$ with $H^{(N)}/N+u.$
By the Gärtner-Ellis theorem $\mu_{\beta_{N}}^{(N),u}$ satisfies
a LDP principle for any fixed $u$ and it also follows that the $j$
th marginal of $\mu_{\beta_{N}}^{(N),u}$ converges to $(\mu_{u})^{\otimes j}$
where
\[
\mu_{u}:=d\mathcal{F}_{|u}.
\]
 In the terminology of Kac this says that the sequence $\mu_{\beta_{N}}^{(N),u}$
is \emph{$\mu_{u}-$ chaotic} (see \cite{sn} and references therein).
Now the main extra property that is needed to deduce the LDP for a
fixed $\beta$ is that the whole measure $\mu_{\beta_{N}}^{(N),u}$
is sufficiently close to the corresponding product measure $(\mu_{u})^{\otimes N}$
in an entropic sense. 
\begin{thm}
\label{thm:general ldp}Assume given data $(X,\mu_{0},H^{(N)},\beta_{N})$
as above, such that the corresponding functional $\mathcal{F}(u)$
is Gateaux differentiable. Assume that
\begin{itemize}
\item $\inf_{X^{N}}\frac{H^{(N)}+u}{N}=-\frac{1}{\beta_{N}N}\log Z_{N,\beta_{N}}[u]+o(1)$
for any fixed $u\in C^{0}(X)$
\item $\lim_{N\rightarrow\infty}\frac{D\left((\mu_{u})^{\otimes N},\mu_{\beta_{N}}^{(N),u}\right)}{N\beta_{N}}=0$
for any ``good'' function $u$ in $C^{0}(X)$ i.e. such that $D_{\mu_{0}}(\mu_{u})<\infty$
and $E_{\mu_{0}}(\mu_{u})<\infty$
\item There exists some probability measure $\mu$ such that $D_{\mu_{0}}(\mu)<\infty$
and $E(\mu)<\infty$ and any such measure may be written as a weak
limit of measures $\mu_{u_{j}}$ with $u_{j}$ good such that the
functionals $D_{\mu_{0}}$ and $E$ are continuous along $\mu_{u_{j}}.$
\end{itemize}
Then $\mu_{\beta}^{(N),u}$ satisfies an LDP with speed $\beta N$
and rate functional $F_{\beta}:=E+D_{\mu_{0}}/\beta$ (up to an additive
constant) and 
\begin{equation}
-\lim_{N\rightarrow\infty}\frac{1}{\beta_{N}N}\log Z_{N,\beta_{N}}=\inf_{\mu\in\mathcal{M}_{1}(X)}F_{\beta}(\mu)\label{eq:asympt of part function in thm general ldp}
\end{equation}
\end{thm}
\begin{proof}
The proof is essentially contained in the previous arguments, so we
will only recall the main points. First, fixing a probability measure
$\mu$ on $X$ and a continuous function $u$ we we may, in a similar
manner as in the proof of Theorem \ref{thm:existence of mean energ},
rewrite 
\[
\int_{X^{N}}\frac{H^{(N)}}{N}\mu^{\otimes N}=\left(-\frac{1}{\beta_{N}N}\log Z_{N,\beta_{N}}[u]-\int_{X}u\mu\right)+\frac{D\left((\mu)^{\otimes N},\mu_{\beta_{N}}^{(N),u}\right)}{N\beta_{N}}-\frac{1}{\beta_{N}}D_{\mu_{0}}(\mu)
\]
(assuming that all terms are finite). In particular, if $u$ is good
and $\mu=\mu_{u}$ the rhs converges to $E(\mu),$ as $N\rightarrow\infty.$
We can then proceed exactly as before, using the Gibbs variational
principle to get, for any fixed $\mu,$ by first taking the infimum
over all good $u,$ that \emph{
\[
F_{\beta}^{(N)}(\mu_{\beta}^{(N)})\leq\inf_{u}F_{\beta}(\mu_{u})=\inf_{\mathcal{M}_{1}(X)}F_{\beta}(\mu_{u}),
\]
} where we have used the third assumption in the last equality. As
for the lower bound it is proved exactly as before: if $\mu_{*}$
is a limit point of the first marginal $\mu_{1}^{(N)},$ then, for
any $u\in C^{0}(X)$ we get \emph{
\[
\lim_{N\rightarrow\infty}\inf_{X^{N}}\left(\frac{H^{(N)}+u}{N}-\int u\mu_{*}\right)\leq\liminf_{N\rightarrow\infty}F_{\beta}^{(N)}(\mu_{\beta}^{(N)})
\]
} and we can then conclude that the asymptotics in formula \ref{eq:asympt of part function in thm general ldp}
hold. The large deviation property then follows precisely as before
by replacing $H^{(N)}$ with $H^{(N)}+u.$ 
\end{proof}
The reason that we have invoked the approximation property is that,
in general, the solution of the corresponding mean field type equations
may not be continuous. Anyway, in many cases continuity of the solution
is guaranteed and then the approximation assumption appearing in the
third point above is not needed. It should also be pointed out that
the first assumption may be removed if one instead defines the functional
$\mathcal{F}$ in terms of the infimum over $X^{N}$ (as in formula
\ref{eq:def of f of u as limit of ham}). This is particularly useful
when $\mu_{0}$ is very irregular. To illustrate this we state the
following general result about Coulomb gases whose complete proof
will appear elsewhere.
\begin{thm}
\cite{be-4}Let $\mu_{0}$ be a (Borel) measure on $\C$ and $\phi_{0}$
a continuous function of super logarithmic growth (i.e. $\phi_{0}\geq\log((1+|z|^{2})-C)$
defined on the support of $\mu_{0},$ such that $\int e^{-\beta\phi}\mu_{0}<\infty.$
Set 
\[
H^{(N)}(z_{1},...,z_{N}):=-\frac{1}{N}\sum_{1\leq i\neq j\leq N}\log|z_{i}-z_{j}|+\sum\phi(z_{i})
\]
and 
\[
E(\mu):=-\int\log|x-y|\mu(x)\otimes\mu(y)+\int\phi_{0}\mu
\]
 Then, for any positive number $\beta,$ the law of the empirical
measure of the corresponding Gibbs measure $\mu_{\beta}^{(N)}$satisfies
a LDP with rate functional $E+D_{\mu_{0}}(\mu)/\beta$ (up to an additive
constant) as long as the latter functional is not identically equal
to infinity. 
\end{thm}
We recall that corresponding Gibbs measure $\mu_{\beta}^{(N)}$ may
in this setting be written in terms of the Vandermonde determinant
for polynomials of degree $\leq k:=N-1$ (compare section \ref{sub:Relation-to-the}
for the relation to permanents): 
\[
\mu_{\beta}^{(N)}=\frac{1}{Z_{N,\beta}}|\Delta^{(N)}(z_{1},...,z_{N})|^{2\beta/k}e^{-(\phi_{0}(z_{1})+\cdots+\phi_{0}(z_{N}))}\mu_{0}^{\otimes N}
\]
(this is a special case of the complex geometric framework outlined
i the previous section, where $X$ is the complex projective line
(i.e. the Riemann sphere) $\P^{1},$ viewed as the one-point compactification
of $\C,$ and $L$ is the hyperplane line bundle $\mathcal{O}(1)\rightarrow\P^{1}).$
In this setting the corresponding mean field equations on $\C$ may,
in complex notation, be written as 
\begin{equation}
\frac{i}{2\pi}\partial\bar{\partial}\phi=e^{\beta(\phi-\phi_{0})}\mu_{0},\label{eq:mean field eq in complex plane}
\end{equation}
 where the function $\phi$ (which is automatically subharmonic) satisfies
the normalization condition $\frac{i}{2\pi}\int_{\C}\partial\bar{\partial}\phi=1$
(which for example holds if $\phi=\log|z|^{2}+O(1)$ as $|z|\rightarrow\infty).$ 

In the case when $\mu_{0}$ is the Lebesgue measure supported on a
bounded domain $\Omega$ in the plane, the previous theorem was first
shown in \cite{clmp,k}. But the main point here is that the method
of proof indicated above applies to any measure $\mu_{0}$ with the
property that the corresponding free energy functional $F_{\beta}$
is not identically equal to infinity (which is equivalent to the existence
of a minimizer). For example, any measure $\mu_{0}$ not charging
polar sets will do. In particular, the assumption about the Bernstein-Markov
property of $\mu_{0}$ which appears in the case $\beta=\infty$ \cite{be3,b-l}
is not needed for $\beta$ finite. 

However, in the case of $\beta<0$ stronger assumptions on $\mu_{0}$
are needed. It turns out that an essentially optimal regularity assumption
is that there are positive constants $C$ and $d$ such that the measure
$\mu_{0}$ satisfies 
\[
\mu_{0}(B_{r})\leq Cr^{d},
\]
 for $r$ sufficiently small, for every Euclidean ball of radius $r.$
Under this assumption the corresponding Gibbs measure $\mu_{\beta}^{(N)}$
is well-defined for any $\beta>-d$ and, after passing to a subsequence,
the law of the corresponding empirical measure converges weakly to
a measure concentrated on the set of minimizers of $F_{\beta}$ \cite{be-4}.
In particular situations the uniqueness of such minimizers can be
ensured. See \cite{clmp,k} for the case when $\mu_{0}$ is the Lebesgue
measure supported on a bounded domain $\Omega$ in the plane. The
authors main motivation for studying this setting comes from the relation
to the study of Kähler-Einstein metrics on complex algebraic varieties,
i.e. Kähler metrics with constant Ricci curvature. For example, the
general case referred to above in particular applies to the following
setting which corresponds to the complex geometric setting of of conical
metrics on the Riemann sphere with positive constant curvature:
\begin{thm}
\cite{be-4}Consider $m$ points $p_{1},...p_{m}$ in $\C$ with weights
$w_{1},...,w_{m}\in]0,1[$ such that $d:=2-(w_{1}+\cdots+w_{m})>0$
and set $\phi_{0}=0$ and 
\[
\mu_{0}:=\frac{1}{|z_{1}-p_{1}|^{2w_{1}}\cdots|z_{m}-p_{m}|^{2w_{m}}}\frac{i}{2}dz\wedge d\bar{z}
\]
 Then the following is equivalent: 
\begin{itemize}
\item The equation \ref{eq:mean field eq in complex plane} admits a unique
solution for $\beta=-d$ 
\item the corresponding partition functions $Z_{N,-d}$ are finite for $N>>1$
\end{itemize}
Moreover, in that case the corresponding empirical measure converges
in probability to $\mu_{-d}:=\frac{i}{2\pi}\partial\bar{\partial}\phi_{-d}$
where $\phi_{-d}$ is the solution of equation \ref{eq:mean field eq in complex plane}.
\end{thm}
The relation to complex geometry comes from the well-known fact that
$\phi$ is a solution to the equation appearing in the previous theorem
iff $\omega:=\frac{i}{2\pi}\partial\bar{\partial}\phi$ defines a
Kähler metric on the Riemann sphere $X$ with constant positive curvature
and conical singularities encoded by the effective $\R-$divisor $E:=p_{1}w_{1}+\cdots+p_{m}w_{m},$
i.e. $\omega$ is a conical Kähler-Einstein metric. From this point
of view the previous theorem can be formulated as saying that a one
dimensional \emph{log Fano variety }$(X,E)$ admits a unique conical
Kähler-Einstein metric iff it is \emph{Gibbs stable} in the sense
of \cite{be-3}. This and the higher dimensional setting will be studied
in detail in \cite{be-4}. For the moment we just point out that it
is well-known that a necesserary condition for uniqueness in the previous
theorem is that there are at most \emph{two} points $p_{i}$ or more
precisely, either (i) no points or (ii) precisely two points (and
in the latter case $w_{1}=w_{2}).$ In the first case uniqueness indeed
fails since the equations are invariant under all biholomoprhic maps
of the Riemann sphere (i.e. the Möbius group) and in the second case
invariance holds under the standard $\C^{*}$ action which fixes the
two points. This symmetry can also be seen to be responsable for the
fact that the corresponding partitionfunctions $Z_{N}$ then diverge.
However, there is a way to break the symmetry in order to restore
uniqueness and finiteness in the previous theorem. To this end one
fixes a subharmonic function $\phi_{0}$ on $\C$ with logarithmic
growth and considers the correponding Coulomb gas defined by the weight
$\phi_{0}$ and $\mu_{0}=e^{-d\phi_{0}}\frac{1}{|z_{1}-p_{1}|^{2w_{1}}\cdots|z_{m}-p_{m}|^{2w_{m}}}\frac{i}{2}dz\wedge d\bar{z}.$
Then it can be proved that there is a unique solution $\phi_{\beta}$
to the corresponding mean field type equations for $\beta>-d.$ Moreover,
the corresponding partition functions $Z_{N,-\beta}$ are then finite
and the empirical measure converges in probability towards corresponding
measure $\frac{i}{2\pi}\partial\bar{\partial}\phi_{\beta}.$ In the
particular case when $\phi_{0}$ has circular symmetry this is closely
related to the one-dimensional real Monge-Ampère equations. In fact,
a similar phenomena persists in higher dimensions under toric symmetry.
This is the subject of the the next section where we will outline
the relation between toric Kähler-Einstein metrics and the previous
probabilistic setting of permanents and the real Monge-Ampère equation.

\subsection{\label{sub:K=0000E4hler-Einstein-metrics,-negativ}Toric Kähler-Einstein
metrics, negative $\beta$ and phase transitions }

First consider the following special case of the setting of weighted
measures $(\mu_{0},\phi_{0})$ in $\R^{n}$ (section \ref{sub:Weigted-sets-and}):
given a weight function $\phi_{0}$ on $\R^{n}$ we take the measure
$\mu_{0}$ to be given by $\mu_{\phi_{0}}:=e^{-\phi_{0}}dx.$ Setting
$\gamma=-\beta$ the corresponding Monge-Ampère equation \ref{eq:ma eq weak intro}
may then be written as 
\begin{equation}
MA(\phi)=e^{-(\gamma\phi+(1-\gamma)\phi_{0})}dx,\label{eq:monge-ampere with gamma}
\end{equation}
 which in turn can be written as a \emph{twisted Kähler-Einstein equation}
on the complex torus $\C^{*n}.$ Indeed, let $\mbox{Log }$be the
map from $\C^{*n}$ to $\R^{n}$ defined in section \ref{sub:Relation-to-the}
and set $\varphi:=\mbox{Log}^{*}\phi.$ Then $\omega:=\frac{i}{2\pi}\partial\bar{\partial}\varphi$
defines a Kähler metric on $\C^{*n}$ which satisfies 
\begin{equation}
\mbox{\ensuremath{\mbox{Ric }}}\omega=\gamma\omega+(1-\gamma)\omega_{0},\label{eq:aubin}
\end{equation}
 where $\mbox{\ensuremath{\mbox{Ric }}}\omega$ is the Ricci curvature
of the Kähler metric $\omega,$ represented as two-form. In particular,
for $\gamma=1$ a solution $\omega$ is a bona fide Kähler-Einstein
metric, i.e. a Kähler metric with constant (positive) Ricci curvature
and the corresponding convex function $\phi(x)$ then satisfies the
$\phi_{0}-$independent equation 
\begin{equation}
MA(\phi)=e^{-\phi}dx\label{eq:toric k-e equation}
\end{equation}
When the convex body $P$ is a polytope (containing zero in its interior)
$\omega$ extends to a (singular) Kähler-Einstein metric on the corresponding
toric variety $X_{P}$ compactifying $\C^{*n}$ (see \cite{ber-ber}
and references therein). The most studied case is when $P$ is a \emph{reflexive
Delzant polytope, }which equivalently means that $X_{P}$ is a \emph{Fano
manifold. }Then the equation \ref{eq:aubin} coincides with Aubin's
continuity equation, designed by Aubin to prove the existence of a
Kähler-Einstein metric by deforming $\gamma$ from $\gamma=0$ to
$\gamma=1.$ The existence of solutions for $\gamma$ a sufficiently
small positive number was shown by Aubin. However, as is well-known
there are in general obstructions to the existence of Kähler-Einstein
metric with positive Ricci curvature and according to the fundamental
Yau-Tian-Donaldson conjecture the existence of a Kähler-Einstein metrics
is equivalent to an algebro-geometric notion of stability, called
K-stability. 

Here we will only briefly explain the relation between Aubin's continuity
equation and the probabilistic framework as developed in previous
sections. To this end we assume that the given function $\phi_{0}$
is in $\mathcal{P}_{+}(\R^{n}).$ As shown in \cite{ber-ber} (generalizing
the seminal result of Wang-Zhu concerning the smooth Fano case) the
equation \ref{eq:monge-ampere with gamma} then admits a solution
iff $0$ is the barycenter in $P.$ More precisely, denoting by $R$
the sup over all $\gamma\in]0,1]$ such that the equation admits a
solution in $\mathcal{P}_{+}(\R^{n}),$ it was shown in \cite{ber-ber}
that $R$ is given by the following formula: 
\begin{equation}
R:=\frac{\left\Vert q\right\Vert }{\left\Vert q-b\right\Vert },\label{eq:inv R of conv bod}
\end{equation}
 where $q$ is the point in $\partial P$ where the line segment starting
at $b$ and passing through $0$ meets $\partial P$ (this is a generalization
of a result of Li concerning the smooth Fano case). Moreover, the
corresponding solution is unique for $\gamma<1.$ Interestingly, the
free energy functional $F_{\beta}$ may be identified with \emph{Mabuchi's
K-energy functional} in this setting (or rather its twisted version,
compare \cite{ber-ber,berm2}). To see the connection to the random
point processes considered in the previous section we note that the
processes (i.e. the Gibbs measures) in question are still defined
for negative $\beta$ (i.e. positive $\gamma)$ as long as the corresponding
partition function, which may be written as 
\[
Z_{N,\beta}:=\int_{(\R^{n})^{N}}\left(\mbox{per}(x_{1},...,x_{N})\right)^{-\gamma/k}(e^{-(1-\gamma)\phi_{0}}dx)^{\otimes N},
\]
 is finite. Now the point is that it can be shown that for any fixed
$\gamma<R$ the corresponding partition function $Z_{N,-\gamma}$
is indeed finite for $N$ sufficiently large; more precisely $Z_{N,-\gamma}\leq C_{\gamma}^{N}$
for some constant $C_{\gamma}$ (this can be seen as a toric variant
of the notion of Gibbs stability introduced in \cite{be-4}). Hence,
in the light of Theorem \ref{thm:main for perman intro} one would
then expect that the corresponding empirical measures converge in
probability, as $N\rightarrow\infty,$ towards $\mu_{\beta}:=MA(\phi_{\beta})$
where $\phi_{\beta}$ is the unique solution of the equation \ref{eq:monge-ampere with gamma},
for $\gamma=-\beta.$ However, to extend the previous arguments to
the present setting one needs to handle two issues. First, since we
are assuming that $\phi_{0}$ is in $\mathcal{P}_{+}(\R^{n})$ the
growth assumption on the weight functions (section \ref{sub:Weigted-sets-and})
does not hold any more, but this is only a minor technical point.
More seriously, since $\beta<0$ the argument for the lower bound
in the proof of Theorem \ref{thm:conv of one-point cor measur} is
not valid anymore as it stands. However, this problem can be circumvented
using the Hewitt-Sanders decomposition theorem and the sub-additivity
of the entropy, as in \cite{m-s} (it is also important to know that
$\mu_{\beta}$ is still the \emph{unique} minimizer of the free energy
functional $F_{\beta},$ which is indeed the case \cite{ber-ber}).
Another useful fact is that the first correlation measures are of
the form $e^{-(\gamma\phi_{k}(x)+(1-\gamma)\phi_{0}(x)}dx,$ where
$\phi_{k}(x)$ is in $\mathcal{P}_{+}(\R^{n}),$ as follows immediately
from the Prekopa inequality. This and further relations to the Kähler-Einstein
problem are deferred to \cite{be-4}. Here we will only summarize
the corresponding main results:
\begin{thm}
\cite{be-4}Let $P$ be a convex body containing $0$ in its interior
and $\phi_{0}$ a convex function on $\R^{n}$ such that $\phi_{0}-\phi_{P}$
is bounded. Then, for any $\gamma<R$ the corresponding Gibbs measure
$\mu_{-\gamma}^{(N)}$ is well-defined and the law of its empirical
measure converges in probability, when $N\rightarrow\infty,$ towards
$\mu_{\gamma}:=MA(\phi_{\gamma}),$ where $\phi_{\gamma}$ is a solution
to the equation \ref{eq:monge-ampere with gamma}. 
\end{thm}
One subtle feature of this setting is the presence of \emph{translational
symmetry} at the critical value $\gamma=1$ (assuming that the barycenter
of $P$ vanishes) and the way that it is broken by introducing a weight
$\phi_{0}.$ The point is that for $\gamma<1$ the corresponding mean
field type equations \ref{eq:monge-ampere with gamma} have a unique
solution, while for $\gamma=1,$ there is an $n-$dimensional space
of solutions. This is due to the fact that, in the latter case, the
equations are invariant under the action of $\R^{n}$ by translations.
On the other hand, introducing a weight $\phi_{0}$ brakes this symmetry
and it turns out that the corresponding solutions $\phi_{\gamma}$
tend, when $\gamma\rightarrow1,$ to a particular solution $\phi$
of the equation \ref{eq:toric k-e equation} depending on the choice
of $\phi_{0}.$ The most transparent case is when $P$ as well as
$\phi_{0}$ are symmetric around the origin. i.e. $-P=P$ and $\phi_{0}(-x)=\phi_{0}(x).$
Then we have the following
\begin{cor}
Suppose that the convex body $P$ is symmetric and that $\phi_{0}$
is also symmetric with respect to the origin, i.e. $\phi_{0}(-x)=\phi_{0}(x).$
Then, 
\[
\phi_{\gamma}^{(N)}(x):=-\frac{1}{\gamma}\log\int_{(\R^{n}){}^{N-1}}\frac{1}{Z_{N,-\gamma}}(\mbox{Per}(x,x_{2},...,x_{N}))^{-\gamma/k}(e^{-(1-\gamma)\phi_{0}}dx)^{\otimes N-1},
\]
 converges point-wise, in the double limit where first $N\rightarrow\infty$
and then $\gamma\rightarrow1,$ to the unique solution $\phi$ of
the equation \ref{eq:toric k-e equation} satisfying $\phi(-x)=\phi(x).$
\end{cor}
In the general case, it can be shown that $\phi_{\gamma}$ converges
to the unique solution whose Monge-Ampere measure minimizes the associated
energy functional $E_{\phi_{0}}$ on the solution space. Interestingly,
the $\R^{n}-$symmetry is also responsible for the fact that at the
critical value $\gamma=1$ the corresponding random point processes
are not well-defined (for any $N).$ Indeed, $\mbox{Per}(x_{1},x_{2},...,x_{N})$
is invariant under the diagonal action of $\R^{n}$ and hence the
corresponding partition function $Z_{N,-1}$ diverges.

Let us finally point out that, in general, the critical value $\gamma=R$
can be interpreted as a \emph{second order} \emph{phase transition}
(compare the discussion in the end of \cite{k}, which turns out to
be related to the simplest case of the present setting, namely when
$n=1$ and $P=[-1,1]).$

\subsection{Langevin dynamics}

In this section we will briefly comment on a dynamical version of
Theorem \ref{thm:main for perman intro}, which, from the point of
view of statistical mechanics corresponds to the relaxation to equilibrium
of the corresponding system (often referred to as Langevin dynamics).
It can be seen as a fully non-linear version (when $n>1)$ of McKean's
interacting diffusions \cite{mc1,mc2} which concern the case when
the Hamiltonian is a sum of two-point functions.

For simplicity we consider the weighted setting $(\mu_{0},\phi_{0})$
when $\mu_{0}=dx,$ so that $X=\R^{n}$ and $\phi_{0}$ is thus a
weight function on $\R^{n}$ and we first assume that $\beta>0.$
Then we introduce the following system of Stochastic Differential
Equations (SDE) for $x_{1}(t),....,x_{N}(t)$ viewed as stochastic
processes with values in $\R^{n}:$

\begin{equation}
dx_{i}(t)=-\nabla_{x_{i}}H_{\phi_{0}}^{(N)}(x_{1},...,x_{i},....x_{N})dt+\frac{2}{\beta^{1/2}}dB_{i}(t),\label{eq:stoc gradient flow}
\end{equation}
 where $H_{\phi_{0}}^{(N)}$ is defined by formula \ref{eq:def of hamilt perfurbed by weight}
and $\nabla_{x_{i}}$ denotes the (partial) gradient defined with
respect to the Euclidean metric on $\R^{n}$ and the $B_{i}$s are
$N$ independent standard Brownian motions on $\R^{n}.$ This is thus
a system of Ito diffusions which can be seen as the down-ward stochastic
gradient flow on $X^{N}$ for the $N-$particle Hamiltonian $H_{\phi_{0}}^{(N)}.$
In concrete terms the system is obtained by adding noise to the following
system of ODEs:
\[
\frac{\partial x_{i}}{\partial t}=-\frac{1}{k}\frac{\sum_{\sigma\in S_{N}}p_{\sigma(i)}e^{x_{1}\cdot p_{\sigma(1)}+x_{2}\cdot p_{\sigma(2)}+\cdots+x_{N}\cdot p_{\sigma(N)}}}{\sum_{\sigma\in S_{N}}e^{x_{1}\cdot p_{\sigma(1)}+x_{2}\cdot p_{\sigma(2)}+\cdots+x_{N}\cdot p_{\sigma(N)}}}-\nabla_{x_{i}}\phi_{0}(x_{i})
\]
It seems natural to conjecture that, given the initial condition that
$x_{i}(0)$ be i.i.d variables with law $\rho_{0}dx$ for some fixed,
say smooth and strictly positive probability density $\rho_{0}$ on
$\R^{n},$ the system \ref{eq:stoc gradient flow} converges, when
$N\rightarrow\infty$ (in a sense to be detailed below), to a solution
$\rho_{t}$ of the following deterministic fully non-linear parabolic
system of PDEs with initial data $\rho_{|t=0}=\rho_{0}:$

\begin{equation}
\frac{\partial\rho_{t}}{\partial t}=\frac{1}{\beta}\Delta\rho-\nabla\cdot(\rho\nabla(\phi_{t}-\phi_{0}))\label{eq:system}
\end{equation}
\[
\rho_{t}dx=MA(\phi_{t}),
\]
i.e. $\phi_{t}$ is the unique normalized potential of $\rho_{t}dx$
in the class $\mathcal{P}(\R^{n}),$ which we recall means that $\phi_{t}$
convex and its subgradient image is contained in the given convex
body $P)$ (interestingly, a closely related parabolic system appears
in dynamical meteorology - see \cite{l} and references therein).
More precisely, the converge referred to above should hold in the
following sense: for a fixed time $t,$ the empirical measures $\delta_{N}(x(t)):=\frac{1}{N}\sum_{i=1}^{N}\delta_{x_{i}(t)}$
determined by the the SDEs \ref{eq:stoc gradient flow} converge in
probability to the deterministic measure $\mu_{t}=\rho_{t}dx$. In
fact, if one assumes that the empirical measures $\delta_{N}(x(t))$
converge in probability to \emph{some} deterministic measure $\mu_{t},$
then it can be shown, using Theorem \ref{thm:existence of mean energ}
and the linear Fokker-Planck equations associated to \ref{eq:stoc gradient flow}
(i.e. the corresponding forward Kolmogorov equations), that the density
$\rho_{t}$ of $\mu_{t}$ evolves according to the parabolic PDE \ref{eq:system}.
As is well-known, the general problem of establishing a priori convergence
in probability is essentially equivalent to establishing\emph{ propagation
of chaos} in the sense of \cite{sn} and we leave this problem for
the future. 

In the light of the discussion in the previous section and the connections
to Kähler-Einstein metrics on toric varieties it is also very interesting
to study the case when $\beta<0$ where one would expect that there
exists a global solution to the system \ref{eq:system} in the case
when $\beta\geq-R,$ where $R$ is the invariant of the convex body
$P$ defined by formula \ref{eq:inv R of conv bod}, and that the
convergence statement should hold for $\beta>-R.$ 

Finally, it may be worth pointing out that, inspired by ideas introduced
by Otto (see \cite{vi1,vi2} and references therein), it can be shown
that the parabolic equation \ref{eq:system} is (at least formally)
the down-ward gradient flow for the corresponding free energy functional
$F_{\beta},$ defined with respect to the Wasserstein 2-metric on
the space $\mathcal{M}_{1}(\R^{n}),$ when $\R^{n}$ is equipped with
the Euclidean metric. This observation becomes particularly striking
in the toric setting considered in the previous section, where $F_{\beta}$
may be identified with the Mabuchi K-energy functional. In fact, for
any Kähler manifold $(X,\omega)$ the down-ward gradient flow of the
latter functional, defined with another metric, namely the one defined
by the Mabuchi-Semmes-Donaldson metric on the space of all Kähler
metrics in the Kähler class $[\omega]$ is precisely the \emph{Calabi
flow} which plays a prominent role in Kähler geometry. This also motivates
studying the complex version of the parabolic equation \ref{eq:system},
which is naturally defined on any given Kähler manifold $(X,\omega).$
In fact, the complex version of \ref{eq:system} (for $\beta$ negative)
in the case when $X$ is Riemann sphere is closely related to the
\emph{Keller-Segal system }in $\R^{2},$ which has been extensively
studied in recent years (see for example \cite{b-c-c} and references
therein). It seems also natural to conjecture that the complex version
of the parabolic equation \ref{eq:system} may be obtained as the
large $N-$limit of a systems of SDE's of the form \ref{eq:stoc gradient flow}
obtained by replacing the permanent appearing in the definition of
$H_{\phi_{0}}^{(N)}$ by the corresponding Vandermonde type determinant
(compare \cite{be3,be-3,be,berm2}). But we also leave the study of
this complex story for the future.

\section{\label{sec:Appendix:-the-comparison}Appendix}

\subsection*{A1: Background on optimal transport and its discrete version}

The classical\emph{ assignment problem} (also known as the \emph{bivariate
perfect matching problem} in graph theory) is the problem to, given
an $N$ times $N$ matrix $(c_{ij})$ minimize the functional 
\begin{equation}
\sigma\mapsto\sum_{i=1}^{N}c_{i\sigma(i)}\label{eq:cost of assignement}
\end{equation}
In economical terms we have $N$ workers and $N$ jobs to conduct
and $c_{ij}$ is the cost of assigning work $j$ to a worker $i.$
The problem is to minimize the total cost, if all the every workers
are assigned different jobs, i.e. worker $i$ is assigned the job
$j$ where $j=\sigma(i)$ for some permutation $\sigma\in S_{N}.$ 

The assignment problem relevant to the present paper appears in the
following setting of (discrete) \emph{optimal transport theory}. Consider
two sets $X$ and $P$ in $\R^{n}$ and a given \emph{cost function
}$c(x,p)$ on $X\times P.$ As in the previous sections we denote
by $N$ the number of lattice points in $P,$ i.e. the points in $P_{\Z}:=P\cap\Z^{n}.$
Fix also a configuration $(x_{1},...,x_{N})$ of $N$ points on $X.$
Then we define the \emph{transport cost} from $x_{i}$ to $p_{j}$
as the number $c_{ij}:=c(x_{i},p_{j}).$ Fixing an order $p_{1},...,p_{N}$
of the points in $P_{\Z}$ the problem o\emph{f (discrete) optimal
transport} may then be defined as the corresponding assignment problem,
i.e. as the problem of minimizing \ref{eq:cost of assignement}, for
$c_{ij}:=c(x_{i},p_{j}).$ Concretely, this means that we want to
assign $N$ different points in $P$ to the $N$ given points $x_{i}$
on $X$ in such a way that the corresponding total cost is minimized.
As before we get an asymptotic problem, with $N$ tending to infinity,
by replacing $P_{\Z}$ with $P_{\Z/k}:=P\cap(\Z/k)^{n}$ to get a
sequence of discrete optimal transport problems. When studying asymptotics
it will be convenient to divide the total cost (appearing formula
\ref{eq:cost of assignement}), by the number $N$ of workers to get
the average cost of the work performed. Accordingly, we define the
\emph{(normalized) cost }$C(\sigma)$ by 
\begin{equation}
C(\sigma):=\frac{1}{N}\sum_{i=1}^{N}c(x_{i},p_{\sigma(i)})\label{eq:cost of perm normalized version}
\end{equation}

\subsubsection{\label{sub:Optimal-transport-theory}Optimal transport theory (continuous
version)}

In the classical ``continuous'' setting for optimal transport theory,
as originally introduced by Monge, the given data consist of two probability
measures $\mu$ and $\nu$ on $\R^{n}$ and a \emph{cost function
}$c(x,p)$ on $\R^{n}\times\R^{n}$ (see the monographs \cite{vi1,vi2}
for further background and extensive references). A \emph{transport
map} $T$ is, by definition, a map from $\R^{n}$ to $\R^{n}$ such
that 
\begin{equation}
T_{*}\mu=\nu\label{eq:t pushes mu}
\end{equation}
One defines the \emph{transport cost} of the transport map $T$ as
\[
c(T):=\int_{\R^{n}}c(x,T(x))\mu
\]
and $T$ is said to be an \emph{optimal transport map} if it minimizes
the cost $c(T)$ over all transport maps (i.e. those satisfying the
push-forward formula \ref{eq:t pushes mu}). However, in general such
an optimal transport map may not exist and following Kantorovich one
usually considers a relaxed version of Monge's problem where the transport
map $T$ is replaced with a\emph{ coupling} $\Gamma$ ( between $\mu$
and $\nu)$ i.e. $\Gamma$ is a measure on $\R^{n}\times\R^{n}$ whose
push-forwards to the first and second factor are equal to $\mu$ and
$\nu,$ respectively (such a $\Gamma$ is also called a\emph{ transference
plan}). Its cost is then defined by 
\[
C(\Gamma):=\int_{\R^{n}\times\R^{n}}c(x,p)\Gamma
\]
which is thus the restriction of a \emph{linear} functional to the
space of all couplings (an optimal coupling $\Gamma$ exists under
very general assumptions \cite{vi1,vi2}) Accordingly, fixing $\nu,$
the \emph{optimal total cost }to transport $\mu$ to $\nu$ is defined
by 
\[
C(\mu):=C(\mu,\nu):=\inf_{\Gamma}c(\Gamma),
\]
 where the infimum is taken over all couplings between $\mu$ and
$\nu.$ In particular, any transport map $T$ defines a coupling $\Gamma_{T}:=(I\times T)_{*}\mu$
such that $C(T)=C(\Gamma_{T}).$ 

Assume now that we are given a closed set $X$ in $\R^{n}$ and a
convex body $P$ of unit-mass. We then fix $\nu:=\lambda_{P}$ to
be the Lebesgue measure supported on $P$ (that we will sometimes
also write as $1_{P}dp)$ and consider cost functions $c(x,p)$ of
the form 
\begin{equation}
c_{\phi_{0}}(x,p):=-x\cdot p+\phi_{0}(x)\label{eq:const function for discrete optimal}
\end{equation}
 where $\phi_{0}$ is a given weight function on $X$ (compare section
\ref{sub:Weigted-sets-and}) and denote by $C_{\phi_{0}}(\mu)$ the
correspondning optimal cost functional. It may be decomposed as 
\[
C_{\phi_{0}}(\mu)=C_{0}(\mu)+\int\phi_{0}\mu,
\]
 where $C_{0}$ is the ``unweighted'' cost funcitonal defined with
respect to $c(x,p):=-x\cdot p.$ Since we are only interested in the
dependence of $C_{\phi_{0}}(\mu)$ with respect to $\mu$ we could
also have added any continuous function $\psi_{0}(p)$ to the cost
function $c_{\phi_{0}}(x,p).$ Indeed, this would only shift $C_{\phi_{0}}(\mu)$
by an overall additive constant. A classical case is when $\phi_{0}(t)=\psi_{0}(t)=|t|^{2}/2,$
so that the correspondind cost function is $|x-p|^{2}$ and $C(\mu)$
is hence the Wasserstein $2-$distance between $\mu$ and $\lambda_{P}.$

To see the relation to the discrete setting above we note that for
any given cost function $c(x,p)$ setting $\mu:=\delta(x^{(N)}):=\frac{1}{N}\sum_{i=1}^{N}\delta_{x_{i}}$
for a given configuration $(x_{1},...,x_{N})$ of points on $X$ and
$\nu:=\frac{1}{N}\sum_{i=1}^{N}\delta_{x_{i}}$ clearly gives, 

\begin{equation}
C(\Gamma)=C(\sigma)\label{eq:c cont cost is discrete cost}
\end{equation}
when $\Gamma=\Gamma_{\sigma}=\frac{1}{N}\sum_{i=1}^{N}\delta_{x_{i}}\otimes\delta_{p_{\sigma(i)}}$
for a permutation $\sigma$ in $S_{N}.$ More over, by the Birkhoff-Von
Neumann theorem, when computing the corresponding optimal total cost
$C(\mu,\nu)$ it is enough to minimize over all couplings of the form
$\Gamma_{\sigma},$ i.e. 
\[
C(\mu,\nu)=\inf_{\sigma\in S_{N}}C(\sigma),
\]
In particular, taking $c(x,p)=-|x-p|,$ so that $C(\mu,\nu)$ coincides
with the Wasserstein 1-metric on $\mathcal{M}_{1}(X),$ the previous
equality says that the embedding \ref{eq:def of empricical measure}
is an isometry.

\subsection*{A2: The comparison and domination principles for MA}

In this appendix we will provide proofs of the comparison and domination
principle for the Monge-Ampère operator acting on the function space
$\mathcal{P}(\R^{n})$ associated to a convex body $P$ (following
the notation in section \ref{sub:Setup function spaces in convx}).
These results are without doubt well-known to experts, but for completeness
we have provided proofs that mimic the proofs in the complex setting
(see \cite{begz} and references therein). In fact, the proofs only
use the following basic properties of the real Monge-Ampère operator
$MA:$ 
\begin{itemize}
\item $MA$ is a local operator on $\mathcal{P}(\R^{n}),$ i.e. if $u=v$
on an open set $U$ then $1_{U}MA(u)=1_{U}MA(v).$ 
\item For any $u\in\mathcal{P}_{+}(\R^{n})$ the measure $MA(u)$ is a probability
measure.
\item The space $\mathcal{P}_{+}(\R^{n})$ is closed under the max operation
\end{itemize}
One reason for isolating these properties is that they may be useful
then studying the general space $\mathcal{P}(Y$) of ``ambient potentials''
associated to a Hamiltonian as in section \ref{sub:The-ambient-point},
but we will not go further into this here. We start with a verification
of the second part of the first point above:
\begin{lem}
Let $u$ and $v$ be elements in $\mathcal{P}(\R^{n})$ of maximal
growth, i.e. $u,v$ are in $\mathcal{P}_{+}(\R^{n}),$ which by definition
means that $u-\phi_{P}$ and $v-\phi_{P}$ are bounded. Then 
\[
\int_{\R^{n}}MA(v)=\int_{\R^{n}}MA(u)=1
\]

More generally, if $u$ and $v$ be elements in $\mathcal{P}(\R^{n})$
such that $u\rightarrow\infty$ as $|x|\rightarrow\infty$ and $v\leq u+C$
on $\R^{n}$ for some constant $C,$ then 
\[
\int_{\R^{n}}MA(v)\leq\int_{\R^{n}}MA(u).
\]
\end{lem}
\begin{proof}
The first equality in the lemma can be proved by various means (for
example using the Legendre transform). Here we will instead prove
the more general second inequality only using the fact that $MA$
satisfies a weak form of Stokes theorem. Given a ball $B_{R_{1}}$
of radius $R_{1}$ centered at $0$ there exists a constant $R_{2}>R_{1}$
and constant $\epsilon>0$ and $A>0$ such that $(1-\epsilon)v+A\geq0$
on $B_{R_{1}}$ and $(1-\epsilon)v+A\leq u$ on $B_{R}$ for any $R>R_{2}.$
Hence, setting $\tilde{v}:=\max\{(1-\epsilon)v+A,u\}$ gives $\tilde{v}=v$
on $B_{R_{1}}$and u on $B_{R}$ for $R>R_{1}.$ Hence, $\int_{B_{R_{1}}}MA((1-\epsilon)v)=\int_{B_{R_{1}}}MA(\tilde{v})\leq\int_{B_{R}}MA(\tilde{v}).$
But, by Stokes theorem and a simple approximation argument the latter
integral is equal to $\int_{B_{R}}MA(u)$ and hence $\int_{B_{R_{1}}}MA(v(1-\epsilon))\leq\int_{\R^{n}}MA(u).$
Since, this inequality holds for any $R_{1}>0$ and $\epsilon>0$
this concludes the proof of the lemma. \end{proof}
\begin{prop}
\label{prop:(Comparison-principle)-Let}(Comparison principle) Let
$u$ and $v$ be elements in $\mathcal{P}_{+}(\R^{n})$ (or more generally
elements in $\mathcal{P}(\R^{n})$ of full Monge-Ampère mass). Then
\[
\int_{\{u<v\}}MA(v)\leq\int_{\{u<v\}}MA(u).
\]
\end{prop}
\begin{proof}
Let us first prove the a priori weaker inequality 
\begin{equation}
\int_{\{u<v\}}MA(v)\leq\int_{\{u\leq v\}}MA(u).\label{eq:pf comp pr}
\end{equation}
Since $MA$ is a local operator we have 
\[
1_{\{u<v\}}MA(v)=1_{\{u<v\}}MA(\max(u,v))\leq1_{\{u\leq v\}}MA(\max(u,v)).
\]
 Writing $\{u\leq v\}=\R^{n}-\{u>v\}$ and using locality again hence
gives
\[
1_{\{u<v\}}MA(v)\leq1_{\R^{n}}MA(\max(u,v))-1_{\{u<v\}}MA(u)
\]
Integrating this inequality over $\R^{n}$ and using that, by assumption,
$\int_{\R^{n}}MA(u)=\int_{\R^{n}}MA(v)$ then gives 
\[
\int_{\{u<v\}}MA(v)\leq\int_{\R^{n}}MA(u)-\int_{\{u<v\}}MA(u)
\]
 which hence proves the inequality \ref{eq:pf comp pr}. Finally,
to treat the general case we apply the previous inequality to $u+\delta$
and $v$ giving 
\[
\int_{\{u+\delta<v\}}MA(v)\leq\int_{\{u+\delta\leq v\}}MA(u).
\]
 Finally, letting $\delta\rightarrow0$ and using that the two sequences
of sets $\{u+\delta<v\}$ and $\{u+\delta\leq v\}$ both increase
to $\{u<v\}$ concludes the proof of the proposition.
\end{proof}
Now we can prove the following
\begin{cor}
\label{cor:(Domination-principle)-Let}(Domination principle) Let
$u$ and $v$ be elements in $\mathcal{P}(\R^{n})$ such that $u$
is in $\mathcal{P}_{+}(\R^{n}).$ If $u\geq v$ almost everywhere
with respect to the Monge-Ampère measure $MA(u)$ then $u\geq v$
everywhere on $\R^{n}.$\end{cor}
\begin{proof}
First note that we may as well assume that $v$ is also in $\mathcal{P}_{+}(\R^{n}),$
by replacing $u$ with $\max\{u,v\}.$ In the case when $MA(v)>\delta dx$
for some $\delta>0$ the corollary follows immediately from the previous
proposition. In the general case we simply fix an element $v_{+}$
in $\mathcal{P}(\R^{n})$ such that $MA(v_{+})>\delta dx$ for some
$\delta$ and $v_{+}\leq u$ on $\R^{n}$ (for example, $v_{+}=\log\int_{P}e^{p\cdot x}dp-C$
for $C$ sufficiently large) and apply the previous argument to $u$
and $v_{\epsilon}:=(1-\epsilon)v+\epsilon v_{+}$ for any $\epsilon>0.$
This shows that $v_{\epsilon}\leq u$ on $\R^{n}$ and letting $\epsilon\rightarrow0$
thus concludes the proof. \end{proof}


\begin{thebibliography}{10}
\bibitem{ba}Bakelman, I.J: Convex analysis and non-linear geometric
elliptic equations. Springer-Verlag (1994)

\bibitem{be0}Berman, R.J:\emph{ }Bergman kernels for weighted polynomials
and weighted equilibrium measures of C\textasciicircum{}n\emph{.}
19 pages. Indiana University Mathematics Journal, Volume 58, issue
4, 2009

\bibitem{be}Berman, R.J; Kahler-Einstein metrics emerging from free
fermions and statistical mechanics. J. of High Energy Physics (2011).

\bibitem{ber-bou}Berman, R.J; Boucksom, S: Growth of balls of holomorphic
sections and energy at equilibrium. Invent. Math. 181 (2010), no.
2, 337\textendash{}394.

\bibitem{b-b-w}Berman, R.J:; Boucksom, S; Witt Nyström, D: Fekete
points and equidistribution on complex manifolds. Acta Math. Vol.
207, Issue 1 (2011), 1-27

\bibitem{be3}Berman, R.J: Determinantal point processes and fermions
on complex manifolds: Large deviations and Bosonization. arXiv:0812.4224. 

\bibitem{bbgz}Berman, R.J.: Boucksom, S; Guedj, V; Zeriahi, A: A
variational approach to complex Monge-Ampère equations. Publications
mathématiques de l'IHÉS (to appear). arXiv:0907.4490

\bibitem{berm2}Berman, R.J: A thermodynamical formalism for Monge-Ampere
equations, Moser-Trudinger inequalities and Kahler-Einstein metrics.
arXiv:1011.3976 

\bibitem{ber-ber}Berman, R.J; Berndtsson, B:\emph{ }Real Monge-Ampere
equations and Kahler-Ricci solitons on toric log Fano varieties\emph{.}
arXiv:1209.0996 

\bibitem{be-3}Berman, R.J: A probabilistic approach to Kahler-Einstein
metrics. http://sms.cam.ac.uk/media/1245708. Video and audio from
a workshop on Kahler Geometry in Cambridge 2012, organized by I.Cheltsov
and J.Ross 

\bibitem{be-4}Berman, R.J: A probabilistic approach to Kähler-Einstein
metrics, stability and Coulomb type gases. Article in preparation.

\bibitem{b-c-c}Blanchet, A; Carlen, Eric A.; Carrillo, J. A. Functional
inequalities, thick tails and asymptotics for the critical mass Patlak-Keller-Segel
model. J. Funct. Anal. 262 (2012), no. 5, 2142\textendash{}2230.

\bibitem{b-l}Bloom, T; Levenberg, N. Pluripotential energy and large
deviation. arXiv:1110.6593 

\bibitem{begz}Boucksom, S; Eyssidieux, P; Guedj, V; Zeriahi, A: Monge-Ampère
equations in big cohomology classes. Acta Math. 205 (2010), no. 2,
199\textendash{}262. 

\bibitem{bo}Bovier, A: Statistical Mechanics of Disordered Systems:
A Mathematical Perspective. Cambridge University Press. 

\bibitem{b-e-t}Boucher, C; Ellis, R. S.; Turkington, B: Derivation
of maximum entropy principles in two-dimensional turbulence via large
deviations. (English summary) J. Statist. Phys. 98 (2000), no. 5-6, 

\bibitem{br}Brenier, Y: Polar factorization and monotone rearrangement
of vector-valued functions. Comm. Pure Appl. Math. 44 (1991), no.
4, 375\textendash{}417.

\bibitem{clmp}Caglioti.E; Lions, P-L; Marchioro.C; Pulvirenti.M:
A special class of stationary flows for two-dimensional Euler equations:
a statistical mechanics description. Communications in Mathematical
Physics (1992) Volume 143, Number 3, 501-525

\bibitem{do}Donaldson, S. K. Some numerical results in complex differential
geometry. Pure Appl. Math. Q. 5 (2009), no. 2,

\bibitem{ca0}Caffarelli, L.A: Interior \$W\textasciicircum{}\{2,p\}\$
estimates for solutions of the Monge-Ampère equation. Ann. of Math.
(2) 131 (1990), no. 1, 135\textendash{}150.

\bibitem{ca}Caffarelli, L.A: Some regularity properties of solutions
of Monge Ampère equation. Comm. Pure Appl. Math. 44 (1991), no. 8-9,
965\textendash{}969. 

\bibitem{ca1}Caffarelli, L. A. A localization property of viscosity
solutions to the Monge-Ampère equation and their strict convexity.
Ann. of Math. (2) 131 (1990), no. 1, 129\textendash{}134. 

\bibitem{ca-2}Caffarelli, L.A: The regularity of mappings with a
convex potential. J. Amer. Math. Soc. 5 (1992), no. 1, 99\textendash{}104,

\bibitem{de-ze}Dembo, A; Zeitouni, O: Large deviations techniques
and applications. Jones and Bartlett Publishers, Boston, MA, 1993.
xiv+346 pp.

\bibitem{e-h-t}Ellis, R. S.; Haven, K; Turkington, B: Large deviation
principles and complete equivalence and nonequivalence results for
pure and mixed ensembles. J. Statist. Phys. 101 (2000), no. 5-6, 999\textendash{}1064.

\bibitem{g-mc}Gangbo, W; McCann, R. J.: The geometry of optimal transportation.
Acta Math. 177 (1996), no. 2, 113\textendash{}161. 

\bibitem{gu}Gutierrez, C.E: The Monge-Ampère equation. Progress in
Nonlinear Differential Equations and their Applications, 44. Birkhäuser
Boston, Inc., Boston, MA, 2001. xii+127 pp. ISBN: 0-8176-4177-7 

\bibitem{h-k-p}Hough, J. B.; Krishnapur, M.; Peres, Y.l; Virág, B:
Determinantal processes and independence. Probab. Surv. 3 (2006),
206--229

\bibitem{h-s}Huesmann, M; Sturm, K-T: Optimal Transport from Lebesgue
to Poisson. http://arxiv.org/abs/1012.3845

\bibitem{k}Kiessling M.K.H.: Statistical mechanics of classical particles
with logarithmic interactions, Comm. Pure Appl. Math. 46 (1993), 27-56.

\bibitem{l}Loeper, G: A fully non-linear version of Euler incompressible
equations: the Semi-Geostrophic system. SIAM Journal of Math Analysis
(to appear). 

\bibitem{mc1} McKean, H. P. Jr. A class of Markov processes associated
with nonlinear parabolic equations. Proceedings of the National Academy
of Science 56: 1907-1911, 1966. 

\bibitem{mc2}McKean, H. P. Jr. Propagation of chaos for a class of
nonlinear parabolic equations. Lecture Series in Di\textregistered{}erential
Equations 7: 41-57. Catholic University, Washington, D.C., 1967.

\bibitem{m-s}Messer, J; Spohn, H: Statistical mechanics of the isothermal
Lane-Emden equation. J. Statist. Phys. 29 (1982), no. 3, 561\textendash{}578,

\bibitem{n-o}Negele, J.W; Orland, H: Quantum Many Particle Systems.
Westview Press (1998)

\bibitem{r}Rockafellar, R. T: Convex analysis. Reprint of the 1970
original. Princeton Landmarks in Mathematics. Princeton Paperbacks.
Princeton University Press, Princeton, NJ, 1997. 

\bibitem{sn}Sznitman, A-S: Topics in propagation of chaos. École
d'Été de Probabilités de Saint-Flour XIX\textemdash{}1989, 165\textendash{}251,
Lecture Notes in Math., 1464, Springer, Berlin, 1991

\bibitem{vi1}Villani, C: Topics in optimal transportation, Amer.
Math. Soc., Providence, RI, 2003

\bibitem{vi2}Villani, C: Optimal transport. Old and new. Grundlehren
der Mathematischen Wissenschaften {[}Fundamental Principles of Mathematical
Sciences{]}, 338. Springer-Verlag, Berlin, 2009

\bibitem{w-z}Wang, X; Zhu, X: K\textasciidieresis{}ahler\textendash{}Ricci
solitons on toric manifolds with positive first Chern class, Advances
in Mathematics 188 (2004), 87\textendash{}103.\end{thebibliography}
\end{document}